\documentclass[11pt, reqno]{smfart}
\usepackage{amsmath,amssymb,amsthm,amsfonts,amscd, url, leftidx, xcolor, hyperref, dutchcal}
\usepackage[mathscr]{eucal}
\usepackage{hyperref}
\usepackage[numeric]{amsrefs}
\usepackage{tikz-cd}
\usepackage{todonotes}
\usepackage[toc,page]{appendix}

\usepackage[bbgreekl]{mathbbol}

\DeclareSymbolFontAlphabet{\mathbb}{AMSb}
\DeclareSymbolFontAlphabet{\mathbbl}{bbold}

\input xy
\xyoption{all}

\usepackage{mathptmx} 

\usepackage[scaled=0.90]{helvet} 
\usepackage{courier} 
\normalfont
\usepackage[T1]{fontenc}

\definecolor{cardinal}{rgb}{0.77, 0.12, 0.23}
\newcounter{Todo}
\newcommand{\TODO}[1]{}

\newcommand{\tpoint}[1]{\subsubsection{#1}}

\newcommand{\spoint}{\subsubsection{}}

\newtheorem*{nthm}{Theorem}
\newtheorem*{nlem}{Lemma}
\newtheorem*{nprop}{Proposition}
\newtheorem*{ncor}{Corollary}

\newtheorem*{nclaim}{Claim}
\theoremstyle{definition}
\newtheorem*{define}{Definition}
\newtheorem*{nrem}{Remark}
\theoremstyle{remark}

\numberwithin{equation}{section}

\DeclareMathOperator{\GL}{GL}

\newcommand{\sie}{\widehat{\mf{S}}}

\newcommand{\sJ}{\leftidx^{\ast}J}

\newcommand{\bN}{\mathbf{N}}

\newcommand{\ov}[1]{\overline{#1}}

\newcommand{\mf}[1]{\mathfrak{#1}}
\newcommand{\wh}[1]{\widehat{#1}}
\newcommand{\mc}[1]{\mathcal{{#1}}}
\renewcommand{\mod}{\mathrm{mod}}
\newcommand{\rr}{\rightarrow}

\newcommand{\la}{\langle}
\newcommand{\ra}{\rangle}
\newcommand{\be}[1]{\begin{eqnarray} \label{#1}}
	\newcommand{\ee}{\end{eqnarray}}

\newcommand{\rts}{\mathscr{R}}

\newcommand{\ti}{\mathfrak{h}^+}

\newcommand{\R}{{\mathbb{R}}}
\newcommand{\C}{{\mathbb{C}}}
\newcommand{\A}{{\mathbb{A}}}
\newcommand{\Q}{{\mathbb{Q}}}

\renewcommand{\O}{{\mathcal{O}}}

\newcommand{\bP}{\mathbf{P}}
\newcommand{\bQ}{\mathbf{Q}}

\newcommand{\bG}{\mathbf{G}}
\newcommand{\bA}{\mathbf{A}}

\newcommand{\cen}{\mathrm{cen}}

\newcommand{\vertices}{\mathrm{vert}}
\newcommand{\cham}{\mathrm{cham}}

\newcommand{\alphav}{\check{\alpha}}

\newcommand{\sY}{\mathscr{Y}}

\newcommand{\Ratc}{\mathrm{Rat}}

\newcommand{\hG}{\wh{G}}
\newcommand{\zee}{\mathbb{Z}}
\newcommand{\hA}{\wh{A}}
\newcommand{\hU}{\wh{U}}
\newcommand{\hK}{\wh{K}}
\newcommand{\hgam}{\wh{\Gamma}}
\newcommand{\hp}{P}
\newcommand{\hq}{Q}
\newcommand{\hr}{R}
\newcommand{\hg}{\widehat{G}}

\newcommand{\hw}{\widehat{W}}
\newcommand{\hu}{\widehat{U}}
\newcommand{\hb}{\widehat{B}}

\newcommand{\hP}{P}

\newcommand{\hB}{\widehat{B}}

\newcommand{\dd}{\mathbf{D}}
\DeclareMathOperator{\Lie}{Lie}

\newcommand{\es}{\emptyset}

\DeclareMathOperator{\dep}{depth}
\DeclareMathOperator{\id}{\mathbb{I}}

\newcommand{\End}{\operatorname{End}}
\newcommand{\Aut}{\operatorname{Aut}}

\newcommand\e{\varepsilon}

\renewcommand\L{\Lambda}

\newcommand{\dw}{\dot{w}}

\newcommand\im{\mbox{Im~}}

\renewcommand\Im{\mbox{Im~}}

\newcommand{\cc}{\mathbf{c}}

\newcommand{\lv}{\check{\lambda}}

\newcommand{\av}{\check{a}}

\newcommand{\ratp}{\mathrm{Par}_{\Q}}

\newcommand{\dleq}{\leq_\mathscr{C}}
\newcommand{\dgeq}{\geq_\mathscr{C}}

\newcommand{\hfe}{\widehat{\mathfrak{h}}^e}

\newcommand{\inst}{\mathrm{inst}}
\newcommand{\stdp}{\mathrm{ParStd}_{\mathbb{Q}}}
\newcommand{\cp}{\mathbf{cp}}
\newcommand{\qX}{{_{\Q}}{X}}

\begin{document}
	
	\newcommand{\kin}{K_{\infty}}
	\newcommand{\kino}{K_{\infty}^1}
	\renewcommand{\hw}{\mathbf{v}}
	\newcommand{\conv}{\mathrm{Hull}}
	\newcommand{\gbs}{\mathrm{GBS}}
	\newcommand{\omv}{\check{\omega}}
	\newcommand{\Span}{\mathrm{Span}}
	\newcommand{\hGam}{\widehat{\Gamma}}
	\newcommand{\Xb}{\leftidx_{\Q} \overline{X}}
	
	\newcommand{\hext}{\wh{\mf{h}}^e}
	\newcommand{\tc}{\widehat{\mathscr{T}}}
	\newcommand{\Mor}{\mathrm{Hom}}
	\newcommand{\hf}{\mathfrak{h}}
	\newcommand{\hcf}{\mathfrak{\widehat{h}}}
	\newcommand{\Ac}{\widehat{A}}
	\newcommand{\xg}{X_{\hG^+}}
	\newcommand{\qxg}{\leftidx_{\Q}X_{\hG^+}}
	\newcommand{\bordX}{\leftidx_{\Q}\mathfrak{B}\left(\xg\right)}
	\newcommand{\bordXf}{\leftidx_{\Q}\mathfrak{B}({X_G})}
	
	\newcommand{\hhat}{\widehat{\mathfrak{h}}}
	\newcommand{\hDelta}{\widehat{\Delta}}

	\newcommand{\gammav}{\check{\gamma}}
	\newcommand{\betav}{\check{\beta}}
	\newcommand{\posBor}{\mathscr{B}_+}
	\renewcommand{\ss}{\mathrm{ss}}	
	\newcommand{\sta}[1]{^*{#1}}
	\newcommand{\loopBS}{\bord{\xg}/ \hGam}
	
	\newcommand{\xgss}{\xg^{ss}}
	\newcommand{\cor}{\mathbf{c}}
	\renewcommand{\e}{\mathbf{e}}
	\newcommand{\bhg}{\widehat{\mathbf{G}}}
	
	\newcommand{\bord}[1]{\leftidx_{\Q}\mf{B}\left(#1\right)}
	\renewcommand{\e}{\mathbf{e}}
\newcommand{\hAp}{\hA}

\newcommand{\bsX}{\overline{X^+}}
\renewcommand{\e}{\mathbf{e}}	
\newcommand{\sP}{\leftidx^{\ast}P}
\newcommand{\sQ}{\leftidx^{\ast}Q}

\newcommand{\mq}{M_Q}
\renewcommand{\mp}{M_P}

\newcommand{\cl}{Cl}

\newcommand{\corn}[1]{\cor(#1)}

\newcommand{\Ss}{\mathscr{S}}
\newcommand{\stdmaxp}{\text{StdMax}_{\Q}}
\newcommand{\sspar}{\mathrm{SStdPar}_{\Q}}

\newcommand{\pbord}[1]{\mathscr{X}(#1)}
\newcommand{\cbord}[1]{\mathscr{C}(#1)}

\newcommand{\pr}{\mathrm{pr}}
\newcommand{\vth}{\vartheta}
\newcommand{\vthv}{\check{\vartheta}}
\newcommand{\hfh}{\widehat{\mathfrak{h}}}

\newcommand{\edv}{\check{\lambda}}
\newcommand{\cpsi}{\check{\psi}}
\newcommand{\rhov}{\check{\rho}}

\newcommand{\SL}{\mathrm{SL}}
\newcommand{\SO}{\mathrm{SO}}

\newcommand{\m}{\mathrm{m}}


\DeclareFontFamily{U}{stixbbit}{}
\DeclareFontShape{U}{stixbbit}{m}{it}{<-> stix-mathbbit}{}
\DeclareRobustCommand{\stixdanger}{%
	{\usefont{U}{stixbbit}{m}{it}\symbol{"F6}}%
}

\newcommand{\ueg}{\mmathscr{U}}
\newcommand{\ug}{\mathscr{U}}
\newcommand{\ugo}{\mathscr{U}_o}
\newcommand{\uzgo}{\mathscr{U}_{\zee, o}}
\newcommand{\uzg}{\mathscr{U}_{\zee}}
\newcommand{\ue}{\mathscr{U}}

\newcommand{\wtlat}{\mathscr{L}}
\renewcommand{\k}{\mathsf{K}}

\newcommand{\wts}{\mathrm{wts}}
\newcommand{\bbv}{\mathbbl{v}}

\newcommand{\cbA}{\widehat{\bA}}

\newcommand{\bL}{\mathbf{L}}

\newcommand{\bim}{\mathsf{\Sigma}}
\newcommand{\cbB}{\widehat{\mathbf{B}}}
\newcommand{\cbU}{\widehat{\mathbf{U}}}
\newcommand{\hT}{\widehat{T}}

\newcommand{\maxratp}{\mathrm{MaxPar}_{\Q}}
\newcommand{\build}{\mathcal{T}_{\Q}^+}

\newcommand{\adb}{\Omega}

	\title{Borel--Serre type bordifications and Canonical Pairs for Loop groups}

	\author{Manish M. Patnaik}
	\email{patnaik@ualberta.ca}

	\author{Punya Plaban Satpathy}
	\email{psatpathy81@gmail.com}

\date{\today}
\begin{abstract}
		To the symmetric space of the (positive half) of a real loop group, we attach a Borel--Serre type bordification and equip it with a Hausdorff topology. The attached boundary, indexed by certain rational parabolics of the loop group, is shown to be homotopic to an affine, rational Tits building. A loop analogue of an arithmetic group is also shown to act continuously on the bordification and its quotient by this action is studied using the reduction theory of H. Garland. While the quotient is no longer compact (as in the Borel--Serre construction from finite-dimensions) we relate the non-compactness to the center of the loop group. We also introduce a notion of semi-stability for loop groups, following works of Harder--Narasimhan, Behrend, and most recently Chaudouard, and use this to describe a partition of our loop symmetric space. This partition is then related to the rational bordification and its quotient. 
	
\end{abstract}

	\maketitle 
	\setcounter{tocdepth}{2}
	\tableofcontents
	\section{Introduction} Let  $G$ be the group of real points of a (split) simple Chevalley group, $\hg$ the Kac--Moody central extension of the $\R((t))$ points of $G$, and $\hg^e$ the extended loop group \textit{i.e.} the semi-direct product of $\hg \rtimes \R^*$, where the second factor acts by  loop rotations.   Let $\hK \subset \hG$ and $\hGam \subset \hG$ be the loop group versions  of the maximal compact subgroup and the subgroup of integral points in $\hG$ respectively (see \S \ref{sec:loop-groups}). As in the work of H. Garland \cite{gar:ihes},  to formulate a reduction theory for $\hGam$, we work within a certain subset $\hG^+ \subset \hG^e$ where the loop rotation acts in a `contractive' \footnote{This means that the action of loop rotation on a one parameter subgroup corresponding to large positive root is by rescaling through a small, positive number. The notion of positive depends on our choice of Borel subgroup, and we take it to contain positive powers of $t.$} manner and consider \be{} \label{xg:quot-intro} \begin{array}{lcr}  \xg:= \hK \setminus \hG^+ & \mbox{ and } &\xg / \hGam \end{array} \ee as loop group analogues of a symmetric and locally symmetric space. The former can be equipped with the structure of a (contractible) topological space on which $\hgam$ acts continuously, and in this paper, we construct a \textit{bordification} $\bord{\xg} \supset \xg$ with the following properties:  
	
	 \begin{itemize}
	\item[a)] $\bord{\xg}=\xg \cup \, \bigsqcup_{P \in \ratp} \e(P)$ where $P$ ranges over certain rational parabolic subgroups of $G$ and $\e(P)$ is constructed from lower rank (finite-dimensional) symmetric spaces and the (infinite-dimensional) pro-unipotent radical of $P$;
	\item[b)] the boundary in $ \bord{\xg}$ is homotopic to the rational (affine) Tits building of  $\hG_{\Q}$; and,
	\item[c)]  the (continuous) $\hGam$ action on $\xg$ extends to a continuous action of $\hG_{\Q}$ on $\bord{\xg}.$
	\end{itemize} 
	
		\noindent If $\hG$ is replaced by a finite-dimensional reductive  group, the analogous construction is known as the Borel--Serre bordification, see  \cite{bs-announcement, borel-serre}. Moreover, it is show in \textit{op. cit.} that the quotient of this bordification by an arithmetic subgroup is a compact Hausdorff space, namely the Borel--Serre or Siegel\footnote{See \cite{siegel} and \cite[footnote (2), p.1158]{bs-announcement}} compactification. In our  setting, $\bord{\xg}/\hGam$ is no longer compact due to the presence of the central extension of the loop group. With future applications to the cohomology of $\hGam$ in mind, we do \textit{not} want to further compactify this central direction and so what we obtain might be called a Borel--Serre partial compactification.

		In the main body, we consider and relate two approaches to $\bordX$ and its quotient.  Our first approach, the `attachment method' of Borel--Ji (see \cite{ji-macpherson}, \cite[Ch. 10]{borel-ji}), adds in certain boundary pieces to $\xg$ `by hand' and then specifies a topology on the resulting bordification $\bordX$ using a family of convergent sequences.  The properties of the quotient are probed using  the reduction theory of H. Garland.   The second approach, the `partition method,' we learnt from the (unpublished) manuscript of G. Harder  \cite{harder-buch}, \textit{cf.} Raghunathan \cite{ragh}, Leuzinger \cites{leu, leu:poly}, Saper \cite{saper} for closely related ideas. The basis of this method is a partition of $\xg$ obtained using the concept of semi-stability, a notion which we introduce here following Chaudouard \cite{chau} and whose origin is the work of Harder and Narasimhan  \cite{harder-narasimhan} on destabilizing parabolics from the moduli theory of vector bundles on curves  (see Remark \ref{intro-remark-references-HN}.) Our partition has pieces indexed by parabolics, where the part corresponding to $\hG$ is called the semi-stable locus and those corresponding to proper parabolics $P$ are referred to as $P$-semi-stable ends. In finite dimensions, the semi-stable locus is invariant under the action of the arithmetic group and has a  compact quotient under this action, so the non-compactness is due to the $P$--semi-stable ends. In the loop setting, though the semi-stable locus does not have compact quotient (due again to the center), we do obtain a natural fibration describing the $P$--semi-stable ends in terms of lower-rank symmetric spaces and certain non-compact toral directions. Our construction of $\bord{\xg} / \hGam$ amounts to a  completion of these toral ends.

		 We plan to use the results here for further investigations into the theory of automorphic forms on loop groups and also to begin the study of cohomology of $\hGam$. The length of this paper notwithstanding, we have aimed to work in the simplest setting that could also reveal the novel features of loop symmetric spaces.  So we are always dealing with $\Q$, but we do expect our methods, when combined with existing techniques \textit{cf.} \cite{gar:lg2}, will yield results for adelic groups over general global fields. In the function field setting, it would be interesting to compare our results to the compactifications in the theory of bundles on an algebraic surface (note that our non-compact central direction is related to the second Chern class in that theory, an observation due to Kapranov \cites{kapranov, pat:thesis}). We also restrict ourselves to just working with the analogue of the `full modular group' $\hGam$ (which is essentially $\hG_{\zee}$), but our techniques should carry over to other `arithmetic' sugroups of $\hG$ once this theory is developed.

In the remainder of this introduction, we discuss in some detail the finite-dimensional picture that motivated us in \S \ref{intro:sub-partition-method}, focusing especially on the case of $G= \SL_2$ where things be visualized easily. Though some aspects of our treatment for finite-dimensional groups may be novel, our main new results on loop groups are described in \S \ref{sub:intro-loop-groups}.

\subsection{The finite dimensional case} 
\label{intro:sub-partition-method}

	\tpoint{Notation} \label{notation:conjugation} For a group $G$ and subset $A \subset G$, we write \be{} \begin{array}{lccr} \label{conj} \leftidx^xA := xAx^{-1} & \mbox{ and } & A^x = x^{-1}Ax &  \mbox{ for } x \in G. \end{array} \ee If $A= \{ y \}$ is a singleton set, we write $\leftidx^x y$ and $y^x$ for the corresponding conjugations. Also, for a finite set $S$, we write $|S|$ to denote the number of elements in the set S.

	\tpoint{Basic setup } \label{subsub:intro-basic-setup} Let $G:= \bG(\R)$ be a Chevalley group with maximal compact subgroup $K$ and let $\Gamma$ be the stabilizer of the natural lattice in the representation defining $G$. As the real and rational rank $\ell > 1 $, the quotient of $X_G:= K \setminus G$ by $\Gamma$ will be non-compact.  Fix $H$ a maximal torus of $G$ with $A$ its connected component of the identity. Then $A \cong \R_{>0}^{\ell}$ and we write $\mf{h}= \Lie(A)$. Choose $B$ a Borel subgroup whose unipotent radical will be denoted as $U$. 
	
	  The Iwasawa decomposition states $G = K \times A \times U$; for $g \in G,$ write (uniquely) $g = k_g a_g u_g$ with respect to this decomposition. So $X_G = A \times U$, and use the natural topologies on $A$ (a torus) and $U$ (a unipotent group) to  give $X_G$ a natural (contractible) topology with respect to which we have a continuous map  $H_B: X_G \rr \mf{h}$ that sends $x \mapsto \log (a_g)$.  Define  $\stdp$ and $\stdp'$ as the set of parabolics (resp. proper parabolics) containing $B$. We also write $\ratp$ and $\ratp'$ for the set of rational (resp. proper rational) parabolics, each one of which is conjugate under $G_{\Q}$ to a member of $\stdp$ (resp. $\stdp'$).  For $P \in \ratp'$ we have $P$-horospherical decompositions (see \cite[III.1]{borel-ji}) \be{} \begin{array}{lcr} P = M_P \times A_P \times U_P & \mbox{ and } & X_G = X_{M_P} \times A_P \times U_P \end{array} \ee respectively, where $A_P$ is the connected component of the center of the Levi, $M_P$ is the semi-simple part of the Levi,  $U_P$ is unipotent radical of $P$, and $X_{M_P}$ is a lower rank symmetric space attached to $M_P$. Denote by $H_P: X_G  \rr  \mf{h}_P:= \Lie(A_P)$ the (continuous) logarithm of the $A_P$ coordinate.  Let $\Delta_P$ denote the simple roots of $P$ respectively. Then for $t > 0$  we define 
	  \be{} \label{def:AP-intro} A_{P, t} := \{ h \in A_P \mid \la \alpha, H_P(h) \ra < \log \, t \mbox{ for } \alpha \in \Delta_P \} \cong (0, t)^d  \ee where $d = |\Delta_P|$.
	  We also define, given $\Omega \subset \R$ the subset $U_{P, \Omega}$ obtained by restricting the root-group coordinates (in some fixed order) to lie in $\Omega$.  The first fundamental theorem of reduction theory (see \cite[Thm. 13.1 and 16.7]{borel:red}) states that if $t$ is sufficiently small ($t < 2 / \sqrt{3}$ works) and $\Omega$ sufficiently large ($\Omega \supset \{ x \in \R\mid |x| \leq 1/2 \}$ works), then every $x \in X_G$ has a representative, after right translation by $\Gamma$, in the subset $\mf{S}_{t, \Omega} = A_{B, t} \times U_{B, \Omega}.$ The theorem of parabolic transformations (see \cite[Thm. 12.6]{borel:red}) further studies  the intersections $\mf{S}_{t, \Omega} \gamma \cap \mf{S}_{t, \Omega}$ for $\gamma \in \Gamma$, a typical result being that for $t$ sufficiently small: this set can be non-empty only if $\gamma \in \Gamma \cap P$ for $P \in \stdp'$.
	  
	  \tpoint{Example: $G= \SL_2$}  \label{intro-ex-sl2-part1} Suppose now that $G=\SL_2(\R)$, $K= \SO(2)$, and $\Gamma = SL_2(\zee)$. There is a bijection $X_G  \stackrel{1:1}{\longrightarrow} \mathbb{H}$ sending $g \mapsto g^{-1}.i$, and where the action of $G$ on $\mathbb{H}$ is via fractional linear transformations.  Let $B$ and $U$  denote the subgroup of upper-triangular  and unipotent matrices in $G$ respectively. Defining $A = \{ \big(\begin{smallmatrix}
	  	a & 0\\
	  	0 & a^{-1}
	  \end{smallmatrix}\big) | a \in \R^+\}$, the Iwasawa or $B$-horospherical coordinates $X_G =  A \times U $ can be expressed as follows: $z= x+iy \in \mathbb{H}$ is identified with $(a(z), n(z)) \in A \times U$ where  \be{}\begin{array}{lcr}  a(z) = \begin{pmatrix} 
	  		y^{-1/2} & 0 \\
	  		0 & y^{1/2} 
	  	\end{pmatrix} \in A & \mbox{ and } & n(z) = \begin{pmatrix} 
	  		1 & -x \\
	  		0 & 1 
	  	\end{pmatrix} \in U  \end{array}.\ee Hence if $H:= \begin{pmatrix} 1 & 0 \\ 0 & -1 \end{pmatrix}$ is a chosen generator of $\mf{h}=\Lie(A)$,  the map $H_B: X_G \rr \mf{h}$ becomes \be{}  \begin{array}{lcr} H_B: \mathbb{H} \rr \mf{h}_B, & \mbox{ where } &  z \mapsto \log(a(z)) = -\frac{1}{2} \log( \im(z) ) H. \end{array} \ee

  	\noindent In this case, $\stdp= \{ B, G\}$ and $\ratp$ is in bijection with \be{} \Gamma / \Gamma \cap B = \SL_2(\Q)/ B(\Q) \cong \mathbb{P}^1(\mathbb{Q}). \ee
  	  To $\infty=[1:0]$ corresponds to the coset of the identity, denoted as $e$, and the corresponding parabolic is just $B=B_e$. 
  	  To an element  $m \in \mathbb{Q}:= \mathbb{P}^1(\mathbb{Q}) \setminus \{ \infty \},$ say $m=[r:s]$ with $r, s \in \zee, s \neq 0,$ and $(r, s)=1$, we attach a matrix $\gamma_{m}:= \begin{pmatrix} a & b \\ s & -r \end{pmatrix}$ for some \footnote{The $a, b$ are not unique, but the class of $\gamma_{m}$ in $\Gamma/ \Gamma \cap B$ does not depend on the choice} $a, b \in \zee$ and rational parabolic $B_m:= B^{\gamma_m}$. 
  	  Writing $A_m$ and $U_m$ for corresponding torus and unipotent radical in $B_m$, if $z \in X_G= \mathbb{H}$ has coordinates $(a_m(z), n_m(z))$ with respect to $X_G = A_m \times U_m$, then under a standard identification $A_m \cong \R_{>0}$ and $U_m \cong \R$, we can regard $a_m(z)$ as the inverse of the radius of the unique circle in $\mathbb{H}$ tangent to $m$ and passing through $z$, while $n_m(z)$ can be regarded as the (other) point of intersection with the $x$-axis with unique geodesic passing through $m$ and $z.$ The same geometric picture holds for $\infty$ as well, as the reader can easily check.

	   \tpoint{Semi-stability} Returning to the general setting, let $\rho$ denote the sum of all the positive roots of $G$, $\rho(P)$ the corresponding sum of the all positive roots of $M_P,$ and $\rho_P= \rho - \rho(P).$ For each $x \in X$, we define the \textit{degree of instability} \be{} \label{intro:deg-inst} \deg_{\inst}(x):= \min\limits_{P \in \stdp', \delta \in G_{\Q}/ P_{\Q}} \la \rho_P, H_P(x \delta) \ra, \ee where the minimum is taken first over proper, standard parabolics $P$, and then, for a fixed $P,$ over points in the flag variety $G_{\Q}/P_{\Q}$. 

	   \begin{nthm}(Theorem of Canonical Pairs)  \label{thm:intro-canonical-pairs} \begin{enumerate}
	   		\item For every $x \in X_G$, the minimum defining $\deg_{\inst}(x)$ can be achieved. 
	   		\item Either $\deg_{\inst}(x) \geq 0$ in which case we say that $x$ is semi-stable, or $\deg_{\inst}(x)<0$ and among all the $(P, \delta)$ realizing the minimum, there is a unique pair which maximizes $P$ with respect to inclusion. This is called the canonical pair of $x$ and written $\cp(x) = (P, \delta)$
	   		\item A pair $(P, \delta)$ with $P \in \stdp'$ and $\delta \in G_{\Q}/P_{\Q}$ is a canonical pair for $x$ if and only if \begin{enumerate}
	   			\item The projection of $x \delta$ onto $X_{M_P}$ is semi-stable for $M_P$
	   			\item The projection of $x \delta$ onto $A_P$ lies in $A_{P, 1}$, the set defined in \eqref{def:AP-intro}.
	   		\end{enumerate}
	   		
	   	\end{enumerate}
   	
   	 \end{nthm}
	\begin{nrem}	\label{intro-remark-references-HN}
	In this formulation, the results are essentially due to Chaudouard \cite[\S2]{chau} (written adelically), and motivated by results in the function field setting, due to Harder--Narasimhan \cite{harder-narasimhan} for vector bundles and to Behrend \cite{beh} for $G$-bundles (in this last context, a new proof was also found by Drinfeld--Gaitsgory \cite{drin:gai} and Schieder \cite{schieder}). For $GL_n(\R)$ and in the language of lattices, this result is contained in the work of Stuhler and Grayson \cite{stuhler, grayson}; details of the comparison of Chaudouard's definition and the lattice picture are being worked out, see T. Nguyen\footnote{MSc thesis, 2025, University of Alberta}. \end{nrem}

	Write $X^{ss}_G$ for the semi-stable elements and $X_G(P, \delta):= \{ x \in X_G \mid \cp(x) = (P, \delta) \}$ for each $P \in \stdp'$ and $\delta \in G_{\Q}/P_{\Q}.$ With this notation, we obtain a partition \be{} \begin{array}{lcr} X_G=X_G^{ss} \sqcup \, \bigsqcup\limits_{P \in \stdp} \, \qX_G(P) & \mbox{ where } & \qX_G(P) = \sqcup_{\delta \in G_{\Q}/P_{\Q}} \, X_G(P, \delta) \end{array}. \ee
	
	\noindent Finally, we remark that part (3) of the theorem gives us a fibration \be{} \label{intro-fibration} f_{(P, \delta)}: X(P, \delta) \rr X_{M_P}^{ss} \times A_{P, 1}, \ee with fibers homeomorphic to $U_P$.  If $\delta=e$, then we have an even simpler product structure: \be{} \label{XP:e} X_G(P,e) = X_{M_P}^{ss} \times A_{P,1} \times U_P. \ee 

	\tpoint{Example: $\SL_2$-continued} \label{intro-ex-sl2-part2} Keeping the notation of \S \ref{intro-ex-sl2-part1}, we also denote by $\alpha: \mf{h} \rr \C$ the unique linear map such that $\alpha(H) = 2$ and set $\rho:=  \alpha/2.$ One may check  \be{} \begin{array}{lcr} \deg_{\inst}(x):= \min\limits_{\gamma \in \Gamma/\Gamma \cap B}  \la \rho, H_B(x \gamma)  \ra  = \min\limits_{\gamma \in \Gamma/\Gamma \cap B}  \Im ( \gamma^{-1} z ) & \mbox{ for } & x^{-1}.i=z\in \mathbb{H} \end{array}. \ee 	Hence we have that 
	\be{} X_G^{ss} &=& \{z  \in \mathbb{H} | \, \im (\gamma^{-1} z) \leq 1  \text{ for all }  \gamma \in \Gamma/ \Gamma \cap B \}   \\  &=& \{z \in \mathbb{H} | \, \im(z) \leq 1 \} \setminus \bigcup\limits_{\gamma \in \Gamma \cap B \setminus \Gamma ,\gamma \neq e} \{z \in \mathbb{H}| \,  \text{Im}(\gamma z) > 1 \}.  \ee
 In terms of the identification of $\Gamma / \Gamma \cap B$ with $\mathbb{P}^1(\Q)$ used in \S \ref{intro-ex-sl2-part1}, we may describe this set a bit more explicitly as follows. Let $m \in \mathbb{P}^1(\Q) \setminus \{ \infty \}$, with $m= [r:s]$, $s \neq 0$ and with $(r,s)=1$, for $\gamma_m$ the element defined in \S \ref{intro-ex-sl2-part1}, we find that \be{} \begin{array}{lcr} \im(\gamma_{m} z ) > \epsilon^{-1} & \mbox{if and only if } & \big(\frac{\epsilon}{2s^2}\big)^2> \big(x -\frac{r}{s}\big)^2+  \big(y - \frac{\epsilon}{2s^2})^2 \end{array}. \ee Hence, if we define $S(m)$ to be the circle in $\mathbb{H}$ tangent to the real axis at the point $m$ with radius $\frac{1}{2s^2},$ and $D(m)$ for the (open) disc this circle bounds, and also set \be{}
\begin{array}{lcr} S(\infty,  1):= \{ z \in \mathbb{H} \mid \im(z) =  1 \} & \mbox{ and } & D(\infty):= \{ z \in \mathbb{H} \mid \im(z)>  1 \} \end{array} \ee  for the corresponding objects at $\infty$, we have \be{} \begin{array}{lcr} X^{ss}_G = X_G \setminus \, \cup_{m \in \mathbb{P}^1_{\Q}} D(m), & X_G(B) = \cup_{m \in \mathbb{P}^1_{\Q}} D(m), & \mbox{ and }X_G(B, m) = D(m). \end{array} \ee

\newcommand{\pone}{\mathbb{P}^1(\Q)}

	\tpoint{Completion of $A_P$}  Let $P \in \ratp$ and $d= |\Delta_P|$. Then we have a homeomorphism $A_P \cong (\R_{>0})^{d}$ and a natural completion is obtained by taking the closure in $\R^{d}$, i.e. considering $(\R_{\geq 0})^{d}.$ We can express this in more Lie theoretic terms as follows. If $P \subset Q$ are two (standard) parabolics, then $P$ defines a parabolic $\sta{P} \subset M_Q$ which has a Langlands decomposition  $^*P = M_{^*P} \times  A_{^*P} \times  U_{^*P}$ \be{} \begin{array}{lcr} M_P =M_{^*P} , & U_P= U_Q \, U(Q)_{^*P} & \mbox{ and } A_P = A_Q A(Q)_{^*P}, \end{array} \ee where $U(Q)$ and $A(Q)$ denote the unipotent radical and connected center of  $\sta{B} \subset M_Q.$ Letting \be{} \label{intro:finite-bord-Ap}\cor(A_P) :=  A_P \ \  \cup \ \ \sqcup_{Q\supset P} \, \, A(Q)_{^*P}, \ee one can introduce a topology on this set, described naturally in Lie theoretic terms, which renders $\cor(A_P) $ homeomorphic to $(\R_{\geq 0})^{d}$ extending the homeomorphism $A_P \cong (\R_{> 0})^{d}$. We remark that $\cl_{\corn{A_P}} \, A_{P, 1},$  the closure in $\corn{A_P}$ of the set $A_{P, 1}$ from \eqref{def:AP-intro}  is now compact and homeomorphic to $[0, 1]^d$. In addition to this closure, we shall also consider the interior of this closure, namely  \be{} \corn{A_{P, 1}}  = A_{P, 1} \, \sqcup \, \bigsqcup_{Q \supset P} A(Q)_{\sta P, 1} \cong [0, 1)^d. \ee

		
	\tpoint{Bordification}  For each $Q \in \ratp'$, introduce the \textit{boundary component}  $\e(Q) = X_{M_Q} \times U_Q$ as well as the \textit{bordification} of $X_G$ as \be{} \bord{X_G} = X_G \sqcup \, \, \bigsqcup_{Q \in \ratp'} \e(Q). \ee We refer to the \textit{corner attached to P} as the subset   \be{} \cor(P) := X \cup \, \sqcup_{Q\supset P} \, \e(Q). \ee

	\begin{nthm} \begin{enumerate}
			\item \cite[Theorem 7.8]{borel-serre} There is a natural way to introduce a Hausdorff topology on $\bord{X_G}$ so that the inclusion $X_G \hookrightarrow \bord{X_G}$ extends to an open embedding from $X_{M_P} \times \corn{A_P} \times U_P$ into $\cor(P)$, where the former set has the product topology. 
			\item \cite[Thm. 8.5.2]{borel-serre} The boundary of $\bord{X_G},$ namely $\bord{X_G} \setminus X_G = \sqcup_{P\in \ratp'} \e(P)$, has the same homotopy type as the rational Tits building attached to $G_{\Q}$.

		\end{enumerate}
\end{nthm}

	\tpoint{Bordifications and semi-stability} \label{subsub:intro-bord} 
	For $P \in \stdp'$, let us write $\cl_{\corn{P}} \,  X_G(P, e)$ for the closure of $X(P, e)$ in the corner $\corn{P}$ and $\pbord{P}$ for the interior of this closure. By translation through $G_{\Q}$ we can also define $\pbord{S}$ for any $S \in \ratp'$. For $P \in \stdp', \delta \in G_{\Q}/P_{\Q}$ and $S = P^{\delta}$, let us write $X_G(S):= X_G(P, \delta)$ and now denote the fibration in \eqref{intro-fibration} as $f_S: X_G(S) \rr X^{ss}_{M_P} \times A_{P, 1}.$ If we write $i$ for the natural inclusions $X(S) \hookrightarrow \pbord{S}$ and also $X_{M_P} \times A_P \hookrightarrow X_{M_P} \times \corn{A_P},$ then one can extend $f_S$ to a map $\mf{f}_S$ making the following diagram commute 
\be{} 
\begin{tikzcd} 
	 X_G(S) \arrow[d, "f_S"] \arrow[r] & \pbord{S} \arrow[d, "\mf{f}_S"] \\  
	 X^{ss}_{M_P} \times A_{P, 1} \arrow[r, "i"] & X^{ss}_{M_P} \times \corn{A_{P, 1}} 
\end{tikzcd} \ee
	
\noindent If  $S=P \in \stdp'$ then the fibration degenerates to a product, and in fact, \be{} \pbord{P} &=& X_{M_P}^{ss} \times \corn{A_{P, (0, 1)}} \times U_P \\
			&=& X_G(P, e)  \, \sqcup \, \bigsqcup_{Q\supset P} X_{M_Q}(\sta{P}, e) \times U_Q. \ee 
	\noindent In other words $\pbord{P}$ is obtained from $X(P,e)$ by completing in the torus direction $A_P$ using the fibration \eqref{intro-fibration};  a similar description holds for a general $S$. Since we also have \be{} \label{bord:pbord-intro} \bord{X_G} = X_G^{ss} \, \sqcup \, \bigsqcup_{S \in \ratp'} \, \pbord{S}, \ee generalizing the partition of $X_G$, this gives us an alternative way of thinking about $\bord{X_G}$.

 	\tpoint{Compactness and separability} \label{intro:compactness} Turning to the quotient $\bordXf/ \Gamma,$ one has

 	\begin{nthm} \label{thm:intro-compactness} The natural action of $G_{\Q}$ on $X_G$ extends to a continuous, proper action on $\bord{X_G}$. The quotient $\bord{X_G}/\Gamma$ is a compact, Hausdorff space. \end{nthm}

 	The Hausdorff property follows from the properness of the action. The properness in turn can be reduced to Siegel finiteness \cite[]{borel:red}:  that for $t, \Omega$ as above, $\mf{S}_{t, \Omega} \gamma \cap \mf{S}_{t, \Omega}$ is non-empty for only finitely many $\gamma \in \Gamma.$  As for the compactness, let us give two (related) arguments.
 	
 	The argument of \cite{borel-serre} proceeds as follows: pick a Siegel set $\mf{S}_{t, \Omega}$ as in \S\ref{subsub:intro-basic-setup} and let $\Ss_{t, \Omega} = \corn{A_{B, t}} \times U_{\Omega}$ be the interior of its closure in $\corn{B}$. Notice that $\corn{A_{B, t}}$ is relatively compact in $\corn{A_B}$.  Now, one can show that $\cup_{\gamma \in G_{\Q}} \Ss_{t, \Omega} \gamma = \bord{\xg}$, so that the image of $\cl_{\corn{P}}{ \Ss_{t, \Omega}}$ is compact and maps onto $\bord{X_G}/\Gamma$.  Alternatively, we first note that \eqref{bord:pbord-intro} and Theorem \ref{thm:intro-canonical-pairs} imply  \be{} \bord{X_G}/ \Gamma = X_G^{ss}/\Gamma \sqcup \bigsqcup_{P \in \stdp'} \pbord{P}/ \Gamma_{P}, \ee and we then try to argue each piece of the above union is relatively compact. As for the semi-stable part,  let  $\mathfrak{S}_{t,\Omega} ^{ss} = X_G^{ss} \cap \mathfrak{S}_{t,\Omega}$ denote the semi-stable points in the Siegel set.  If $x \in \mathfrak{S}_{t,\Omega}^{ss} $, then from the Siegel set condition $ \la \alpha,H_B(x) \ra < \log t $ for all $\alpha \in \Delta_B$ and, from the semi-stability assumption, $\la \lambda ,H_B(x) \ra \geq 0$ for all $ \lambda \in \hat{\Delta}_B,$ the set of fundamental weights of $G$. Hence
 	\be{} \begin{array}{lcr} \mathfrak{S}_{t,\Omega} ^{ss} \subset A^{1}_{B,t} \times U_{B, \Omega} & \mbox{where} & 
 		A^{1}_{B,t} = \{a \in A_{B,t} \mid a^{\lambda} \geq 1, a^{\alpha} \leq  t \mbox{ for all } \lambda \in \hat{\Delta}_B, \alpha \in \Delta_B \}. \end{array} \ee Using essentially the fact that $\lambda \in \hat{\Delta}_B$ is a \textit{non-negative} linear combination of the elements from $\Delta_B$ (see \cite[Lemma 1.1]{ragh}), we can show that $A^{1}_{B,t}$ and hence $A^{1}_{B,t} \times U_{B, \Omega}$ is a relatively compact subset of $X_G$ and hence $X^{ss}_G/\Gamma$ will be compact. 
 		
 		On the other hand, for $P \in \stdp'$  we have seen that $X(P, e) = X^{ss}_{M_P} \times A_{P, 1} \times U_P$, so that  $X(P, e) /\Gamma_P = X^{ss}_{M_P} / \Gamma_{M_P} \times A_{P, 1} \times U_P/\Gamma_{U_P}$, and $\pbord{P}/ \Gamma_P$ is obtained has a similar description obtained by replacing $A_{P, 1}$ with the compact piece $\corn{A_{P, 1}}$.  Now $X^{ss}_{M_P}/\Gamma_{M_P}$ is compact by the same argument as in the last paragraph and it is easy to see that $U_P/ \Gamma_{U_P}$ is also compact. The desired relative compactness follows.

	\tpoint{Example: $\bord{X_{\SL_2}}$} \label{intro-ex-sl2-part3} We continue with Example \ref{intro-ex-sl2-part2}. From the description of $X_G^{ss}$ given there, it is clear that $X_G^{ss} \cap \mf{S}_{t, \Omega}$ is relatively compact. As for $\pbord{P_m}$ with $m \in \mathbb{P}^1_{\Q}$, first note that $A_{B, 1}$ is homeomorphic to $(0, 1)$ and $\corn{A_{B, 1}} \cong [0, 1)$. Hence, the set $\pbord{P_m}$ is obtained from $X(P_m):=X(B, \gamma_m)$ by adding in a copy of $\R \cong U_{P_m}$ This can be visualized as either lying on the $x$ axis (if $m \neq \infty$) or lying parallel to it at an infinite height (if $m=\infty$). If $m \neq \infty$, a sequence of points $x_k \in X(P_m)$ achieves a limit along this newly added line, if when we write $P_m$-coordinates $x_k=(a_k, u_k)$ the $a_k \rr 0$ and $u_k \rr u_{\infty} \in \R$.  As for the quotient of $\sqcup_{m \in \mathbb{P}^1(\Q)} \pbord{P_m}$ by $\Gamma$, it may be easiest to consider the case $m=\infty$, in which case $\Gamma_{P_m}=\Gamma_B$ just acts by horizontal translations (even on the new line added) and so  $\pbord{B}/ \Gamma_B$ is a cylinder. 

	\subsection{The case of loop groups} \label{sub:intro-loop-groups}

	\tpoint{Basic setup} We assume the connected torus $\hA^e$ in $\hG^e$ has rank $\ell+2$, where $\ell$ is the rank of the underlying finite-dimensional group $G$. We write $I= \{ 1, \ldots, \ell+1\}.$ Let $\hfh^e:= \Lie(\hA^e)$, and note that $\hfh^e = \R \dd \oplus \mf{h} \oplus \R \cc$ where $\mf{h}$ is the Cartan of  $\mf{g}=\Lie(G)$, $\cc$ is a central element, and $\dd$ is the degree derivation. Writing $\av_{\ell+1}$ for the new simple coroot, the subalgebra $\hfh^e$ has two natural bases: in terms of the simple coroots $\{ \av_, \ldots, \av_{\ell+1}, \dd\}$, which allows one to relate it to $\mf{h}$; and in terms of fundamental coweights $\{ \edv_1, \ldots, \edv_{\ell+1}, \cc \}$ which allows one to relate it to parabolic subgroups.

	As mentioned earlier, we need to work within a subset $\hG^+ \subset \hG^e$. This amount to working within $\hfh^+ \subset \hfh^e$, the Tits cone, which is concretely described as  \be{} \hfh^+ &:=& \cup_{r >0}  \{ H - r \dd \mid H \in \Span_{\R} \{ \av_i \}_{i \in I} \} \\ 
&=& \{ \sum_{i \in I} p_i \edv_i + m \cc \mid \sum_{i \in I} p_i d_i < 0  \},\ee for $d_i$ certain positive constants depending on the root system. At the group level, we often write $\eta(s):= \exp(-r \dd)$ and denote by $A_{\cen}$ the subgroup corresponding to $\cc$ and $A= \exp(\mf{h}).$ Then \be{} \label{intro:A-std-coords} \hA^+ := \bigcup_{0 < s < 1} A_{\cen} \times A \times \eta(s) \label{intro:A-cowt-coords}
 \cong   \{ (s_1, \ldots, s_{\ell+1}; z ) \in (\R_{>0})^{\ell+1} \times \R_{>0} \mid  s_1^{d_1} \cdots s_{\ell+1}^{d_{\ell+1}} < 1 \}. \ee Sometimes one fixes a value of $r$ (or $s$) and write $\hfh^r$ or $\hA^r$ for the corresponding slice.

	The Iwasawa decomposition states $\hG^+= \hK \times \hA^+ \times \hU$ for $\hA^+$ a subset of a finite-dimensional torus $\hA^e$ and $\hU$ a pro-unipotent group.  Hence $\xg = \hK \setminus \hG = \hA^+ \times \hU$ with the structure of a (contractible) topological space using the natural product topologies on the two factors (note that we have, in this way, avoided ever topologizing $\hG^+$). Moreover one can show that $\hGam$ acts continuously in this topology.  In fact, for a fixed $r$, the group $\hGam$ acts on the slice $\hG^r = \hK \times \hA^r \times \hU$, and so many of our results can be formulated just in this fixed $r$ context.  This is the point of view found in the works of H. Garland, but we prefer to work with all the slices for $r>0$ as it makes the formulation of some statements, which involve taking a limit to large $r$, a bit easier to state. On the other hand, it creates an extra non-compact parameter when $r \rr 0$ (since we do not allow $r=0$ in $\hG^+$). We do not attempt to add in a boundary piece corresponding to the $r=0$ case, since in the applications which we have in mind, it is the $r \rr \infty$ sitaution that is of interest. 
	
	Finally, let us emphasize that as in \cite{gar:ihes}, we work with a `complete' Kac--Moody group (as opposed to the minimal one of Tits), since our (pro-unipotent) topology can be specified by putting an explicit set of `Iwahori--Matsumoto' type coordinates on $\hU$. The quotient $\hU/ \hGam \cap \hU$ is then easily seen to be compact. Since we only complete the group in `one direction' (so the opposite unipotent $U^-$ is not pro-group), we need to restrict to the class of \textit{positive} parabolics, denoted $\ratp$, all of which are conjugate to certain standard parabolics, denoted $\stdp$ and indexed by subsets of $I$.

	\tpoint{Reduction theory of H. Garland} 	 In \cite{gar:ihes}, the Siegel set $\sie_{t, \Omega}:= \hA^+_t \times \hU_{\Omega}$ is considered, where $ \hA^+_t:= \{ h \in \hA^+ \mid h^{a_i} < t \} $ and where $\hU_{\Omega}$ is constructed using the Iwahori--Matsumoto coordinates mentioned above. Note that $\hU_{\Omega}$ is relatively compact. The projection of $\hA^+_t$ onto $\hA_{\cen}$ is unconstrained, but  on the other hand it has projection onto $\hA$ which is actually a \textit{bounded} set. This is the familiar observation from the theory of affine Weyl groups: chambers get replaced by alcoves.

	 The reduction theory of H. Garland can now be stated
	 \begin{nthm} \label{intro-thm-garland}  Fix $t < \frac{\sqrt{3}}{2}$ and $\Omega \supset [-1/2, 1/2]$ a subset of $\R$.  
	 	\begin{enumerate}
	 		\item \cite[Thm. 19.3]{gar:ihes} For any $x \in \xg$, there exists $\gamma \in \hGam$ such that  $x \gamma \in \sie_{t, \Omega}$.
	 		\item \cite[Thm. 21.6]{gar:ihes} Setting $\sie^r_{t, \Omega}:= \sie_{t, \Omega} \cap \hG^r$, there exists $r_0 > 0$ so that if $r  > r_0$ and $\sie^r_{t, \Omega} \gamma \cap \sie^r_{t, \Omega} \neq \emptyset$, then $\gamma \in \hGam \cap \widehat{P}$ for some (standard) proper parabolic $\hP$. 
	 	\end{enumerate}
	 	
	 \end{nthm}	
	
	\begin{nrem} Notice that the second statement is quite different from the finite-dimensional picture, and it led Garland to suggest that for $r \gg 0$, \textit{all} the points in $\xg^r$ should lie close to the `Borel--Serre boundary' when this space was eventually constructed. We will return to this point a bit later in \S \ref{intro-Arthur-comments}. \end{nrem}

	\tpoint{On the corners $\corn{\hA^+_P}$} The first step in our construction of $\bord{\xg}$ is to find an analogue of the bordification $\cor(\hA_P^+)$ for $\widehat{P} \in \ratp'$. Suppose $P=B$ is a Borel, then $\hA^+_P= \hA^+ \subset \R_{>0}^{\ell+2}.$  If we were to just take a standard completion of $\hA^+$ in $\mathbb{R}_{\geq 0}^{\ell+2}$, we obtain too many boundary faces (as there are only $2^{\ell+1}$ standard  parabolics), so  instead, mimicing \eqref{intro:finite-bord-Ap}, we now \textit{define} \be{} \cor(\hA_P^+) := \hA_P^+ \ \ \cup \ \ \bigsqcup_{Q\supset P} \; A(Q)_{^*P}, \ee and then introduce on this set a natural class of convergent sequences which specifies a topology. Alternatively, in the coweight coordinates \eqref{intro:A-cowt-coords},   $\hA^+ \subset (\R_{>0})^{\ell+1} \times \R_{>0}$ and one can show the completion $\corn{\hA^+}$ is obtained by allowing some of the first $\ell+1$-coordinates to go to $0$ (and not imposing any conditions on the last coordinate). One can also consider the sets $\hA^+_t$ and the interior of its closure in $\corn{\hA^+}$ which we denote as $\corn{\hA^+_t}$. Unlike the finite-dimensional case, this set is non-compact for two reasons: first, there is no bound on the central direction; and second, one does not have a limit as $r \rr 0,$ i.e. in the coweight coordinates we do not allow $s_1^{d_1} \cdots s_{\ell+1}^{d_{\ell+1}}=1$.

	\tpoint{On the bordification $\bord{\xg}$}  Following the `uniform construction' of Borel--Ji \cite[\S III.9]{borel-ji}, we first define the set \be{} \begin{array}{lcr} \bord{\xg}:= \xg \ \ \cup \ \ \bigsqcup_{P\in \ratp} \; \e(P) & \mbox{ where } & \e(P) = X_{M_P} \times \hU_P \end{array}, \ee and then introduce a topology on this set using a family of convergent sequences. If we define, as before, the corner attached to $P$ as the subset \be{} \label{intro:corner-def} \corn{P} = \xg \sqcup \, \, \sqcup_{Q \supset P} \e(Q), \ee we require that all of the convergence happens within a corner, so to speak. More precisely,

	\begin{nthm} \label{thm:intro-bord-loop-main} 
		\begin{enumerate}
			\item (`Embedding extension Lemma', Lemma \ref{lem:embedding-extension-lemma}, cf. \cite[Prop III.9.5]{borel-ji}) For each proper parabolic $P \in \ratp$, the inclusion $\xg \hookrightarrow \bord{\xg}$ extends to an open embedding $\iota_P: X_{M_P} \times \corn{\hA^+_P} \times \hU_P \rr \bord{\xg}$ where the former set on the left is equipped with the product topology, and where the image of $\iota_P$ is equal to $\corn{P}$.
			  
			\item (Proposition \ref{prop:Hausdorff} and Theorem \ref{thm:boundary-homotopy}) The space $\bord{\xg}$ is a Hausdorff and its boundary is homotopic to the affine rational Tits building attached to $\hG_{\Q}$  
			\item (Proposition \ref{prop:cont-Rg}) The action of $\hG_{\Q}$ on $\xg$ extends to a continuous action on $\bord{\xg}$.
		\end{enumerate}
	\end{nthm}	
		
	\tpoint{Affine Canonical Pairs} The theory of affine canonical pairs is the subject of \S \ref{sec:canonical-pairs}, and uses an affine analogue of the definition of degree of instability as in  \eqref{intro:deg-inst} as a starting point. Some care is required to define and relate $\rho$ and $\rho_P$, since in the affine context $\rho$ is no longer one-half the sum of all positive roots (this would be infinite). This is explained in \S \ref{subsub:rho-p-q}.

	 \begin{nthm}\label{thm:intro-canonical-pairs-affine} 
	 	\begin{enumerate}
		 	\item (Proposition \ref{prop:existence-of-minima}) The quantity $\deg_{\inst}(x):= \min_{P, \gamma \in \hG_{\Q}/P_{\Q}} \la \rho_P, H_B(x \gamma) \ra$ exists for each $x \in \xg$.
		 	\item (Theorem \ref{thm:canonical-pair}) If $\deg_{\inst}(x) < 0$, there exists a unique pair $(P, \delta)$ with $P$ a standard, proper parabolic and $\delta \in \hG_{\Q}/P_{\Q}$ such that $\deg_{\inst}(x) = \la \rho_P, H_B(x \delta) \ra$ and such that $P$ is maximal with respect to inclusion among all such pairs. This is called the canonical pair attached to $x$, and written as before $\cp(x)=(P, \delta)$.
		 	\item (Lemma \ref{lem:can-pair-crit} and Corollary \ref{cor:can-pair-criterion}) The pair $(P, \delta)$ is a canonical pair for $x$ if and only if
		 		\begin{enumerate}
		 			\item The projection of $x \delta $ onto $X_{M_P}$ is semi-stable for $M_P$
		 			\item The projection of $x  \delta$ onto $\hA^+_P$ lies in \be{} \label{affine-Ap-intro} \hA^+_{P, (0, 1)} = \{ a \in \hA^+_P \mid a^{\alpha} < 1 \, \mbox{ for } \alpha \in \Delta_P \mbox{ and } a^{\rho} < 1 \}. \ee
		 		\end{enumerate}
		\end{enumerate}
	\end{nthm}	
	
	\begin{nrem}
		\begin{enumerate}
			\item The existence of the minimum  follows from the techniques developed by Garland in his study of the Maass--Selberg relations for Eisenstein series, see \cite{gmp, gar:ms-2}.  Note that the notion of being semi-stable, i.e. that $\deg_{\inst}(x) \geq 0$, is not invariant under multiplication by elements from the central subgroup $\hA_{\cen} \subset \hA^+$.

			\item  The usual proofs in the finite-dimensional literature for uniqueness use Langlands combinatorial lemma (see \cite[Lemma 1.8.2]{friday-morning} and \cite[Proposition 2.4.1]{chau}) or, what probably amounts to the same thing, a distance minimizing argument, see Behrend \cite[Prop. 3.13]{beh} or Casselman \cite[Prop 5.2]{cass:partitions}. We found that a (naive?) generalization of Langlands combinatorial lemma fails in the affine setting, so instead our proof of uniqueness uses a workaround relying on two things: a) properties of (affine) orthogonal families (a notion that in finite dimensions is due to Arthur introduced, but that also implicitly used by Garland in the affine setting, see \cite{gar:ms-2}); and (b) an `inequality' relating two possible canonical pairs  \ref{prop:w-ineq} that was inspired by the work of Drinfeld--Gaitsgory \cite{drin:gai} and Schieder \cite[Thm. 4.2]{schieder} in the geometric context. It would be interesting to take our work as a starting point and find analogues of the Langlands combinatorial lemma or a Behrend-type uniqueness theorem for affine complementary polyhedron.
			
			\item Note the extra condition in \eqref{affine-Ap-intro} as opposed to the finite-dimensional picture described in Theorem \ref{thm:intro-canonical-pairs} and \eqref{def:AP-intro}. In finite dimensions, the condition $a^{\rho}<1$ follows from the condition $a^{\alpha} <1$ for $\alpha \in \Delta_P$, since $\rho$ is a (positive) sum of simple roots. This is no longer true in the affine case. Note that this again tells us that the $P$-unstable ends depend on the central direction (since $\rho$ does).
			
		\end{enumerate}
		
		\end{nrem} 
	 
	 \noindent As in finite dimensions, the above theorem allows us to define subsets $\xg(P, \delta)$ for $P \in \stdp'$ and $\delta \in \hG_{\Q}/P_{\Q}$  which fit into a  partition
	 \be{} \begin{array}{lcr} \xg = \xg^{ss} \cup \, \bigsqcup\limits_{P\in \stdp'} \qxg(P) & \mbox { where } \qxg(P):= \bigsqcup\limits_{\delta \in \hG_{\Q}/P_{\Q}} \xg(P, \delta) \end{array}. \ee Part (3) of the above theorem also gives us a fibration \be{}  \label{affine:f-P} f_{(P, \delta)}: \xg(P, \delta) \rr X^{ss}_{M_P} \times \hA^+_{P, (0, 1)}, \ee and as in the finite-dimensional case \eqref{XP:e},  for $P \in \stdp'$, one has  \be{} \label{affine-Xp-e} \xg(P, e) = X^{ss}_{M_P} \times \hA^+_{P, (0, 1)} \times \hU_P. \ee
	   
	  \tpoint{Relations to partitions in the theory of automorphic forms} \label{intro-Arthur-comments} In his work on trace formula, Arthur introduced a partition of the space $X_G$ (see \cite[Lemma 6.4]{arthur:trace}) attached to a truncation parameter $T$ (for an exposition see \cite[\S 3.6]{friday-morning}).  We have not studied the relation between the semi-stability partition and Arthur's (for $GL_n$ this can likely be deduced from the results in \cite{cass:partitions}), but note the following. The uniqueness in Theorem \ref{thm:intro-canonical-pairs-affine}(2) shows that if $x \in \xg(P,e)$ and $x \gamma \in \xg(P,e)$ for $\gamma \in \hGam$, then in fact $\gamma \in \hGam_P$. A similar result holds for Arthur's partition (see \cite[Lemma 3.6.1]{friday-morning}) when $T$ is `sufficiently regular' suggesting a comparison  in this regime (see \cite[Remark 2.5.2]{chau}).
	  
	  Assuming such a relation holds, let us further note the following: in the Arthur--Langlands partition, starting from $x \in X_G$, one first translates $x$ into a Siegel set and then, based on which walls the point $x$ is near (more formally, one uses the `Langlands partition' of $\mf{h}$ (see \cite[Ch. IV, \S 6]{borel:wallach}), one decides the parabolic to assign $x.$  Now we know (from an adaptation of \cite[\S 19]{gar:ihes} to the finite-dimensional case) that the process of moving $x$ into the Siegel set involves only \textit{one} minimization over the highest weight representation $V_{\rho}$, i.e. it just uses the Borel. On the other hand, our partition using semi-stability uses several minimizations, over all standard parabolics $P$ and the representations corresonding to $\rho_P.$ This suggests that perhaps there is a way to rephrase the notion of semi-stability entirely in terms of the Borel. In the finite-dimensional setting, and over function fields, this indeed seems possible, see \cite[Lemma 7.3.2]{drin:gai}.
	 
	  Now we turn to the loop setting. As remarked after Theorem \ref{intro-thm-garland}, if we follow the scheme that leads to the Arthur--Langlands partition described in the previous paragraph but also take into account the form of Garland's theory of parabolic transformations, we are led to expect that for $r$ large, the `core' of the partition, or that portion of $\xg$ \textit{not} assigned to a proper parabolic $P \subset \hG$, should vanish. In other words, we might suspect that $\xg^{ss} \cap \hg^r$ should shrink for $r$ large. We have not been able to prove a precise result along these lines. In fact, for any $r >0$ we produce elements in \S \ref{subsub:construction-semistable points} $\xg^{ss} \cap \hG^r.$ These points do however lie in the closure of the non-semistable locus leading us to speculate whether every point of $X_{\hG^r} \cap \xg^{ss}$ satisfies this property for large $r$.

	\tpoint{On arithmetic quotients}  	Turning to the quotients $\xg/ \hGam$ and $\bord{\xg}/\hgam$, for each $S \in \ratp'$, we may again define  $\pbord{S} \subset \bord{\xg}$ and obtain a partition \be{} \bord{\xg} = \xg^{ss} \sqcup \, \bigsqcup_{S \in \ratp'} \pbord{S}. \ee Each $\pbord{S}$ is equipped with a fibration with fibers homeomorphic to $\hU_P$ \be{} \mf{f}_S: \bord{\xg} \rr  X^{ss}_{M_P} \times \corn{\hA^+_{P, (0, 1)}} \ee extending \eqref{affine:f-P}, and if $P \in \stdp'$ this fibration degenerates to a product.  We also have  \be{} \bord{\xg}/\hGam = \xg^{ss}/ \hGam \, \, \cup \, \, \bigsqcup_{P\in \stdp'} \, \pbord{P} / \hGam_P. \ee Note that none of these pieces is now compact, though we can show that the boundary of $\bord{\xg}$ has compact quotient, see Prop \ref{prop:compact-boundary}.  For example, 
	for $P$ standard, $\pbord{P} / \hGam_P$ is homeomorphic to $X^{ss}_{M_P}/ \Gamma_{M_P} \times \corn{\hA^+_{P, (0, 1)}} \times \hU_P/\hGam_{U_P}$, and whereas the first and third factors are compact, the middle term is again non-compact. Imposing certain restrictions, one can again formulate a compactness result (see Proposition \ref{prop:compact-interior}) for all of $\bord{\xg}/\hGam$.
	
	\tpoint{Loose Ends} We collect here a few loose ends, which we hope to address at a future time.
	
	\begin{enumerate} 
		\item Our results on separability of the quotient $\bord{\xg}/\hGam$ are only conditional. They rest on a certain Siegel-type finiteness result that is straightforward for finite dimensional unipotent groups and whose extension to the pro-unipotent case we not yet established. 
		
		\item In finite dimensions, Saper \cite{saper} has described a $\Gamma$-equivariant `tiling' (see Definition 2.1 of \textit{op. cit.}) of the bordified space $\bord{X_G}$ with tiles $\mc{T}_S$ indexed by $S \in \ratp$. Among these, there is a distinguished central tile, denoted $\mc{T}_G$, which is a closed, $\Gamma$-invariant, codimension-zero submanifold with corners contained in $X_G$ that is compact modulo $\Gamma$. Further, the boundary faces of this central tile give rise to the other tiles in the decomposition. Each such tile is constructed by applying the geodesic action of Borel--Serre (see \cite[\S 3.2]{borel-serre}) corresponding to a suitable cone in $\cc(A_S)$ of the associated boundary face. As such, one finds a resemblance to our construction of the partition $\pbord{S} \subset \bord{X_G}$ described in \S \ref{subsub:intro-bord}. Saper also established the existence of a $\Gamma$-equivariant \textit{retraction} from $\bord{X_G}$ to $\mc{T}_G$ \cite[Theorem 6.1]{saper} and indicated a number of interesting applications of his result, one of which is that the quotient $\mc{T}_G/\Gamma$ serves as a realization of the Borel–Serre compactification inside of $X_G$. It would be interesting to investigate whether a similar $\hgam$-equivariant retraction might exist in our setting.		

		\item In the work of Arthur \cite{arthur:trace} and also in Saper \cite{saper}, one finds partitions and tilings indexed to certain \textit{truncation parameters} or \textit{tiling parameters} respectively. We have refrained from introducing corresponding \textit{semi-stability parameters} as they introduce another layer of notational complexity and were not needed for our present purposes. They can however be put in via a straightforward manner, and it seems that to make a comparison between our work and that of Arthur or Saper one may actually need this. Also, as Harder points out \cite[p.58]{harder-buch}, it is useful to have spaces slightly larger than the exact $P$-semistable loci (which can be defined through the introduction of these parameters) when studying the cohomology of arithmetic groups.
		
	\end{enumerate}

	\subsection{Acknowledgements} 
	
	We would like to thank Howard Garland for helpful discussions about this paper and related topics. The reader will quickly see how his ideas permeate many of the central facets of this work, and we thank him for his encouragement and for generously sharing his insights with us. We also owe a debt to Bill Casselman who both brought to our attention the work of P.-H. Chaudouard \cite{chau} on canonical pairs and made us aware of the similarity between Arthur's orthogonal families and Behrend's complementary polyhedron. We have also benefited from his article \cite{cass:partitions} and thank him for a discussion about that paper and related matters. We would also like to thank Lizhen Ji for his valuable comments on an earlier version of this manuscript.

	Both authors were supported by the M.V. Subbarao Professorship in Number Theory and  NSERC Grant RGPIN-2019-06112 while this paper was in preparation. During the final stages of writing, M.P. was supported by NSERC RGPIN-2025-05292. P.S. also gratefully acknowledges his time as a visitor at the School of Mathematical Sciences, NISER Bhubaneswar, during the final stages of this work, and thanks Brundaban Sahu and Sudhir Pujahari for facilitating the visit.

	\part{ Loop groups and their arithmetic quotients} 
	
	\section{Affine Kac-Moody algebras, \textit{aka} loop algebras}   \label{sec:affine-Lie-algebras}

	\newcommand{\As}{\mathsf{A}}
	\newcommand{\dv}{\check{d}}
	\newcommand{\deltav}{\check{\delta}}
	\newcommand{\hfg}{\widehat{\mf{g}}}
	
	\subsection{Affine Lie Algebras: basic constructions} 
	\label{sub:affine-basic}
	\tpoint{Notations on GCMs}  \label{gcm} For a natural number $\ell$ and set \be{} \label{I:I_o} \begin{array}{lcr} I=\{ 1, \ldots, \ell+1 \}& \mbox{ and } & I_o= \{ 1, \ldots, \ell \} \end{array}. \ee  Let $\As=(a_{ij})_{i, j \in I}$ be an indecomposable generalized Cartan matrix (GCM) of affine type (see \cite[\S 1.1 and \S 4.3]{kac}). Note that we have $a_{ij} \leq 0$ for $i, j \in I$ with $i \neq j$ and also $a_{ii}=2$ for all $i \in I$.
	
	As $\As$ is of affine type, its null space is one dimensional and we write $\delta = (d_1, \ldots, d_{\ell+1})$ for the unique vector with integral, relatively prime, and positive entries in this space (cf. \cite[Theorem 4.8 (b)]{kac}). Similarly, the transpose ${}^t\As$ is again an indecomposable GCM of affine type and we define an analogous vector $\deltav=(\dv_1, \ldots, \dv_{\ell+1})$ in its null space (cf. \cite[Ch. 6]{kac} ). Note that one always has $\dv_{\ell+1}=1$. If $d_{\ell+1}=1$, the matrix $\As$ is said to be of \emph{untwisted} affine type and we shall restrict ourseleves to this case henceforth.
	
	 \label{gcm:o} For $J \subsetneq I$, set $\As(J)=(a_{ij})_{i, j \in J};$ if $J \neq I$, then $\As(J)$ is a (usual) Cartan matrix. It is of \textit{finite type} (i.e. positive definite).

	\tpoint{Lie algebras attached to $\As$ } \label{subsub:loop}  To any \textit{finite or affine} GCM $\As$ as above\footnote{It is automatically \emph{symmetrizable}, cf. \cite[Lemma 4.6]{kac}}, one attaches a Kac-Moody Lie algebra $\mf{g}(\As)$ over $\R$ as follows (cf. \cite[\S 2.1]{kac}): it has generators $e_i, f_i, \av_i \, (i \in I)$ subject to the relations: \be{km-rel} \begin{array}{lcccr} [\av_i, \av_j]=0, & [e_i, f_j] = \delta_{ij} \av_i, & [ \av_i, e_j] = a_{ij} e_j, & \ [\av_i, f_j] = - a_{ij} f_j & \text{ for } i, j \in I  \\ 
		(\mathrm{ad}\, e_i)^{1-a_{ij}}  e_j = 0, &(\mathrm{ad}\, f_i)^{1-a_{ij}}  f_j = 0 & & \text{ for } i, j \in I, \, i \neq j.    \end{array} \ee
	Let  $\mf{h}(\As)$ be the linear span of $\av_i \, (i \in I)$; it is an abelian Lie subalgebra of $\mf{g}(\As)$ and $\mf{n}_+(\As)$ (resp. $\mf{n}_-$) be the Lie algebras generated by $e_i, i \in I$ (resp. $f_i, j \in I$). 	Unless otherwise mentioned, we fix $\As$ an untwisted affine Cartan matrix as above and $\hfg:= \mf{g}(\As)$ and $\mf{g}:= \mf{g}(\As(I_o))$, $\widehat{\mf{h}}:= \mf{h}(\As)$, $\mf{h}:= \mf{h}(\As(I_o))$, \textit{etc.}  In other words, the `hat' will refer to the affine version of the corresponding unadorned construction. We identify $\mf{g} \subset \hfg$ as a Lie subalgebra in the natural way.

	\tpoint{On the structure of $\hfh$} \label{subsub:long-coroots} Defining  $\vthv:= \sum_{i \in I_o} \dv_i \av_i  \in \mf{h}$ one verifies using (\ref{km-rel})  that \be{}\label{def:cc} \cc := \sum_{i \in I} \dv_i \av_i= \av_{\ell+1} + \vthv \in \hfh \ee  is a central element in $\hfg.$
	Under the identification $\mf{g} \subset \hfg$, we have $\mf{h}  \subset \hfh$, and in  fact, \be{h:ho} \hfh = \mf{h} \oplus \R \av_{\ell+1} = \mf{h} \oplus \R \cc.\ee 
	
	 \noindent Write $\la \cdot, \cdot \ra: \hfh  \times (\hfh)^{\ast} \rr \R$ for the dual pairing and define $a_i \in (\hfh)^{\ast} \, (i \in I)$ by specifying \be{a_i} \la h, a_i \ra e_i = [h, e_i] \text{ for } h \in \mf{h}. \ee Very often, we adopt the `functional' notation and write \be{} \label{functional:pair}  \lambda(h):= \la h, \lambda \ra \text{ for } h \in \hfh \text{ and } \lambda \in (\hfh)^\ast. \ee 	
	 Thus we have $a_{ij} = \la \av_i, a_j \ra$ for $ i, j \in I.$ Note the elements $a_i$ are not linearly independent:  indeed, if we define, with a slight abuse of notation, the element in $\hfh^*$ called the \textit{minimal imaginary root} \be{} \label{def:delta} \begin{array}{lcr}  \delta= a_{\ell+1}+ \vth &\mbox{ with } &\vth  = \sum_{i \in I} d_i a_i, \end{array} \ee then $\la h, \delta \ra = 0$ for all $h \in \hfh.$ To remedy this issue, we need to extend $\hfh$ futher as we now explain.

	\tpoint{Extended affine Lie algebras $\hfg^e$}  \label{subsub:der}  For a tuple $(n_i)_{i \in I}$ of non-negative (resp. non-positive integers), let $\mf{g}( \{ n_i \}) \subset \mf{g}$  be the space spanned by \be{e-comm} [e_{i_1}, [ e_{i_2}, \ldots, [e_{i_{r-1}}, e_{i_r} ] \ldots ] \, ] \ \ (resp., \,  [f_{i_1}, [ f_{i_2}, \ldots, [f_{i_{r-1}}, f_{i_r} ] \ldots ] \, ] \ee where each $e_j$ (resp. $f_j$) occurs $| n_j |$-times in the above expressions. Define the derivation $\dd_{\ell+1}:= \dd$ of $\hfg$ by requiring $\dd$ act as the scalar $n_{\ell+1}$ on $\mf{g}( \{n_i\})$, and set \be{g-ext} \hfg^e = \hfg \rtimes \R \dd \text{ and } \hfh^e = \hfh \oplus \R \dd = \R \cc \oplus \mf{h} \oplus \R \dd. \ee We continue to write $\la \cdot, \cdot \ra: \hfh^e  \times (\hfh^e)^{\ast} \rr \R$ for the dual pairing and extend $a_i$ above to elements of $(\hfh^e)^{\ast}$ using (\ref{a_i}). One can check that $\la \dd, a_i \ra =0$ for $i \in I_o$ and $ \la \dd, a_{\ell+1} \ra = 1. $ So our inclusion of $\dd$ ensures that the  $\{ a_i \}_{i \in I}$ are now linearly independent elements of $(\hfh^e)^{\ast}.$ Defining the \textit{fundamental weights} $\lambda_j \in (\hfh^e)^{\ast}$ for $j \in I$  by  \be{Lambda:j} \la \av_i , \lambda_j \ra = \delta_{ij} \text{ and } \la \dd, \lambda_j \ra =0, \ee we find that the elements $\delta, \lambda_1, \ldots, \lambda_{\ell+1}$ form a basis of $(\hfh^e)^{\ast}$ dual to the basis $\dd, \av_1, \ldots, \av_{\ell+1}$ of $\hfh^e$.

	\tpoint{Kac form on $\hfh^e$} \label{subsub:bilin} Set $\epsilon_{i}:= d_i \, (\dv_i)^{-1} \text{ for } i \in I$ and  define a symmetric, non-degenerate (see \cite[Lemma 2.1b]{kac})) bilinear form $(\cdot, \cdot)$ on $\hfh^e$ by  \be{sym:frm} (\av_i, h) = \la h, a_i \ra \epsilon_i \text{ for } i \in I,\, h \in \mf{h} \text{ and }  (\dd,\dd) = 0 . \ee Then $(\cdot, \cdot)$  induces an isomorphism $\nu: \hfh^e \stackrel{\sim}{\longrightarrow} (\hfh^{e})^\ast.$ Note \be{in:aij} \nu(\av_i) = \epsilon_i a_i  , \, \nu(\cc) = \delta, \, \text{ and } \nu(\dd) = d_{\ell+1} \Lambda_{\ell+1} \ee With respect to the decomposition of $\hfh^e$, the form $(\cdot, \cdot)$ satisfies \be{ext-form} && \begin{array}{lcr} (\cc, \av_i) =0  \, \, (i \in I) , & (\cc, \cc)=0 , & (\cc, \dd)=d_{\ell+1} , \end{array} \\  && 
	(\av_i, \dd) = 0 \, (i \in I_{o}), \ \ \text{ and }  
	(\av_{\ell+1}, \dd) = 1. \ee  Denote again by $(\cdot, \cdot)$ the induced bilinear form on $(\hfh^e)^{\ast}.$ One has  $(a_i, a_j) = \epsilon_ja_{ij}/\epsilon_i \text{ for } i, j \in I$, 
	\be{invt-form-dual}  (a_i, \lambda_{\ell+1}) &=& 0, \, (i \in I_o) ; \ \  (a_{\ell+1}, \lambda_{\ell+1}) =1 ; \ \ (\lambda_{\ell+1}, \lambda_{\ell+1})=0; \ee and also $(\delta, x) =0$ unless $x$ is a non-zero multiple of $\lambda_{\ell+1}$ in which case $(\delta, \lambda_{\ell+1})=1$.
	
	\newcommand{\dcox}{\check{\mathsf{h}}}

	\newcommand{\rtl}{\mathscr{Q}}
	
	\tpoint{Roots} \label{subsub:roots} For each $\varphi \in (\hfh^e)^{\ast}$ let $\hfg^{\varphi}:= \{ x \in \hfg \mid [ h, x] = \la h, \varphi \ra x \text{ for all } h \in \hfh^e \}.$ The set of all non-zero $\varphi$ such that $\hfg^{\varphi} \neq 0$ will be called the roots of $\hfg^e$ \footnote{ We informally refer to these as the roots of $\hfg$, though without the inclusion of $\dd$, this is not a very useful notion}  and denoted by $\rts$. It is an affine root system (in the language of \cite[\S 2]{mac:eta}) with a base of simple roots $\Delta:= \{ a_1, \ldots, a_{\ell+1} \} \subset \rts$.  The root lattice is defined as \be{} \label{root-lattice} \rtl:= \mathrm{\Span}_{\zee} \rts,\ee and one defines $\rts_+$, the set of positive roots, as the set of non-negative integral linear combinations of $\Delta$, and $\rts_- = - \rts_+.$  Let $\rts_o:= \rts(I_o)$ denote the set of roots $\varphi$ such that $\hfg^{\varphi} \subset \mf{g}$. Note that $\{ a_1, \ldots, a_{\ell} \}$ form a set of simple roots of $\mf{g}$, and so postive and negative roots for $\mf{g}$ can be defined as the restriction of the corresponding notions from $\hfg$. More explicitly, $\rts= \rts_{re} \sqcup \rts_{im}$ where $\rts_{re}$ denotes the real roots and $\rts_{im}$ denotes the imaginary roots. The imaginary roots are easy to describe: $\rts_{im}  = \{ n \delta \mid n \in \zee \setminus \{ 0 \} \}$ and as for the real roots, we have $\rts_{re} = \{ \alpha + n \delta \mid \alpha \in \rts_o, n \in \zee \}.$ 
	
		\newcommand{\rtld}{\check{\mathscr{Q}}}
	\newcommand{\rtlv}{\check{\rts}}
	\tpoint{Coroots} \label{subsub:coroots} For $a_i \in \Delta$ we have defined the coroot $\av_i \in \hfh.$ We extend this and associate a coroot $a^{\vee} \in \mf{h}$ to  $a \in \rts$ as follows. Define $x_{a_i}:= \epsilon_i^{-1} \av_i$ for $i \in I$; and for $\varphi:= \sum_{i \in I} c_i a_i \in \rts_{re}$, set $x_{\varphi}:= \sum_{i \in I} c_i x_{a_i}.$ Then let  $ \check{a}:= \frac{2}{(a, a)}x_a \text{ for } a \in \rts_{re} ;$ if $a = a_i$, the two meanings of $\av_i$ coincide. 

	\newcommand{\reg}{\mathrm{reg}}

	\tpoint{The element $\rho$} \label{wts:def} \label{subsub:weights}

	We have already introduced the fundamental weights  $\lambda_i \, (i \in I)$ in (\ref{Lambda:j})  and the minimal, imaginary root $\delta$ at the end of \S \ref{subsub:loop}. Define also \be{} \label{def:rho} \rho :=  \sum_{i \in I}  \lambda_i. \ee It satisfies $\la \av_i, \rho \ra =1$ for all $i \in I$, $\la \dd, \rho \ra =0.$

	\tpoint{Coweights} \label{subsub:coweights} 
	
	Note that $(\mf{h}^e)^*$ has a basis $\{ a_i,  \lambda_{\ell+1} \}_{i \in I}$ and a non-degenerate form $(\cdot, \cdot)$  described in \S \ref{subsub:bilin}. 
	This allows us to construct a dual basis, say $\{ \edv_i, \cpsi \}_{i \in I}$ of $(\mf{h}^e)^*$, characterized using the isomorphim $\nu: \hfh^e \rr (\hfh^e)^*$ by  \be{} \begin{array}{lccccr} \left( \nu(\edv_i) ,   a_j \right) = \delta_{ij} , & ( \nu(\edv_i), \lambda_{\ell+1})=0, & (\nu(\psi), a_i)=0, & \mbox{ and } (\nu(\psi), \lambda_{\ell+1}) = 1 & \mbox{ for } i, j \in I. \end{array}. \ee 
	One can check:  $\check{\psi}= \cc$ , $ \hfh^e \cong  \R \cc \oplus \Span_{\R} \{ \lv_1, \ldots, \lv_{\ell+1} \}$, and $\edv_{\ell+1}$ is just $\dd$ (see Example \ref{example:sl2-epsilon}).

	\tpoint{The subalgebras $\hfh(J)$} \label{subsub:h(J)-coweights} Suppose now that $J \subsetneq I$. Then we let \be{} \hfh(J):= \hfh^+(J) = \Span_{\R} \{ \av_j \}_{j \in J}. \ee It it a Cartan subalgebra of the (finite dimensional) Lie algebra $\mf{g}(\As(J))$, and $\Delta(J):= \{ a_j \}_{j \in J}$ will form a basis of simple roots. Using this, one may define the collection of fundamental coweights $\omv^J_j$ for $j \in J$ which satisfy $\la a_k, \omv^J_r \ra = \delta_{kr}$ for $k, r \in J$. It is well known (see \cite[Lemma 1.1]{ragh} for a simple proof) that $\{\omv_j^J\}_{j \in J}$ form a basis of $\hfh(J)$ as well, and that $\omv^J_k$ can be written as a non-negative linear combinatation of $\{ \av_j\}_{j \in J}$. Similarly we can define the collection of fundamental weights $\omega^J_j$ for $j \in J$, characterized by the condition that $\la \av_d, \omega^J_i \ra = \delta_{di}$ for $d, i \in J$. Again, it is known each of which can be expressed as a non-negative linear combination of $\Delta(J)$, say \be{} \label{omega_J} \omega^J_d = \sum_{j \in J} c_j a_j \mbox{ with } c_j \geq 0 \mbox{ for } d \in J.  \ee Now the above is a relation between functionals on $\hfh(J)$, but the right hand side naturally as defining a functional on $\hfh^e$. It is \textit{not} the extension by $0$ of the left hand side from $\hfh(J)$ to $\hfh^+$.\be{} \label{omega_J} \omega^J_d = \sum_{j \in J} c_j a_j \mbox{ with } c_j \geq 0, \ee where the above is regarded as a relation between functionals on $\hfh(J)$. On the other hand, we can view the right hand side naturally as defining a functional on $\hfh^e$

	\tpoint{An affine positivity Lemma} The following is a straightforward generalization to the affine setting of Raghunathan \cite[Lemma 1.1]{ragh} (see also \cite[Lemma 4.9]{schieder}).
	
	\begin{nlem} \label{lem:ragh} Let $d \in I$. Then for any $J \subsetneq I$, there exists a decomposition  \be{} \label{ragh:dec}  \lambda_d = \sum_{j \in J} c_j a_j +  \mu,  \ee where
		$c_j$ are non-negative (rational) numbers and $\mu \in (\hfe)^*$ is such that $\la \mu, \av_i \ra \geq 0$ for all $i \in I$ and additionally $ \la \mu, \av_i \ra =0$ for $i \in J.$   \end{nlem} 
	
	\begin{proof} If $d \notin J$, we are done by choosing all $c_j=0$ and $\mu= \lambda_d$. So, assume that $d \in J$. Keeping the notation from the previous paragraph and as mentiond there, the difference $\mu:= \lambda_d - \omega^J_d$, regarded as a fucntional of $\hfh^e$, is in general non-zero. By construction,  $\la  \lambda_d - \omega^J_d, \av_i \ra =0$ for $i \in J.$ On the other hand, if  $i \notin J$, then as  $d \in J$, we have $\la \lambda_d, \av_i \ra =0.$ Hence we have  \be{}  \la \mu, \av_i \ra = \la \lambda_d - \omega^J_d , \av_i \ra = - \la \omega^J_d, \av_i \ra = - \la \sum_{j \in J} m_j a_j, \av_i \ra \geq 0, \ee since $\la a_j, \av_i \ra \leq 0$ for $j \in J, i \notin J.$  \end{proof}

	\tpoint{Example: $J=I_{o}$} Let us write (in the notation of \S \ref{subsub:h(J)-coweights} $\omega^{o}_i:= \omega^{I_o}_i$ for $i \in I_o$. For $j \in I_o$, we have $\lambda_j -\omega^{\ell+1}_j$ restricts to $0$ on $\Span\{ \av_1, \ldots, \av_{\ell}, \dd \}$, so this difference is determined by its value on $\av_{\ell+1}$, and this can be computed, using \eqref{def:cc} as follows: \be{} \la \lambda_j-\omega^{o}_j, \av_{\ell+1} \ra  =  -\omega^{o}_j(\av_{\ell+1})= \dv_j \mbox { if } j \in I_o. \ee

\begin{nlem} \label{lem:rho-vs-rho-classical} For $j \in I_o$,  $\lambda_j =\omega^{o}_j + \dv_j \lambda_{o}.$ Letting $\dcox:= 1 + \sum_{i=1}^{\ell} \check{d}_i$ and $\rho_o:=\omega^{o}_1+ \cdots +\omega^{o}_{\ell}$,  \be{} \rho = \rho_{o} +  \dcox \, \lambda_{\ell+1}.\ee   \end{nlem} 

	\renewcommand{\tt}{\mathsf{t}}
\newcommand{\affW}{\widehat{W}}

	\tpoint{Affine Weyl groups} \label{subsub:Weyl-group}  Let $W:= W(\As)$ be Weyl group associated to $\As$ (see \cite[\S 2.1.3]{pat-pus} and references therein). The latter is a Coxeter group whose generators we shall write as $S:= \{ s_i \}_{i \in I}$. Denote by $\ell: W \rr \mathbb{N}$ the length function. \label{s:pb-weyl} For each subset $J \subset I$, we define $W(J)$ as the group generated by $\{ s_i \}_{i \in J}$. If $J \subsetneq I$, then $W(J):= W(\As(J))$ is finite and contains an unique element $w_J$ of maximal length. There exist (see \cite[Lemma 1.3.3]{friday-morning}) a unique set $W^J \subset W$ of coset representatives for $W/ W(J)$ which have minimal length:  \be{W^J} \ell(w_1 w_2) =\ell(w_1) + \ell(w_2) \text{ for any } w_1 \in W^J, w_2 \in W(J). \ee For any $w \in W^J,$ we also have $w (a_j) > 0 \text{ for all } j \in J. $ In fact this condition is equivalent to the length-minimizing property of representatives in $W^J$. 

	For an untwisted $\As$, we usually write $\affW:= W(\As)$ and call this the affine Weyl group. The group $\affW$ may be viewed as a  subgroup of $\End_{\R}(\hfh^e)$ generated by the reflections:  \be{w:act:h}  s_i (h) = h - \la h, a_i  \ra \av_i  \text{ for } h \in \hfh^e \text{ and all } i \in I. \ee

	\newcommand{\con}{\mathscr{C}}
	\tpoint{The Tits cone $\hfh^+$} \label{subsub:tits-cone}  Define the antidominant cone $\mathscr{C} \subset \hfh^e$ as $\con:= \{ X \in \hfh^e \mid \la a_i, X \ra < 0 \}$. The \emph{Tits cone} $\hfh^+ \subset \hfh^e$ is defined as  $\ti := \bigcup_{w \in \affW} w(\con)$. We have (see  \cite[Proposition 1.9(a)]{kac-pet})   
	 \be{ti:exp}\label{def:hfh-r} \begin{array}{lcr} \ti=  \{ \lv\in \hfh^e | \la \lv, \delta  \ra < 0 \}. & \mbox{ For a fixed } r \mbox { set }  &   \hfh^r := \{ X \in \hfh^e \mid \la \delta, X \ra = -r \}. \end{array} \ee This subset is $\affW$-invariant. We can describe these sets more explicitly both in terms of the coroot basis $\{ \av_i, \dd \}_{i \in I}$ and the coweight basis $\{ \edv_i, \cc \}_{i \in I}$. Indeed, using the properties in \S\ref{subsub:loop} and \S \ref{subsub:der}, 
	 \be{} \label{tits-cone-desc} \hfh^+ &=& \{ c_1 \av_1 + \ldots c_{\ell+1} \av_{\ell+1} - r \dd \} = \{ H - r \dd \mid H \in \hfh, \, r > 0 \} \\ &=&  \{ \sum_{i \in I} p_i \edv_i + m \cc \mid \sum_{i \in I} p_i d_i < 0  \}. \ee

\newcommand{\mult}{\mathsf{m}}

\subsection{Enveloping algebras, Representations, and Hermitian structures.} \label{sec:representations}

\tpoint{Chevalley basis and $\zee$-forms} \label{subsub:enveloping-algebras} 
\label{subsub:chev-basis} Starting from a Chevalley basis of $\mf{g}:= \mf{g}(\As_o)$, Garland has constructed a corresponding basis for $\hfg$ also satisfying integrality properties as in the finite-dimensional case (see \cite[\S3 and Thm. 4.12]{gar:la}). The elements of this basis will be denoted as  \be{chev:aff} \Psi= \{ \xi_a \}_{a \in \rts_{re} } \sqcup \{ \xi_i(n) \}_{i \in I_o, n \neq 0} \sqcup \{ h_1, \ldots, h_{\ell+1} \} \ee where $h_i = \av_i$, $\xi_a \in \hfg^{a}$ for $a \in \rts_{re}$ and $\{ \xi(n) \}_{i \in I_o}$ is a basis of $\hfg^{n \delta}.$ Let $\hfg_{\zee}$ the $\zee$-span (which is a subalgebra) of $\hfg$. We also note that by direct computation that $\dd$ preserves $\mf{g}_{\zee}$ so that we may form again the semi-direct product $\mf{g}_{\zee}^e = \mf{g}_{\zee} \rtimes \zee \dd$.

For a (real) Lie algebra $\mf{s}$, let $\ue(\mf{s})$ denote its universal enveloping algebra.  Garland has  introduced a $\zee$-form $\ue_{\zee}(\hfg)$ of $\ue(\hfg)$ in \cite[\S 5]{gar:la} and provided for it a natural $\zee$-\textit{module} basis. Note that $\ue_{\zee}(\hfg)$ is stable under multiplication by the divided powers $\xi_a^{(n)}:= \frac{\xi_a^n}{n!}$ for $a \in \rts_{re}, n \geq 0$ and that $\ue_{\zee}(\hfg) \cap \hfg = \hfg_{\zee}$. This will be important when we want to attach groups to $\hfg$.

\tpoint{Representations $V^{\lambda}$} \label{subsub:V-lambda}

An element $\lambda \in (\hfh^e)^*$ will be said to be a \textit{dominant weight} if $\la \av_i, \lambda \ra \geq 0$ for all $i \in I$. It is said to be \textit{integral} if all of these values are integers, and is called \textit{normal} if at least one of them is non-zero (i.e. it is not a multiple of $\delta$ ). Write $\Lambda_+$ (resp. $\Lambda$) for the set of all dominant (resp. integral) weights.

Given any normal $\lambda \in \Lambda_+$, there exists a \emph{irreducible} highest weight module with highest weight $\lambda$ denoted $V^{\lambda}:=V^{\lambda}_{\R}$.  Moreover, \cite[Lemma 10.4]{gar:la}, for any $ v \in V^{\lambda}$, there exists a positive integer $r$ (depending on $v$) such that $e_i^r. v = f_i^r. v =0$ for $i \in I$. We write $\wts_{\lv}:= \wts(V^{\lv})$ for the set of weights of $V^{\lv}$ and $\Xi_{\lv}$ for the corresponding weight lattice. 

 Note that $\rtl \subset \Lambda,$ and in fact, we have (see \cite[Lemma 15.2]{gar:ihes})  $\rtl \subset \Xi_{\lambda} \subset \Lambda$ and in fact $\Xi_{\lambda} \subset \Lambda$ is of finite index. Moreover, if $\mu$ is any other normal element, there exists from \cite[Lemma 15.7]{gar:ihes} a positive integer $m$ such that $ \Xi_{m \mu} \subset \Xi_{\lambda}. $ In fact (see \cite[(15.19)]{gar:ihes}), if $\lambda_1, \lambda_2$ are two normal elements such that $\Xi_{\lambda_2} \subset \Xi_{\lambda_1}$, then there exists a strictly positive integer $m$ such that $\lambda_2 - m \lambda_1 \in \rtl$ and $m \lambda_1 \in \Xi_{\lambda_2}$. 

\tpoint{$\zee$-forms of $V^{\lambda}$} \label{subsub:hw-lambda-chev}

A basis $\Upsilon$ for a representation of $\hfg$ is said to be \emph{admissible} if it consists of weight vectors (for the action of $\hfh^e$). Such an admissible basis $\Upsilon= \{ \hw_0, \hw_1, \ldots \}$ of $V^{\lambda}$ is said to be \emph{coherently ordered} if whenever $\hw_i \in V^{\lambda}_{\mu}$ and $\hw_j \in V^{\lambda}_{\mu'}$ are such that $i \leq j$, then $\dep(\mu) \leq \dep(\mu')$ where $ \label{depth} \dep(\mu) = \la \mu, \dd \ra.$ From  \cite[Theorem 11.3]{gar:la}, there exists an admissible basis $\Upsilon$ of $V^{\lambda}$ such that if $V_{\zee}^{\lambda}$ denotes the $\zee$-span of $\Upsilon$, then $V_{\zee}^{\lambda}$ is invariant under the integral form $\uzg$. Such a basis is called an integral, admissible basis.

We shall also refer to $V_{\zee}^{\lambda}$ as a $\zee$ or \textit{Chevalley form for the representation} $V^{\lambda}$. Note that for any weight $\nu$ of $V^{\lambda}$, if $V^{\lambda}_{\nu}$ denotes the corresponding weight space, then $V_{\zee}^{\lambda} \cap V^{\lambda}_{\nu}$ is spanned by  vectors in $\Upsilon,$ \textit{ i.e.} $V^{\lambda}_{\zee} = \oplus_{\nu} ( V^{\lambda}_{\nu} \cap V_{\zee}^{\lambda}).$ A highest weight vector  $\hw_{\lambda} \in \Upsilon \cap V^{\lambda}_{\lambda}$ is called a \emph{primitive highest weight vector}.  For $\k$ any field, write $V^{\lambda}_{\k}:= \k \otimes_{\zee} V^{\lambda}_{\zee}$.

\tpoint{Hermitian structure} \label{subsub:hermitian}

Fix a $\zee$-form $V^{\lambda}_{\zee}$ as in \S \ref{subsub:hw-lambda-chev}. Then we recall  (see \cite[Thm. 12.1]{gar:la}) that $V^{\lambda}:= V^{\lambda}_{\R}$ admits a \emph{positive-definite}  inner product $\{ \cdot , \cdot \}$ (the non-trivial part of the above is to show that $\{ \cdot, \cdot \}$ is positive-definite without using integration) such that we have \begin{enumerate} 
	\item $\{ v, w \} = 0$ for any $v \in V^{\lambda}_{\mu, \R}$ and $w \in V^{\lambda}_{\mu', \R}$, for $\mu, \mu' \in \wts_{\lambda}$ with $\mu \neq \mu'$;
	\item $\{ v, w \} \in \zee$ for any $v, w \in V^{\lambda}_{\zee}$;
	\item $\{ \hw_{\lambda}, \hw_{\lambda} \} =1$ where $\hw_{\lambda}$ is a primitive highest weight vector.
	\item $\{ \xi_a v, w \} = \{ v, \xi_{-a} w \}$ for $v, w \in V^\lambda_{\R}$ and $a \in \rts_{re}$ and $\xi_a \in \Psi$ (see \eqref{chev:aff}) \end{enumerate}

\subsection{Parabolic subsets}  \label{sub:parabolic-subsets}

\newcommand{\pb}{\mathscr{P}}
\newcommand{\rs}{\rts}
\newcommand{\bor}{\mathscr{B}}
\newcommand{\sub}{\mathscr{X}}
\newcommand{\ps}[1]{\mathscr{#1}}

\newcommand{\spb}{\leftidx^{*}\pb}
\newcommand{\sqb}{\leftidx^{*}\qb}
\newcommand{\qb}{\mathscr{Q}}

\newcommand{\lh}{\mf{h}}
\newcommand{\simp}{\Delta}

\tpoint{Parabolic, closed, Borel subsets.} \label{subsub:parabolic-defs} A subset of the roots (of $\hfg^e$), say $\sub \subset \rs,$ is called \emph{closed} if $a, b \in \sub$ and $a + b \in \rs$ implies that $a +b \in \sub$. A subset $\sub \subset \rs$ is said to be \emph{parabolic} if $\sub \cup - \sub = \rs.$ A parabolic subset $\sub \subset \rs$ such that $\sub \cap - \sub = \emptyset$ will be called a \textit{Borel subset}. For any parabolic subset $\sub \subset \rs$, we define its semi-simple and nilpotent parts as \be{Psi:um} \sub^s:= \sub \cap (- \sub)    & \text{ and } &\sub^n:= \sub \setminus \sub_s \ee respectively. Also, for $\sub \subset \rs$, let $\la \sub \ra$ (resp. $ \la \sub \ra _+$, resp. $\la \sub \ra_-$) denote the set of all $\zee$-linear (resp. all $\zee_{\geq 0}$, resp. all $\zee_{\leq 0}$) combinations of elements from $\sub$ which also lie in $\rs$. Finally we say that two parabolic subsets $\sub_1, \sub_2$ are \textit{conjugate} if there exists $w \in \affW$ such that $w \sub_1 = \sub_2$. 

\newcommand{\lb}{\left[}
\newcommand{\rb}{\right]}
\tpoint{Standard parabolic subsets} \label{subsub:standard-parabolic} For $J \subset I$, write $\simp(J):= \{ a_i \in \simp \mid i \in J \}$. Write $\left[ \simp(J) \right]_{\pm} \subset \rts$ for the corresponding linear combination of roots with positive/negative integer coefficients, and define the standard parabolic subsets $\pb_{J}:= \rts_+\sqcup \lb \simp(J) \rb_- = \lb \simp \rb_+ \sqcup \, \lb \simp(J) \rb_-,$ and one may verify that \be{Psi:the-u} \pb_{J}^n &=& \{ a \in \rts \mid a = \sum_{j \in J} m_j a_j \text{ where } m_j \geq 0  \,  \text{ and } \,  m_j > 0 \, \text{ for } j \in J \} \\ \pb_{J}^s &=& \{ a \in \rts \mid a = \sum_{j \in J} m_j a_j \text{ where } m_j \in \zee \} = \lb \simp(J) \rb . \ee  In the case that $J$ is equal to  $J_i:= I \setminus \{ i \}$ for some $i \in I$, we simply write $\pb_i$ in place of $\pb_{J_i}$, and hence write $\pb_i^n$ and $\pb_i^s$ \textit{etc.} These are the \textit{maximal parabolic} subsets. For any $J \subsetneq I$, the set $\pb_J^s$ are the roots of a finite-type root system attached to the finite-type Cartan matrix $\mathsf{A}_J$.

\tpoint{Classification of parabolic subsets} For finite-type root systems, one can again define the notions of parabolic subsets, etc. and the notations above match those in \cite[Chap VI, \S 7]{bour}. Write $\As_o$ for a finite-type Cartan matrix with Lie algebra $\mf{g}_o$ having roots and chosen base of simple roots $\rts_o$ and $\simp_o$ respectively, the latter indexed by a set $I_o.$  The classification of parabolic subsets in finite-type asserts (see \cite[Chap. VI, \S 1.7 Prop. 20]{bour}) that every Borel subset of $\rts_o$  is conjugate to $\lb \simp_o \rb$, and, in fact, every parabolic subset $\pb_o \subset \rts_o$ is conjugate to $ \lb \simp_o \rb_+ \,  \sqcup \, \lb J_o \rb_- \,$ for $J_o \subset I_o$. 

The classification of parabolic root systems in affine type is more involved than in finite type, see  \cite[\S2]{kac-jk1}. For the present purposes, it suffices to focus exclusively on the following class of \textit{positive} \footnote{Replacing $\simp$ with $- \simp$ in the definition yields a class of negative paraboilc subsets; there are, however, also parabolic subsets which are neither positive or negative.} parabolic subsets (or just positive parabolic subsets for short) which are defined to be parabolic subsets also satisfying $ | \pb \cap \lb \simp \rb_- | < \infty. $ Based on \emph{op. cit.}, we conclude

\begin{nprop}   \label{prop:base-pb-aff}
Let $\pb \subset \rs$ be a positive parabolic subset. Then there exists $J \subset I$ such that $\pb$ is conjugate to $\pb_{J}:= \rts_+ \sqcup \lb \simp(J) \rb_-$.  

\end{nprop}

\newcommand{\f}{\mathrm{fin}}

\spoint \label{subsub:parabolic-correspondence} Let $\pb \subset \rs$ be a positive parabolic subset with decomposition $\pb = \pb^s \sqcup \pb^n$. For $\qb \subset \pb$  and $\spb \subset \pb^s$ arbitrary subsets, we define subsets of $\pb^s$ and $\pb$ respectively, \be{r:i} \begin{array}{lcr} r(\qb) = \qb \cap \pb^s &  \text{ and } & i(\spb) = \spb \cup \pb^n. \end{array} \ee

\begin{nprop} \label{sub-par} The maps $r$ and $i$ above induce a bijective correspondence \be{bij:pb} \{ \text{positive parabolic subsets } \qb \subset \pb \} \longleftrightarrow \{ \text{parabolic subsets } \spb \subset \pb^s \}. \ee Under this correspondence, we have \be{} r(\qb)^s = \qb^s \cap \pb^s = \qb^s & \text{  and  }&  r(\qb)^n = \qb^n \cap \pb^n  \\ i (\spb)^s = (\spb)^s & \text{ and } &  i(\spb)^n  = (\spb)^n \cup \pb^n .\ee 
\end{nprop}

\newcommand{\chir}{\chi^r}
\renewcommand{\emptyset}{\varnothing}
\newcommand{\vn}{\emptyset}

\tpoint{Parabolic decompositions of the Cartan } \label{subsub:cartan-decompositions}  For each $J \subset I$ and $r \in \R$, define  \be{} \label{def hatrJ} \hfh(J):= \hfh(\es, J) &:=& \Span \{ \av_i \mid i \in J \} \subset \hcf \text{ and } \\ \label{def:a_P} (\hcf^r)_J &:=&   \{ H' \in \hcf \mid  \la a_i ,  H' - r \dd \ra  =0 \text{ for } i \in J  \}. \ee When $J =\es$, then we set $\hfh(\es, \es):= \{ 0 \}$. Whereas $\hfh(\es, J) \subset \hfh$  is a subspace, $ (\hcf^r)_J \subset \hfh$ is just a set in general. Note that $ (\hfh^r)_J \supset (\hfh^r)_K$ when $J \subset K.$  Using the above definitions, we can define a  decomposition of $\hfh$ (not $\hfh^r$!).

\begin{nlem}  \cite[Lemma 9.1]{gar:lg2} \label{lem:garland-dec} For any $r \in \R$ and $J \subsetneq I$, one has a decomposition  \be{}  \label{hr:dec}  \hcf  = \hcf({\es, J}) + (\hfh^r)_J \ee with uniqueness of expression, \textit{i.e.} given $X \in \hfh$ and $r \in \R$, we may write  as $X = Y + Z$ with $Y \in \hfh(\es, J)$ and $Z \in (\hfh^r)_J$ in a unique way. \end{nlem}

\begin{nrem} 
	\begin{enumerate}
		\item Although the notation may not indicate it, both terms $Y$ and $Z$ depend on $r$ (see the example in \S \ref{subsub:example-sl2} below). To indicate this, for $X \in \hfh$, sometimes we write \be{} \label{J-dec-with-r} X - r \dd = (Y - r \dd) +  Z \ee
		\item The Lemma fails when $J = I$ in which case $\hfh(\es, I)= \hfh$.  
		
	\end{enumerate}
	
	The decomposition takes a simple form in terms of the coweight basis $\{ \lv_i, \cc \}$. Indeed, write $H \in \hfh^+$ as $H = \sum_{i \in I} c_i \edv_i + f \cc$  with all $c_i $ and $f \in \R.$ Then for $J \subsetneq I$ decomposing as above $H = H(J) + H_J$ with $\la H_J, a_i \ra =0$ for $i \in J$ and $H(J)$ in the span of $\{ \av_i \}$ with $i \in J$, if  \be{} \label{H(J):p-c} H(J) = \sum_{j \in J} p_j \omv^J_j \mbox{ then we must have } p_j=c_j. \ee
	Indeed, this follows from the properties of $\edv_j$: for $ j \in J$ \be{} c_j = \la H, a_j \ra = \la H(J) + H_J, a_j \ra = \la H(J), a_j \ra = \la \sum_{j \in J} p_j \omv^J_j, a_j \ra = p_j. \ee In other words, using the coweight coordinates, we have a simple way of performing the decomposition in the Lemma, as the coefficients of $H(J)$ are obtained from the coefficeints of $H$.

 \end{nrem}  Using the uniqueness in the Lemma, we can define natural projections \be{} \label{projection-r} \begin{array}{lcr} \pr^r(J): \hfh^r \rr \hcf({\emptyset, J}) & \text{ and } & \pr_J^r: \hfh^r \rr (\hfh^r)_J \end{array}. \ee  Setting  \be{} \label{h:plus} \begin{array}{lcr} \hfh^+:= \{ H - r \dd \mid H \in \hfh, r > 0 \}    & \text{ and } & (\hfh^+)_J:= \cup_{r > 0} (\hfh^r)_J, \end{array}  \ee the above maps can be put together to give projections  \be{} \label{projection:+} \begin{array}{lcr} \pr^+(J): \hfh^+ \rr \hcf({\emptyset, J}) & \text{ and } & \pr_J^+: \hfh^r \rr (\hfh^+)_J \end{array}. \ee 
 
It will be important for us sometimes to keep track of the degree derivation component, so we will also introduce the subset \be{} \label{hplus:d} \hfh^+_J:= \{ H - r \dd \mbox{ where } H \in \hfh, r > 0 \mid \la a_i, H - r \dd \ra =0 \mbox{ for } i \in J \}. \ee If we fix a value of $r > 0$, we write the corresponding subset as $\hfh^r_J$. Note that \be{} \hfh^r = \hfh^r_J + \hfh(\es, J) \ee follows at once from Lemma \ref{lem:garland-dec}.

\tpoint{Example: $\widehat{sl}(2)$}  \label{subsub:example-sl2} Let us consider a simple example for affine $\mf{sl}(2)$, i.e. $\As = \begin{pmatrix} 2 & -2 \\ -2 & 2 \end{pmatrix}$ so that $\hfh^e = \Span_{\R}\{\av_1, \av_2, \dd \}.$ Writing $\av:= \av_1$ and using the fact that $\av_2 = - \av + \cc$, every element in $\hfh^e$ has a unique decomposition \be{} X = X(m, x, r):= m \cc + x \av - r \dd = m \av_2 + ( m + x ) \av_1 -r \dd  \mbox{ for some } m, x, r \in \R. \ee 

\begin{itemize}
	\item Suppose $J = \{ a_1 \}$. Then the decomposition \eqref{hr:dec} amounts to grouping \be{} \label{sl2:a1}  n \av_1 + m \av_2 = \underbrace{(n-m) \av}_{\hfh(\es, I)}  +  \underbrace{m \cc}_{(\hfh^r)_I} . \ee This decomposition does not depend on $r,$ but we have $a_1(m \cc -r \dd)=0$ and $a_2 ( m \cc - r \dd ) = -r $. 
	\item Suppose now that $J = \{ a_2 \}.$ For $p, q \in \R$ we have $a_2( p \av_1 + q \av_2 ) = 2 (q - p ) $ so that \be{} (\hfh^r)_{\{ a_2\} } = \{ p \av_1 + q \av_2  \mid 2 (q-p) = r \}. \ee Hence, the decomposition \eqref{hr:dec} asserts that \be{} X(m, x, r) = - ( x + r /2 ) \av_2  + \left[ (m+ x) \av_1 + (m + x+r/2 ) \av_2 - r \dd \right]  \ee since $- (x+r/2) \av_2 \in \hfh(\es, \{ a_2\})$ and $(m+x)\av_1 + (m + x+r/2 ) \av_2  \in  (\hfh^r)_{\{a_2\}}.$ Note that \be{} a_1 ( (m+ x) \av_1 + (m + x+r/2 ) \av_2 - r \dd  ) = -r  . \ee 
	\item If $J= \{ a_1 , a_2 \}$ then $\hfh(\es, \{ a_1, a_2 \}) = \hfh$ and the only element in $(\hfh^r)_{ \{ a_1, a_2\} }= 0.$
	\item If $J = \es$, then $\hfh(\es, \es)= 0$ and $(\hfh^r)_{\es}= \hfh$. 
\end{itemize} 

\label{example:sl2-epsilon}

Alternatively, we can work with the coweight basis. First we determine $\edv_1, \edv_2, \cpsi$ to satisfy  \be{} && \begin{array}{lcr} ( \nu(\edv_1), a_i ) = \delta_{1, i} & \text{ and } &  (\nu(\edv_1), \lambda_{\ell+1}) = 0 \end{array} \\ && \begin{array}{lcr} ( \nu(\edv_2), a_i ) = \delta_{2, i} & \text{ and } & (\nu(\edv_1), \lambda_{\ell+1}) = 0 \end{array} \\ && 
\begin{array}{lcr} ( \nu(\cpsi), a_i ) = 0 & \text{ and } & (\nu(\cpsi), \lambda_{\ell+1}) = 1 \end{array}. \ee We can then solve these (uniquely) to obtain \be{} \begin{array}{lcr} \edv_1 = \frac{1}{2} \av_1 + \dd  , & \edv_2 = \dd, & \mbox{ and } \cpsi = \av_1 + \av_2 = \cc. \end{array} \ee
From here, we note that if $H= c_1 \edv_1 + c_2 \edv_2 + f \psi$, then \be{} H = (c_1/2+f) \av_1 + f \av_2 + (c_1 + c_2) \dd = (c_1/2) \av_1 + f \cc + (c_1 + c_2) \dd  \ee Using the above, we also conclude that \be{} && \begin{array}{lcr} H = \underbrace{\frac{1}{2} c_1 \av_1}_{H(J)} + \underbrace{ ( f \cc + (c_1+c_2) \dd )}_{H_{J}} = c_1 \omv^J  + H_J & \mbox{ when } & J = \{ a_1 \} \end{array} \\ && \begin{array}{lcr} H = \underbrace{\frac{1}{2} c_2 \av_2}_{H(J)} + \underbrace{( c_1/2 +f )\av_1 + ( f - c_2/2) \av_2 + (c_1+c_2) \dd )}_{H_{J}} = c_1 \omv^J  + H_J & \mbox{ when } & J = \{ a_2 \} \end{array}.  \ee

\tpoint{Relative variant} \label{subsub:relative-variants} Fix now $J \subset K  \subsetneq I$ and define the sub\textit{space} \be{}\label{def:h:J:K} \hfh(J,K) = \{H \in \hfh(\emptyset,K) | \, a_i(H) = 0  \mbox{ for } i \in J\}. \ee 

\begin{nlem} Suppose $J \subset K_1 \subset K_2 \subsetneq I$. Then we have a direct sum decomposition \be{} \label{h:J-K1-K2} \hfh(J, K_2)  = \hfh(J, K_1) \oplus \hfh(K_1, K_2). \ee For a fixed $r$ and $J \subset K \subsetneq I,$ we also have \be{} \label{dec:hrJ:K}(\hfh^r)_{J} =  \hfh(J,K)+(\hfh^r)_{K} \ee with uniqueness of expression in the above decomposition.
\end{nlem}
\begin{proof} For the first assertion,  notice that since $K_2 \subsetneq I$ we are in the setting of finite-dimensional root systems.  As for (2), first note that if $H_1 \in \hfh(J,K)$ and $H_2 \in (\hfh^r)_{K}$, then we have 
	\be{} a_i(H_1 +H_2 -r\dd)  
			= a_i(H_1)+a_i(H_2-r\dd)  = 0 \mbox { for } i \in J, \ee
	and hence we have the inclusion $\hfh(J,K) + (\hfh^r)_{K} \subset (\hfh^r)_J.$ To prove the reverse inclusion, use the splitting \eqref{hr:dec} $\hfh = \hfh(\emptyset,K)+ (\hfh^r)_{K}$ to write $H \in (\hfh^r)_{J} \subset \hfh$ as  $H = H_1 +H_2$ with $H_1 \in \hfh(\emptyset,K)$ and $H_2 \in (\hfh^r)_{K}$. Then from the definition of $(\hfh^r)_J$ and $(\hfh^r)_K$ and since $J \subset K$ we have, \be{} a_i(H-r\dd) = a_i(H_2-r\dd) = 0  \mbox{ for }  i \in J.\ee From this we get $H_1 \in \hfh(J,K)$, since we may write \be{} a_i(H_1) = a_i(H-r\dd) - a_i(H_2-r\dd) = 0  \mbox{ for } i \in J. \ee  Hence,  $(\hfh^r)_{J} \subset \hfh(J,K)+ (\hfh^r)_{K}$. The uniqueness assertion follows from that asserted in \eqref{hr:dec}.
\end{proof}
 \spoint \label{subsub:compatability} The following compatibility result will be used often in the sequel. 
 
\begin{nlem} \label{lem:relative-lie-algebra}  Let $J \subset K \subsetneq I$. Then we have a direct sum decomposition
	 \be{} \label{eq:relative-lie-algebra}  \begin{array}{lcr} \hcf(\emptyset, K) = \hcf(K)_{\leftidx^{\ast} J} \oplus \hcf(\emptyset, J) & \text{ where } & \end{array}  \\ \hcf(K)_{\leftidx^{\ast} J} := \{ h \in \hcf({\emptyset, K}) \mid a_i(h)=0, i \in J \} = \hfh(J, K) .  \ee   Denoting by $\pr^K(J): \hcf(\emptyset, K) \rr \hcf(\emptyset, J)$ the projection induced by the decomposition (\ref{eq:relative-lie-algebra}), the following diagram commutes: 
	\be{} \begin{tikzcd} \hcf^r \arrow[d, "\pr^r(J)"] \arrow[r, "\pr^r(K)"] & \hcf(\emptyset, K) \arrow[ld, "\pr^K(J)"] \\ \hcf(\emptyset, J)  \end{tikzcd} \ee \end{nlem}
\begin{proof}Let $X \in \hcf-rD$ and write $X=Y_J + Z_J$ with $Y_J \in \hcf_J^r$ and $Z_J \in \hcf(\es, J)$ and similarly $X= Y_K + Z_K$. Then $\pi_K^r(X)= Z_K$ and $\pi_J^r(X)= Z_J$.  Use the Lemma \ref{lem:relative-lie-algebra} to write $Z_K = Y_{\leftidx^{\ast} J} + Z'_J$, so that $$ X = \left( Y_K + Y_{\leftidx^{\ast} J} \right)+ Z'_J.$$  For $i \in J \subset K$ we have $$ a_i( Y_K + Y_{\leftidx^{\ast} J} ) = a_i(Y_K) = 0 $$ and so by the uniqueness of the decomposition above, we are done. 	\end{proof}

\subsection{The corner $\cor(\hfh^+)$: construction } \label{sub:bord-h} 

In this section, we will explain a method of enlarging the space $\hfh^+$ to a space $\cor(\hfh^+)$ by attaching various corners. We also introduce a topology on $\cor(\hfh^+)$ and explain some basic properties of this space. To motivate our constructions, we first explain how the construction works in finite type.

\tpoint{Warm up: finite-dimensional corners} 
Let $\As_o$ be a finite type $\ell \times \ell$ Cartan matrix with corresponding Lie algebra  $\mf{g}_o$ and Cartan subalgebra $\mf{h}_o \cong \R^{\ell}.$ Denote by $\av_i, i \in I_{o}$ the simple coroots and $\omv_i, i \in I_{o}$ the fundamental coweights. Both $\{ \av_i \}$ and $\{ \omv_i \}$ are bases of $\mf{h}_o$. For $J \subset I_{o}$, let $\mf{h}(\es, J)$, the finite dimensional analogue of \eqref{def:h:J:K}, be the span of $\av_j, j \in J$; it can be identified with as a Cartan subalgebra of minimal parabolic of the Lie algebra attached to the cartan matrix $\As_{o}(J)$. We shall also write $\{ \omv^J_j \}$ for the basis of $\mf{h}(\es, J)$ consisting of fundamental weights, i.e. $\la \omv^J_j, a_k \ra = \delta_{jk}$ for $j, k \in J$ and $a_k$ a simple root of $\mf{g}(\As_{o}(J))$ which can, in a natural way, be identified with a simple root of $\mf{g}(\As_o)$ that is denoted by the same letter.  

Let $A_o:= (\R_{>0})^{\ell}$ and consider the exponential map (which is a homeomorphism)  \be{} \label{c(A_o):dec} \begin{array}{lcr} \exp: \mf{h}_o \rr A_o, & \mbox{ sending } & \sum_{i \in I_o} \, c_i \omv_i \mapsto (e^{c_i})_{i \in I_o} \end{array}. \ee The closure of $A_o$ in $\R^{\ell}$ is equal to $(\R_{\geq 0})^{\ell}$; we shall denote it by $\cor(A_o)$ and it is easy to see  \be{} \label{finite-corner} \cor(A_o) = A_o \bigcup \, \sqcup_{J \subsetneq I_{o}} \, A_{o}(\es, J) \times o_J \ee where the symbol $A_{o}(\es, J) \times o_J$ designates the set \be{} \{ (s_1, \ldots, s_\ell) \mid s_k=0, k \notin J \mbox{ and } s_j > 0 \mbox{ for } j \in J  \}. \ee Each $A_o(\es, J) \times o_J$ is homeomorphic to $\mf{h}_o(\es, J)$ under a map,  denoted  $\exp_J,$ and which sends $H(J) \in \mf{h}_o(\es, J)$, say $H(J)=\sum_{i \in J} c_i \omv^J_i$ to the tuple $(s_k)$ where $s_i=e^{c_i}$ for $i \in J$ and $s_i=0$ otherwise.  If we now define the \textit{corner attached to $\mf{h}_o$} as \be{} \cor(\mf{h}_o):= \mf{h}_o \, \bigcup \, \sqcup_{J \subsetneq I_o} \mf{h}(\es, J), \ee then we obtain a set theoretic bijection  \be{} \label{c(exp)} \cor(\exp): \cor(\mf{h}_o) \rr \cor(A_o) \ee  obtained by patching together the various $\exp_J$.

\tpoint{Topology on $\cor(A_o$)}  The topology on $\cor(A_o)$ is the one inherited from $\R^{\ell}$, but we would like to describe what form it takes on $\cor(\mf{h}_o)$ under the isomorphism \eqref{c(exp)}. To do this, let us rexamine the topology on $(\R_{\geq 0})^{\ell}$ starting from the decomposition $\cor(A_o) =(\R_{\geq 0})^{\ell}$ given in  \eqref{c(A_o):dec}. First, on the subsets $A_o$ and $A_{o}(\es, J) \times o_J$ for $J \subsetneq I_{o}$, we impose the natural (Euclidean) topology. Then, we specify how a sequence of elements from one of these pieces converges to another piece, \textit{i.e.} we need to specify how a sequence $ \{ h_n \}$ with $h_n \in A_o$ converges to an element $h_{\infty} \in A_{o}(\es, J) \times o_J$,  how a sequence $ \{ h_{n} \}$ with $h_{n} \in A_o(\es, J) \times o_J$ converges to an element $h_{\infty} \in A_{o}(\es, K) \times o_K$ for $J \supset K$, etc. We leave it as an exercise in the notion of Moore--Smith convergence (see \S \ref{sub:moore-smith} ) that imposing the following convergence class will  recover the (subspace) topology on $\cor(A_o)$ obtained from the Euclidean topology on $\R_{\geq0}^{\ell}$ : 
\begin{enumerate}
	\item A sequence $h_n = (s_{n, 1}, \ldots, s_{n, \ell})$ in $A_o$ converges to $h_{\infty}= (s_{\infty, i} )_{i \in I} \in A_o(\es, J) \times o_J$, i.e. $s_{\infty, j} >0$ for $j \in J$ and is zero otherwise, if $s_{n, j} \rr s_{\infty, j} >0$ for $j \in J$ and $s_{n, k} \rr 0$ for $k \notin J$; equivalently, writing $c_{n, i}= \log s_{n, i}$, we have that $c_{n, j} \rr c_{\infty, j}:=\log(s_{\infty, j}) \in \R $ and $c_{n, k} \rr - \infty$ for $k \notin J$.
	
	\item For $I_{o} \supsetneq J \supset K$,  a sequence $h_n = (s_{n, 1}, \ldots, s_{n, \ell}) \in A_o(\es, J) \times o_J$ converges to an element $h_{\infty}:= (s_{\infty, d}) \in A_o(\es, K) \times o_K$ (i.e., $s_{\infty, d} > 0$ for $d \in K$ and is zero otherwise) if $s_{n, i} \rr s_{\infty, i}$ for $i \in I.$ 
\end{enumerate}

We can now translate this definition of the topology via convergent sequences over to $\mf{h}_o$ by explaining how a sequence of points in $\mf{h}_o$ converges to one in $\mf{h}_o(\es, J),$ how a sequence in $\mf{h}_o(\es, J)$ converges to $\mf{h}_o(\es, K)$ with $K \subset J$, etc.  Writing  $H_n \in \mf{h}_o$ as $H_n = \sum_{i \in I} c_{n, i} \omv_i$, we will say that $H_n$ converges to $H_{\infty} \in \mf{h}_o(\es, J)$  if $c_{n, i} \rr - \infty$ for $i \notin J$ and $c_{n, j} $ converge to finite values, say $c_{\infty, j} $ for $j \in J$; in this case, we also must have \be{} \label{H:inf-def} H_{\infty}:= \sum_{i \in J} c_{\infty, j} \omv^J_j, \ee where we now  a switch from $\omv_j$ to $\omv_j^J$ as is required from the definition of \eqref{c(exp)}. A similar procedure specifies how elements $H_n \in \mf{h}_o(\es, J)$ converges to an element $H_{\infty} \in \mf{h}_o(\es, K)$ if $K \subset J$. With this topology, we can verify that $\cor(\exp)$ is a homeomorphism. In particular, since $\cor(A_o)$ is clearly homotopic to its interior $A_o$, the same relation holds between $\cor(\mf{h}_o)$ and $\mf{h}_o$. 

We actually need a slight reformulation of this procedure of topologizing $\corn{\hfh_o}$ which is useful in the affine case.  To state it, note in analogy with Lemma \ref{lem:garland-dec} and the comments after it, 
 \begin{nclaim} \label{claim:conv-dec} Let $J \subset I_o.$ Every $H \in \mf{h}_o$ can be written uniquely as $H = H(J) + H_J$ for $H(J) \in \mf{h}_o(\es, J)= \Span \{ \av_i , i \in J \}$ and $H_J$ is such that $\la a_i, H_J \ra =0$ for $i \in J$. Moreover, if we write $H = \sum_{i \in I} \, c_i \omv_i$ with $c_i \in \R$ then $H(J) = \sum_{i \in J} c_i \omv^J_i.$ \end{nclaim}

Now using the claim, we obtain the desired reformulation for convergence:   $H_n \in \mf{h}_o$ converges to an element $H_{\infty} \in \mf{h}_o(\es, J)$ if when we decompose $H_n= H_{n, J} + H_n(J)$, then \begin{itemize} \item $H_n(J) \rr H_{\infty}$ in the metric topology of the element of $\mf{h}_o(\es, J)$; and \item $\la a_i, H_{n, J} \ra \rr - \infty$ for each $i \notin J$. \end{itemize} 

The first point is where we have used the claim: the limit of $H_n(J)$ in the metric topology is exactly what was required from \eqref{H:inf-def}. As for the second condition: formerly we required that $c_{n, i} \rr - \infty$ for all $i \notin J.$ But since we also required that $c_{n, j}$ are convergent for each $j \in J$, we actually have that $\la a_i, H_{n}(J) \ra$ will again be convergent even for each $i \in I.$ Thus  
\be{}  \la a_i, H_{n, J} \ra = \la a_i, H_n - H_n(J) \ra = c_{n, i} - \la a_i, H_{n}(J) \ra  \rr - \infty  \mbox{ if and only if } c_{n, i} \rr -\infty. \ee  
 
A similar approach can be taken to describe the convergence of $H_n \in \mf{h}_o(\es, J)$ to an element in $H_{\infty} \in \mf{h}_o(\es, K)$ with $K \subset J$ using the (analogues of) the relative decompositions of \S\ref{subsub:relative-variants} and \S \ref{subsub:compatability}. We leave this to the reader.

 \tpoint{Affine corners}  \label{subsub:affine-corners} Recall the basis $\{ \edv_i, \cpsi \}_{i \in I}$ of $(\mf{h}^e)^*$ from \S \ref{subsub:coweights} and write $H \in \mf{h}^e$ as \be{} \label{H:ep} H = c_1 \edv_1 + \cdots c_{\ell+1} \edv_{\ell+1} + f \cpsi, \mbox{ with all } c_i \mbox{ and } f \in \R. \ee 
 Then, setting $\hA^e:= (\R_{>0})^{\ell+2}$, we  define \be{} \begin{array}{lcr} \label{affine-exp} \exp: \hfh^e \rr \hA^e &\mbox{ by requiring} & H \mapsto (e^{c_1}, \ldots, e^{c_{\ell+1}}, e^f) \end{array}. \ee 
 We have singled out the set the Tits cone  $\hfh^+:= \{ H \in \hfh^e \mid \la H, \delta \ra < 0 \}$ in \eqref{ti:exp}. Under $\exp$, it maps (bijectively) into the set  \be{} \label{A+:cond} \hA^+:= \{ (s_1, \ldots, s_{\ell+1}, t) \mid t > 0, s_i >0 \mbox{ for } i \in I \mbox{ and }  s_1^{d_1} \cdots s_{\ell+1}^{d_{\ell+1}} <1   \}. \ee Denote by $\exp$ the restriction to the open set $\hfh^+ \subset \hfh^e$ and note that $\exp$ is a homeomorphism onto its image, where $(\R_{>0})^{\ell+2}$ is equipped with the metric topology.  
 
 Unlike the finite-dimensional case, we do not specify the corner $\cor(\hA^+)$ just by taking a closure in $\R_{\geq0}^{\ell+2}$ since we only allow each of the $s_i$ coordinates to approach $0$, but do not want to impose any conditions on the $t$-coordinate (we will make this more precise below). So, we take as our starting point the description \eqref{finite-corner} and introduce the  \textit{set}   \be{} \label{cor(A+):dec} \cor(\hA^+) := \hA^+ \, \bigcup \, \ \ \sqcup_{J \subsetneq I} \ \ \hA^+(\es, J) \times \hat{o}_J \ee where $\hA^+(\es, J) \times \hat{o}_J $ consist of the points in $(s_i) \in (\R_{\geq0})^{\ell} $ where $s_j=0$ for $j \in J$; we keep the notation $\hat{o}_J$ just as a placeholder for the `missing' central direction (the boundary pieces we add always have this extra `drop' in dimension as we omit the central direction in each of them).

To describe $\cor(\hA^+)$ at the Lie algebra level, start by observing $\hA^+(\es, J) \times \hat{o}_J$ is homeomorphic to $\hfh^+(\es, J)$ under the map $\exp(J)$ that sends $H(J)  \in \hfh^+(\es, J), $ written as $H(J) = \sum_{i \in J} c_j \omega^J_j,$ to  \be{} \label{ajcoord} \exp(H(J) ) \times \hat{o}_J= ( s_i)_{i \in I} \times \hat{o}_J \mbox{ with } s_i = e^{c_i} \mbox { if } i \in J, s_i=0 \mbox{ otherwise}. \ee If we give $A^+(\es, J) \times \hat{o}_J$ the subspace topology (inherited from $\R_{\geq0}^{\ell+1}$), the map $\exp_J$ is clearly a homeomorphism, and if we now define 
 \be{} \label{corh} \cor(\mf{h}^+) =  \mf{h}^+ \bigcup \, \sqcup_{J \subsetneq I} \mf{h}^+(J),  \ee then we can patch $\exp$ together with all of the $\exp_J$ to get a bijection \be{} \cor(\exp): \cor(\mf{h}^+) \rr \cor(A^+) \ee which is a homeomorphism on each of the pieces indexed by $J$.  The next task is to describe the topology on $\cor(\hA^+)$ or equivalently $\cor(\hfh^+)$.
 
\tpoint{The topology on $\cor(\hfh^+)$} Recall that we impose the metric topology on $\hfh^+$ and $\hfh^+(\es, J)$ for each $J \subsetneq I$, so we just need to explain how convergence works. We follow the Moore--Smith approach (see \S \ref{}) and introduce the following convergence class $\mathscr{C}$:
\begin{itemize}
	\item[$\mathscr{C}_{int}$]: (\textit{convergence from the interior}) Suppose we are given a sequence of elements $\{ X_n - r_n \dd \mid  X_n \in \hfh, r_n > 0 \}$ together with $Y_{ \infty} \in \hfh(\es, J).$ Then we say that $(X_n - r_n \dd, Y_{\infty})$ lies in our convergence class if when we decompose $X_n = Y_n + Z_n$ as in \eqref{hr:dec} with $Y_n \in \hfh(\es, J)$ and $Z_n \in (\hfh^{+})_J$ we have \begin{itemize}
		\item $Y_n \rr Y_{\infty}$ in the metric toplogy of $\hfh(\es, J)$; and 
		\item $a_i(Z_n - r_n \dd) \rr - \infty$ for each $i \notin J$. \end{itemize}
	\item[$\mathscr{C}_{bd}$]: (\textit{convergence along the boundary}) Suppose $I \supsetneq K \supset J$. We say that $\{ Y_n(K) , Y_{\infty}(J) \} $  with $Y_n(K) \in \hfh(\es, K)$ and $Y_{\infty}(J) \in \hfh(\es, J)$ lies in our convergence class if with respect to the decomposition \eqref{eq:relative-lie-algebra} $Y_n(K) = W_n + Q_n$ with $W_n \in \hfh(K)_{\leftidx^{\ast}J }$ and $Q_n \in \hfh(\es, J)$ we have \begin{itemize}
		\item   $Q_n \rr Y_{\infty}(J)$ in the metric topology of $\hfh(\es, J)$; and
		\item  $a_j(W_n) \rr - \infty$ for $j \in K \setminus J$. \end{itemize}

\end{itemize}

 One needs to check that the above defines a well-defined convergence class, and hence defines a topology on $\cor(\hfh^+).$ This is a somewhat lengthy verification and will be carried out in \S \ref{sub:verification-conv-class}. 
 
 \tpoint{An alternate take on the topology of $\cor(\hfh^+)$} \label{subsub:top-h-alternate} Alternatively, one can introduce the topology on $\cor(\hfh^+)$ using the coweight basis coordinates as follows (see the remarks after Lemma \ref{lem:garland-dec})
 
 \newcommand{\cclass}{\mathscr{C}}
 
 \begin{nlem} \label{lem:two-points-of-view-convergence} 
 	Then the convergence class $\mathscr{C}$ is equivalent to the following one:

 			\begin{itemize}
 				\item[$\cclass'_{int}$:] Let $H_n \in \hfh^+$ and $H_{\infty} \in \hfh(\es, J)$ for $J \subsetneq I$. Suppose we write as in \eqref{H:ep} \be{} \begin{array}{lcr} H_n = \sum_{i \in I} c_{n,i} \edv_i + f_n \cc & \mbox{ and } & H_{\infty} = \sum_{j \in J} c_{\infty, j} \omv^J_j \end{array}. \ee   Then  $H_n \rr H_{\infty}$ lies in the convergence class if $c_{n, i} \rr - \infty$ for $i \notin J$ and $c_{n, i} \rr c_{\infty, i}$ for $i \in J.$ 
 				\item[$\cclass'_{bd}$:] Suppose $Y_n(K) \in \hfh(\es, K)$ and $Y_{\infty}(J) \in \hfh(\es, J)$ with $K \supset J \subsetneq I$. Then $Y_{n}(\es, K) \rr Y_{\infty}(J)$ means that when we write $Y_{n}(\es, K) = \sum_{k \in K} c_{n, k} \omv_k^K$ and $Y_{\infty}(J) = \sum_{j \in J} c_{\infty, j} \omv_j^J$ then we must have $c_{n, k} \rr - \infty$ for $k \notin J$ and $c_{n, j} \rr c_{\infty, j}$ for $j \in J$. 			
 			\end{itemize}	

 \end{nlem} 
  
 \begin{proof} Suppose that $H_n \rr H_{\infty}$ lies in $\cclass_{int}$. Then by definition, when we decompose $H_n = H_{n, J}+ H_n(J)$ as in Lemma \ref{lem:garland-dec}, then have $\la a_i, H_{n, J} \ra \rr - \infty$ for $i \notin J$ and $H_n(J) \rr H_{\infty}$. Suppose we write $H_{\infty}$ as $\sum_{j \in J} c_{\infty} \omv_j^J$. Then using \eqref{H(J):p-c} applied to $H_{n, J}$ we obtain $H_{n, J} = \sum_{j \in J} c_{n, j} \omv_j^J$. Hence for $H_{n, J} \rr H_{\infty}$ in the metric topology, we must have $c_{n, j} \rr c_{\infty, j}$ which is one of the requirements to be in $\cclass'_{int}$.  As for the other requirement to be in $\cclass'_{int}$, we have noted that $\la a_i, H_{n, J} \ra \rr -\infty$. But we can compute this also as \be{} \la a_i, H_{n, J} \ra = \la a_i, H_n - H_n(J) \ra = c_{n, i} - \la a_i, H_n(J) \ra \mbox{ for } i \notin J. \ee Now since $H_n(J) \rr H_{\infty}$ in the metric topology and since $a_i$ is a linear functional on $\hfh^+(\es, J)$ we must also have $\la a_i, H_n (J) \ra \rr \la a_i, H_{\infty}\ra$, i.e. to a finite quanity. Hence for  $\la a_i, H_{n, J} \ra \rr -\infty$ we need to have $c_{n, i} \rr - \infty$ for $ i \notin J$. The converse, namely that if $H_n \rr H_{\infty}$ lies in $\cclass'_{int}$, then it lies in $\cclass_{int}$ follows from the same argument. 
 	
 	We leave the equivalence of $\cclass_{bd}$ and $\cclass'_{bd}$ to the reader.

 \end{proof} 
 
 We can also phrase this at the level of $\hA^+$ as well, and do so, we equip $\hA^+ \subset (\R_{>0})^{\ell+1} \times \R_{>0}$ with coordinates $(s_i; t)_{i \in I}$ as in \eqref{A+:cond}, and $\hA^+(\es, J) \cong (\R_{>0})^{\ell}$ with coordinates $(y_j)_{j \in J}$ as in \eqref{ajcoord}.
 
 \begin{nprop} \label{prop:top-on-cor(A)} The topology on $\cor(\hA^+)$ is described in terms of the following convergence class:
 	\begin{itemize}
 		\item A sequence of elements $h_n \in \hA^+$ with coordinate $(s_{n, i}; t)$ converges to $h_{\infty} \in \hA^+(\es, J)$ with coordinates $(s_{\infty, j})$ if $s_{n, i} \rr 0$ for $i \notin J$ and $s_{n, j} \rr s_{\infty, j}$ for $j \in J$.
 		\item A sequence of elements $h_n(K) \in \hA^+(\es, K) \cong (\R_{>0})^{|K|}$ with coordinate $(s_{n, k})$ converges to $h_{\infty}(J) \in \hA^+(\es, J) \cong (\R_{>0})^{|J|}$ with coordinates $(s_{\infty, j})$ if $s_{n, k} \rr 0$ for $k \in K \setminus J$ and $s_{n, j} \rr s_{\infty, j}$ for $j \in J$.

 	\end{itemize}

 \end{nprop}

 \tpoint{Example: $\widehat{sl}_2$} Consider the setup of Example \ref{example:sl2-epsilon}. The condition of $H \in \mf{h}^+$ is that $c_1 + c_2 < 0$. Under the map $\mf{h}^+ \rr (\R_{>0})^3$ which sends $(c_1, c_2, f) \mapsto (s_1, s_2, t)$, this carves out $A^+$ as the set $s_1 s_2 < 1$ with $s_1, s_2, t > 0$. Now there are three possibilites which we allow in  $\cor(\hA^+)$: \begin{enumerate}
 	\item $s_1 \rr 0$ but $s_2$ (and $t$) remain bounded, \textit{i.e.} $c_1 \rr - \infty$, $c_2, f$ remain bounded; 
 	\item $s_2 \rr 0$ but $s_1$ (and $t$) remain bounded, \textit{i.e.} $c_2 \rr - \infty$ and $c_2, f$ remain bounded
 	\item    $s_1, s_2 \rr 0$ (and  $t$) remain bounded, \textit{i.e.} $c_1, c_2 \rr - \infty$ but $f$ remains bounded  \end{enumerate}
 	
 One can check these are just the three ways in which we can have a sequence from $\hfh^+$ converge to $\hfh^+(\es, \{ a_1 \})$, 	$\hfh^+(\es, \{ a_2 \})$, and $\hfh^+(\es, \{ a_1, a_2 \})$ respectively.

\tpoint{Variant: $\cor(\hfh_J^+)$} For each $J \subsetneq I$, we can also consider the subset of $\hfh^e$ \be{} \hfh_J^+:= \{ H - r \dd \in \hfh^+ \mid H \in \hfh, r > 0, a_i(H - r \dd )=0 \text{ for } i \in J \}. \ee If $J= \es$, this is the set $\hfh^+$ we have previously considered.  \textit{Warning:} the above is a subset of $\hfh^e$ unlike the similarly named set $(\hfh^+)_J$ which is contained in $\hfh$. We can follow the same construction as above to attach a corner to $\hfh^+_J$, or alternatively we can observe

 \begin{nprop} \label{prop:closure-hJ} For $J \subsetneq I$, the interior of the closure of $\hfh^+_J \subset \hfh^+$ in $\cor(\hfh^+)$ is equal to   
 \be{} \label{pc:AJ}
 \cor(\hfh^+_J):= \hfh_{J}^+ \ \ \cup \ \  \bigsqcup_{J \subset K \subsetneq I} \hfh(J,K) =  \hfh_{J}^+ \ \ \cup \ \  \bigsqcup_{J \subset K \subsetneq I} \hfh(K)_{\sta{J}},
 \ee equipped with the subspace topology inherited from $\cor(\hfh^+)$.  \end{nprop} 
 \begin{proof} It suffices to observe that if a sequence $\{X_n\} \subset \hfh^+_J$ converges to $Z_{\infty} \in \hfh(\emptyset,K)$, then by definition of $\hfh^+_J$, we must have  $a_i(X_n) = 0 $ for $i \in J$. On the other hand, from the description in Lemma \ref{lem:two-points-of-view-convergence}, we must have $a_i(X_n) \rr -\infty$ for $i \notin K.$ Hence we must have
 	$(I \setminus K) \cap J  = \emptyset$, and so $ K \supset J$. We also need to verify that every point in the set \eqref{pc:AJ} arises in the closure; we leave this as an exercise again based on the description of Lemma \ref{lem:two-points-of-view-convergence}.
 	\end{proof}

\subsection{The corner $\cor(\hfh^+)$: further topological properties} \label{sub:top-corner} 

Having introduce the topological space $\cor(\hfh^+)$, or its homeomorphic image, $\cor(\hA^+)$, in the previous section, we now proceed to describe some of its finer  features.

 \tpoint{A basis for the topology of $\cor(\hfh^+)$} The bais $\{ \edv_i; \cc \}$ of $\hfh^e$ gives us a topological embedding 
\be{} \label{ev}\bold{ev}: \hfh^+ &\longrightarrow& \R^{\ell+1} \times \R , \\ H = f \cc + \sum_i c_i \edv_i &\mapsto& (a_1(H),a_2(H),...,a_{l+1}(H); \lambda_{l+1}(H))=(c_1, \ldots, c_{\ell+1}; f ). \ee
For an open set $U \subset \R^{l+2}$, $\bold{ev}^{-1}(U)$ is open in $\hfh^+$ and hence also in $\cor(\hfh^+)$. These will be open sets in the interior $\hfh^+$, and we next need to describe open sets that meet the boundary. 

\newcommand{\op}{\mathbcal{B}}
\newcommand{\cu}{\mathbcal{U}}

So, let $x(J) \in \hfh(\es, J)$, $J \subsetneq I$. As explained above, the map $p_J: \hfh^+ \rr \hfh^+(\es, J)$ which sends $\sum_i c_i \edv_i + f \cc$ to $\sum_{j \in J} c_i \omv_J^i$ is continuous in the topology of $\cor(\hfh^+)$. For $\epsilon >0$ and $T> 0$, define \be{} \begin{array}{lcr} \op(x(J); \epsilon, T ) := p_J^{-1} \left( \op(x(J), \epsilon) \right) \cap \hfh_{J, T}^+  & \mbox{ where } & \hfh_{J, T}^+ = \{ H \in \hfh^+ \mid \la a_i, H \ra < -T \mbox{ for } i \in J \}. \end{array} \ee and $\op(x(J), \epsilon)$ denotes an open set in the metric topology of $\hfh^+(\es, J)$.  This is an open set of $\hfh^+$.

 Similarly, if $K \supset J$, the map $p^K_J: \hfh(\es, K) \rr \hfh(\es, J)$ which sends $H(K) = \sum_{k \in K} c_k \omv_k^K$ to $\sum_{j \in J} c_j \omv_j^J$ is continuous, and so the set \be{} \label{open:K}  \op(K)(x(J); \epsilon, T):= (p^K_J)^{-1}(\op(x(J), \epsilon)) \cap \hfh(J, K)_T \ee  where $\hfh(J, K)_T= \{ H(K) \in \hfh(\es, K) \mid \la a_i, H(K) \ra < -T \mbox{ for } i \in K \setminus J \} $ is open in $\hfh(\es, K)$ (not in the corner $\cor(\hfh^+)$ though).   Putting these together, we define  \be{} \cu(x(J); \epsilon, T) := \op(x(J); \epsilon, T ) \cup \, \sqcup_{K \supset J} \op(K)(x(J); \epsilon, T). \ee 
 It can be seen to be an open set in $\corn{\hfh^+}$ containing $x(J)$, and in fact we have the following.

 \begin{nprop} \label{basis bordh}
     \begin{enumerate}
     	\item The space $\cor(\hfh^+)$ is Hausdorff. 
     	\item Let $\mathcal{B}$ denote any collection of $\cor(\hfh^+)$ containing 
     		\begin{enumerate}
         		\item $\bold{ev}^{-1}(U)$, where $U \subset \R^{\ell+1} \times \R$ ranges over a basis of the latter space;
         		\item $\cu(x(J); \epsilon, T)$ with $x(J) \in \hfh(\es, J)$ for $J \subsetneq I$ and $\epsilon >0, T >0$.
     		\end{enumerate}
     			Then $\mathcal{B}$ is a basis for the topology of $\cor(\hfh^+)$. Moreover, $\cor(\hfh^+)$ is first countable, i.e. each point has a countable basis. 
     \end{enumerate}
 \end{nprop}
 \begin{proof}
 		Since the topology is defined by a convergence class, it is enough by Lemma \ref{lem:top-from-seq} to show that a convergent sequence in $\cor(\hfh^+)$ has a unique limit. This is however clear from the description of the convergence class given in Lemma \ref{lem:two-points-of-view-convergence}. A similar argument holds for convergent sequences along the boundary.
 		
 	 It suffices to prove the following: let $x(J) \in \hfh(\es, J)$ and $H_n \rr x(J)$ lie in our convergence class. Then there exists $\epsilon, T$ so that $\cu(x(J); \epsilon, T)$ contains infinitely many $H_n$. Without loss of generality we may assume that all $H_n \in \hfh(\es, K)$ and so we are dealing with elements of type $\cclass_{bd}$ in the terminology of Lemma \ref{lem:two-points-of-view-convergence}.  Now, pick $\epsilon >0$ so that $\op(J)(x(J), \epsilon)$ contains almost all $H_n(J)$ and $T >0$ so that $\la a_i, H_{n, J} \ra< -T$ for almost all $n$. Then $\op(K)(x; \epsilon, T)$ will contain almost all $H_n$ by using the expression \eqref{open:K}. 
 	
 	To obtain a countable basis around a point $x \in \cor(\hfh^+)$,we consider two cases. If $x \in \hfh^+$, this is possible since $\hfh^+$, being a subset of the first countable space $(\R_{>0})^{\ell} \times \R$, is again first countable. If $x(J) \in \hfh(\es, J)$ for some $J \subsetneq I$, then we take $\cu(x(J); 1/{n_1}; n_2)$ as $n_1, n_2$ vary over the integers.  	\end{proof}
  
\newcommand{\ev}{\bold{ev}}

\renewcommand{\e}{\mathbf{e}}

 \tpoint{Homotopy type of $\cor(\hfh^+)$} In this section, we will argue that $\cor(\hfh^+)$ is contractible by showing the identity map $i: \cor(\hfh^+) \rr \cor(\hfh^+)$ is homotopic to the constant map $\e_{\es}$ sending all of $\cor(\hfh^+)$ to the unique element $o_\es \in \cor(\hfh^+)$. In fact, we will show this at the group level for $\cor(\hA^+).$ Denote again by $o_\es \in \cor(\hA^+)$ the unique elemement in $\hA^+(\es, \es)$.
 
 \begin{nprop} \label{prop:contractible-corner-A} The space $\cor(\hA^+)$ is contractible. \end{nprop}
 \begin{proof} Denote as before $i: \cor(\hA^+) \rr \cor(\hA^+)$ the identity map. Recall also the coordinates we have put on $\cor(\hA^+) = \hA^+ \sqcup_{J \subsetneq I} \hA^+(\es, J)$ in \S \ref{subsub:affine-corners}, namely $(s_i; z)_{i \in I}$ on $\hA^+$ and $(y_j)_{j \in J}$ for $\hA^+(\es, J)$.  With this terminology, we will argue that the following map $H: \corn{\hA^+} \times [0, 1] \rr \cor(\hA^+)$ is a homotopy between $i$ and $\e_{\es}$:
 	 \be{} \label{homotopy-H} 
 	H(h,t)=
 	\begin{cases}
 		o_\es  & \text{if } h \in \cor(\hfh^+), t=0 \\
 		(s_1 e^{ \ln t}, \ldots, s_{\ell+1} e^{\ln t}; z) & \text{if } h=(s_1, \ldots, s_{\ell+1}; z) \in \hA^+, 0< t \leq 1\\
 		(y_1 e^{\ln t}, \ldots, y_{|J|} e^{\ln t}) & \text{if } h=(y_1, \ldots, y_{|J|}) \in \hA^+(J),0< t \leq 1, \emptyset \neq J \subsetneq I \\
 		o_\emptyset \, &  \text{if } h  = \es ,0< t\leq 1
 	\end{cases}
 \ee
As it is clear that $H(h, 0)= \e_{\es}$ and $H(h, 1)= i$, it remains to verify that $H$ is continuous.  Since $\cor(\hA^+)$ is first countable, it suffices to check continuity sequentially (see  \cite[Theorem 1.6.14]{Eng89}), i.e. given $(h_n, t_n) \rr (h_{\infty}, t_{\infty})$ in $\cor(A^+) \times [0, 1],$ we want to show that $H(h_n, t_n)  \rr H(h_{\infty}, t_{\infty})$. 

We have several cases to consider. First assume that $h_n$ and $h_{\infty}$ all lie in $\hA^+$. If $t_{\infty} = 0$ and $t_n  = 0 $ for $n \geq n_0$ then we are done as we defined $H(h,0) = \es$ for  $h \in \hA^+$. So, assume $t_n \neq 0$ for $n \gg 0$. Then as we can see by applying formula \eqref{homotopy-H}, the first $\ell+1$-components of $H(h_n, t_n)$ must go to zero. Indeed, the $h_n$ are bounded (they converge to $h_{\infty} \in \hA^+$), and so $H(h_n, t_n) \rr \es$ by definition of the topology on $\cor(\hA^+)$. Next, suppose that $t_{\infty} \neq 0$. Without loss of generality, we can also assume that $t_n \neq 0$ for all $n$ sufficiently large. Writing $h_n:= (a_{1, n}, \ldots, a_{\ell+1,n}; z_n)$ we have assumed $h_{\infty}  = (a_{i, \infty}; z_{\infty})$ where $a_{i, n} \rr a_{i, \infty}$ and similarly $z_n \rr z_{\infty}$. Hence $H(h_n, t_n) =( a_{i, \infty} e^{ \ln t_n}; z_n)$ converges to $(a_{i, \infty} e^{\ln t_{\infty}}; z_{\infty})$ as desired.

Suppose now that $h_n \in \hA^+$ and $h_{\infty} \in \hA^+(J)$ for some $\es \subsetneq J \subsetneq I$. Writing $h_n=(s_{n,i})_{i \in I}$ and $h_{\infty} = (y_{\infty, j})$ we have that $s_{n, i} \rr 0$ for $i \notin J$ and $s_{n, j} \rr y_{\infty, j}$ for $j \in J$. Now we again have two cases. First, if $t_{\infty}=0$, then we may argue as in the previous paragraph. So we assume that $t_{\infty} \in (0, 1]$, and also, without loss of genrality we assume that all $t_n \geq 0$. Then to compute $H(h_n, t_n)$ we use the second formula in \eqref{homotopy-H} and observe both that $s_{n, i} e^{\ln t_n}$ for $i \notin J$ will converge to zero, since the $e^{\ln t_n}$ stay bounded, and that similarly $s_{n, j} e^{\ln t_n} \rr y_{\infty, j} e^{\ln t_{\infty}}$ for $j \in J$. Hence we see the limit converges to $H(h_{\infty}, t_{\infty})$.

All the other cases are handled similarly, and we leave this this to the reader. \end{proof}

\tpoint{Homotopy to the interior} Let $F_{- \rhov}: \corn{\hfh^+} \rr \hfh^+$ be the map which `pushes the boundary' inside by $- \rhov:= \sum_{i \in I} \lv_i$. More explicitly, we define 
	\begin{equation} \label{define:F-into-interior}
		F_{- \rhov}(X)=
		\begin{cases}
			\displaystyle \sum_{i \in I} \left(e^{-|\la X, a_i\ra|}-1\right) \lv_i & \text{if } X \in \hfh^+\\[1em]
			\displaystyle \sum_{i \in J} \left(e^{-|\la X, a_i \ra|}-1\right) \lv_i -\sum_{i \notin J}\lv_i & \text{if } X \in \hfh(J),\ J \subsetneq I
		\end{cases}
	\end{equation}
Notice that $F_{- \rhov}(\e_{\es})= - \rhov$ and also that $F_{- \rhov}$ kills the central direction.  We claim the map $F_{- \rhov}$ is continuous, and can check this using the sequential criterion. 
	
	\begin{itemize}
		\item If $X_n \longrightarrow X_0 $ in $ \hfh^+$, then $\la X_n, a_i \ra  \longrightarrow \la X_0, a_i\ra$ for all $i \in I$, and so $F_{- \rhov}(X_n)$ converges to $F_{- \rhov}(X_0).$
		 
		\item If $X_n \in \hfh^+$ converges to $X_0 \in \hfh(J)$ , then  $\la X_n, a_i \ra \rr -\infty$ for $i \notin J$ and $\la X_n, a_i \ra \rr  \la X_0, a_i \ra $ for $i \in J$. Hence 
		\begin{align*}
			F_{- \rhov}(X_n) &= \sum_{i \in I} \left(e^{-|\la X_n, a_i \ra|}-1\right) \lv_i 
			 =  \sum_{i\in J} \left(e^{-|\la X_n, a_i \ra|}-1\right) \lv_i + \sum_{i\notin J} \left(e^{-|\la X_n, a_i \ra|}-1\right) \lv_i \\
			& \longrightarrow  \sum_{i \in J} \left(e^{-|\la X_0, a_i \ra|}-1\right) \lv_i -\sum_{i \notin J}\lv_i = F_{- \rhov}(x_0).
		\end{align*} The case of convergence from $\hfh(J)$ to $\hfh(K)$ is handled similarly. 

	\end{itemize}

	\begin{nprop} \label{prop:cor-homotopy-interior} Let  $\iota: \hfh^+ \hookrightarrow \corn{\hfh^+}$ be the natural inclusion and $F_{- \rhov}: \corn{\hfh^+} \rr \hfh^+$ the map above. Then $\iota \circ F_{- \rhov}$ and $F_{- \rhov} \circ \iota$ are homotopic to the identity, i.e. $\corn{\hfh^+}$ is homotopic to its interior $\hfh^+$. 	\end{nprop}

\tpoint{Compact subsets of $\cor(\hfh^+)$} \label{subsub:compact-subsets-h} Since the topology of $\cor(\hfh^+)$ is first  countable by proposition \ref{basis bordh}, we know that a closed subset $C \subset \cor(\hfh^+)$ is compact if and only if every sequence in $C$ has a convergent subsequence converging to a point in $C$. Note that $\hfh^+ \subset \hfh$ is itself not closed in the metric topology of the latter, however $ \{ H \in \hfh^+ \mid - \la H, \delta \ra \geq r_0 \} = \cup_{r \geq r_0} \hfh^r$ is closed.  Bearing this in mind, we introduce for $t, r_0 \in \R_{>0}$,  the subsets \be{} \label{hfh:t} \hfh^{+}_{(-\infty,t]}(r_0) &=& \{H \in \hfh^+ \mid \la H, a_i \ra  \leq t , \, i \in I  \} \cap \{ H \in \hfh^+ \mid - \la H, \delta \ra \geq r_0 \}  \\  \label{hfh:J-t} \hfh(J)_{(-\infty,t]} &=& \{H \in \hfh(J) \mid \la H, a_i \ra  \leq t , \, i \in J  \}   \ee as well as following closed subset of  $\cor(\hfh^+)$:
\be{} \corn{\hfh^+; r_0, t} = \hfh^+_{(-\infty,t]}(r_0) \sqcup\bigsqcup_{J \subsetneq I}\hfh(J)_{(-\infty,t]}. \ee The conditions imposed in \eqref{hfh:t}-\eqref{hfh:J-t} do not constraint the central value of the component in $\hfh^+$. Thus, to obtain a compact subset of $\corn{\hfh^+}$ we introduce for numbers $M_0, r_0, t >0$  
\be{} \label{C:m-t-head} \mc{K}(\hfh^+; M_0,r_0, t) = \{H \in \hfh^+_{(-\infty,t]} \mid |\la \Lambda_{\ell+1} , H \ra| \leq M_0\} \cap \{ H \in \hfh^+ \mid - \la H, \delta \ra \geq r_0 \} \ee as well as its extension into the corner $\cor(\hfh^+)$  \be{} \label{C:m-r-t-corner} \corn{\hfh^+; M_0,r_0, t}:= \mc{K} (\hfh; M_0,r_0, t) \,  \sqcup\bigsqcup_{J \subsetneq I}\hfh(J)_{(-\infty,t]} \ee

\begin{nlem} \label{lemma:compact bordh} 
   For every $M_0, r_0, t >0$, the subspace of $ \corn{\hfh^+; M_0,r_0,  t} \subset \corn{\hfh^+}$ 
    is compact. 
\end{nlem}
\begin{proof} It suffices to argue that for any sequence $\{x_n\} \subset \corn{\hfh; M_0, r_0, t}$ has a limit.  By passing to a subsequence if necessary, it is enough to consider the following two cases,
  \begin{enumerate}
      \item The $x_n \in \corn{\hfh^+; M_0,r_0,t} \cap \hfh^+ = \mc{K}(\hfh^+; M_0,r_0, t)$. In this case, writing \be{} x_n = c_{n,1} \edv_1 + \cdots + c_{n, \ell+1} \edv_{\ell+1} + f_n \cc, \ee  our assumption that $x_n \in \mc{K}(\hfh^+; M_0,r_0, t)$ implies that all the $c_{n, i}$ are bounded above and the $f_{n, i}$ are bounded. So we have two possibilites to consider.

      \begin{itemize}
          \item If all the $c_{n, i}$ are also bounded below, then the $\{ x_n \}$ lie within a bounded subset of $\R^{\ell+2}$ and so have a limit (since $\mc{K}(\hfh^+; M_0, r_0,t)$ was closed as well.)
          
          \item Suppose $c_{n, i}$ are bounded for $i \in J$ but unbounded for $i \notin J$ for $J \subsetneq I$. By the definition of convergence given in \S \ref{subsub:top-h-alternate}
                we see that there exists a unique $x_{\infty} \in \hfh(J)_{(-\infty,t]}$ such that $x_n \rr x_{\infty}$.
      \end{itemize}
      \item  Let $x_n$ lie in the some boundary piece $\hfh^+(J)_{(-\infty, t]}$. Writing $x_n = \sum_{j \in J} c_{n, J} \omv^J_j$ our assumptions again imply that each of the $c_{n, J}$ are bounded above. Let $K \subset J$ be such that the sequence $c_{n, k}$ is bounded for $k \in K$ and unbounded for $i \in J \setminus K$. If $K = \es$, then the limit occurs in $\hfh^+(J)_{(-\infty, t]}$ itselt, and if $K \subsetneq J$, then one easily checks the limit occurs in $\hfh^+(K)_{(-\infty, t]}$

  \end{enumerate}  
\end{proof}

\begin{nrem} \label{remark:boundary-compact} Note that the proof shows that the boundary of $\corn{\hfh^+; M_0,r_0, t}$ is already compact. The imposition of a condition, depending on $M_0$, on the central direction is only to make its `head' compact, cf. Proposition \ref{prop:compact-boundary}. \end{nrem}

\subsection{Appendix: Verifying the convergence class conditions on $\cor(\hfh^+)$ } \label{sub:verification-conv-class} 

\begin{nprop} \label{prop:hplus-conv-class} The above defines a well-defined convergence class on $\cor(\hfh^+)$.  \end{nprop}

\begin{proof}  Properties $(1)-(3)$of (\ref{cclass}) are easy to check, so let us verify (4). Suppose given a sequence $X_{m,n} \in \hfh^+$ such that such that for every fixed $n$, the $\{X_{m,n}\}_{m \in \mathbb{N}}$ converges to $Y_n \in \hfh(\es, K)$ for some $K \subsetneq I$, and such that  $\{Y_n\}$ converges to a point $Y_{\infty} \in \hfh(\es,J)$ for some $I \supsetneq K \supset J.$ We would like to show that there exists a function $f(n)$ with $\lim\limits_{n \rr \infty}f(n) = \infty$ such that $(X_{n, f(n)},  Y_{\infty})$ lies in our convergence class.  To verify this, we need to first decompose $X_{m, n}$ as follows:  \be{} X_{m, n} = X_{m, n}(K) + X_{m, n, K} \mbox{ with } X_{m, n}(K) \in \hfh(\es, K), X_{m, n, K} \in \hfh^+_{K} \ee and then further write \be{} X_{m, n}(K) = X_{m, n}(J) + X_{m, n}(K)_{\sJ} \mbox{ where } X_{m, n}(J) \in \hfh(K, J), X_{m, n}(K)_{\sJ} \in \hfh(K)_{\sJ}. \ee Hence, the decomposition of $X_{m,n}$ in $\hfh^+ = \hfh(J) + (\hfh^+)_J$ is given by
	\be{} X_{m, n} = \underbrace{X_{m, n}(J)}_{\hfh(\es, J)} + \underbrace{X_{m, n}(K)_{\sJ} + X_{m, n, K}}_{(\hfh^+)_J}.\ee We shall also write $X_{m, n, J}:= X_{m, n}(K)_{\ast J} + X_{m, n, K}$.

	By our assumptions, we have \be{} \begin{array}{lcr} \lim_{m \rr \infty} X_{m, n}(K) = Y_n & \text{ and } & \lim_{n \rr \infty} \pr^K(J)(Y_n) = Y_{\infty} \end{array}, \ee where $\pr^K(J):\hfh(\es, K) \rr \hfh(\es, J)$, the natural projection, is continuous with  respect to the metric topology. Hence, from Lemma \ref{lem:relative-lie-algebra}, we deduce that \be{} Y_{\infty}  &=&  \lim_{n \rr \infty} \pr^K(J)(Y_n) = \lim_{n \rr \infty} \pr^K(J) \left( \lim_{m \rr \infty}  X_{m, n} \right) = \lim\limits_{n \rr \infty} \lim\limits_{m \rr \infty} \pr^K(J) \circ \pr(K) (X_{m, n}) \\ &=&  \lim_{n \rr \infty} \lim_{m \rr \infty} \pr(J)(X_{m, n}). \ee

	So to verify that $(X_{m, n}; Y_{\infty})$ with $m=f(n)$ (for $f$ an increasing function to be specified) lies in our convergence class, we just need to verify that for $i \notin J$ that \be{} \lim_{n \rr \infty} \lim_{m \rr \infty} a_i(X_{m,n, J}) = \lim_{n \rr \infty} \lim_{m \rr \infty} a_i(X_{m, n}(K)_{\ast J} + X_{m, n, K} )= - \infty, \ee where again $m:= f(n)$. We have two cases to consider.
	
	\begin{enumerate}
		\item For $i\in K \setminus J$, by assumption  $a_i(X_{m, n, K}) =0,$  so the same remains true when $m, n \rr \infty$ as well.  To control $X_{m,n}(K)_{\sJ}$, we first note that the projection induced from \eqref{eq:relative-lie-algebra}, say $\pr^K_J: \hfh(\es, K) \rr \hfh(K)_{\sJ}$, is continuous. So, since for each $n$,  $\lim\limits_{m \rr \infty} X_{m, n}(K) \rr Y_n$  in the metric topology of $\hfh(\es, K)$ we  have, setting $Y_{n, \sJ}:= \pr^K_J(Y_n)$, that 	
		\be{} \label{pr^K_J:1} \pr^K_J( X_{m, n}(K) ) = X_{m, n}(K)_{\sJ} \rr Y_{n, \sJ}. \ee As $a_i(Y_{n, \sJ}) \rr - \infty$ as $n \rr \infty$ and for $i \in K \setminus J$, it follows that the same holds for $X_{m, n}(K)_{\sJ}$ by the continuity of $a_i$. 	 In sum, we have shown  that as $m, n \rr \infty$, \be{} a_i( X_{m, n}(K)_{\sJ} + X_{m, n, K}) \rr - \infty \mbox{ for } i \in K \setminus J. \ee  Note that in this case, it did not particularly matter which choice we make for $f(n)$ (just so long as it increased $\infty$ as $n$ grows large)
		
		\item Suppose now that $i \notin K$. By assumption for each $n$, $a_i(X_{m, n, K}) \rr - \infty$ as $m \rr \infty$. We next need to analyze $ \lim\limits_{n, m \rr \infty} a_i (X_{m, n}(K)_{\sJ}).$ Since $X_{m, n}(K) \rr Y_{n, \ast J}$  and $a_i(Y_{n, \ast J}) \rr - \infty$ for $i \in K \setminus J$, we may conclude that for $n$ sufficiently large,  \be{} \label{Yn,J:cowts} Y_{n, \ast J} = \sum_{i \in K} c_{n, i} \omv^K_i  \mbox{ with }  c_{n, i} \leq 0 .\ee Indeed, by definition $a_i(Y_{n, \ast J}) = c_{n, i}$, and these values are either $0$ (for $i \in J$) or negative for large $n$. Rewriting $Y_{n, \ast J}$ in terms of the coroots $\av_k, k \in K$, we again obtain a linear combination with negative coefficients (since coweights are \textit{positive} linear combinations of coroots, see \S \ref{subsub:h(J)-coweights}). Now since $\la a_i, \av_j \ra \leq 0$ for $i \neq j$, it follows that if $i \notin K$, then $\la a_i, Y_{n, \ast J} \ra \rr + \infty$ as $n$ grows large. 
		
		In other words, for a fixed $i \notin K$ writing $p_{m,n} = a_i(X_{m, n, K})$ and $q_{m,n} =  a_i (X_{m, n}(K)_{\sJ})$ we are in the following scenario: $p_{m,n} \rr -\infty$ for each fixed $n$ as $m  \rr \infty$ and $q_{m,n} \rr b_n$ for each fixed $n$ as $m \rr \infty$ and $q_n \rr \infty$ as $n \rr \infty$. We observe then the following elementary result, which will conclude the argument.
		
		\newcommand{\n}{\mathbb{N}}
		
		\begin{nclaim} There exists a function $f: \n \longrightarrow \n$ such that $f(n) \rr \infty$ as $n \rr \infty$ and 
			\be{} \lim_{n \rr \infty} p_{f(n),n} +q_{f(n),n} =-\infty \ee \end{nclaim} 
		\begin{proof}
			For every $n \in \n$, we can find a $J_n \in \n$ such that $p_{m,n} < -q_n^2 $ if $m > J_n$. Similarly there exists a $K_n \in \n$ such that
			$|q_{m,n}- q_n| < \frac{1}{n}$ for $m > K_n$. So we get,
			\be{} p_{m,n}+q_{m,n} < -q_n^2 +q_n +\frac{1}{n} \mbox{ for } m > \max\{J_n,K_n\} \ee
			We can certainly assume that both $\{J_n\}_{n \in \n}, \{K_n\}_{n \in \n}$ are strictly increasing. Let $f(n) = \max\{J_n,K_n\} +1 $, then by the above inequality and the fact that $q_n \rr \infty$ as $n \rr \infty$ we get  $\displaystyle \lim_{n \rr \infty} p_{f(n),n} +q_{f(n),n} =-\infty$.  
			
	\end{proof}  \end{enumerate} \end{proof}

\section{Loop groups and their reduction theory } \label{sec:loop-groups}

The aim of this section is to introduce the loop groups which will be studied in this paper. We follow the method of Garland \cite{gar:ihes} which constructs the group together with its action on a natural highest weight module (much like the approch of Chevalley in finite type,  see \cite{stein}. We describe a generalization to loop groups of Arthur's notion of orthogonal families in \S \ref{sec:orthogonal-families} explain the reduction theory for loop groups following \cite{gar:ihes} in \S \ref{sec:reduction}.

\subsection{Loop groups: basic definitions} Let $\Phi$ be the Chevalley basis of $\hfg$ constructed in \S \ref{subsub:chev-basis}.	Fix $\lambda$ a normal element in $\Lambda_+$ and fix a representation $V^{\lambda}$ and $\zee$-form $V^{\lambda}_{\zee}$ corresponding to an admissible basis $\Upsilon$ as in \S\ref{subsub:hw-lambda-chev}.  In this section, as $\lambda$ is usually fixed, we generally drop it from our notation, with the exception that the weight lattice of $V$ will still be denoted by $\Xi_{\lambda}.$

\newcommand{\li}[1]{\leftidx^{{#1}}}
\newcommand{\winv}{w^{-1}}
\newcommand{\Qw}{\bQ_{\winv}}
\newcommand{\qw}{Q_{\winv}}
\newcommand{\Pw}{\bP_{w}}
\newcommand{\pw}{P_{w}}
\newcommand{\sqw}{\leftidx^*\qw}
\newcommand{\sQw}{\leftidx{^*}\Qw}
\newcommand{\sPw}{\leftidx{^*}\Pw}
\newcommand{\spw}{\leftidx^*\pw}

\newcommand{\au}{\aut}

\newcommand{\cbG}{\widehat{\mathbf{G}}}

\tpoint{On the elements $\chi_a(s), h_a(s), w_a(s)$} \label{subsub:elements-in-G} For $a \in \rts_{re}$ and $t \in \R$ , we set $\chi_a(t)$ to be the formal exponential $\exp(t \xi_a)$, which defines automorphism on $V$ by \cite[Lemma 7.12]{gar:ihes}. They satisfy \be{} \chi_a(s) \chi_a(t) = \chi_a(s+t) \mbox{ for } a \in \rts_{re}, s, t \in \R. \ee Defining   for  $a \in \rts_{re}$ and  $s \in \R^*$ the elements \be{w,h} \begin{array}{lcr} w_a(s) := \chi_a(s) \chi_{-a}(-s^{-1}) \chi_a(s) & \text{ and } & h_a(s) := w_a(s) w_a(1)^{-1} \end{array}, \ee we find that if $\mu \in \wts_{\lambda}$ , $v \in V_{\mu},$ and $s \in \k^*$ then from \cite[Lemma 11.2]{gar:ihes} \be{ha:act} h_a(s) . v = s^{ \la \av, \mu \ra} v. \ee  Let $\bN$ be the subgroup generated by $w_a(s), s \in \R^*, a \in \rts_{re}$, and let us now set \be{} \label{wi} \dw_i:= \dw_{a_i}:= w_{a_i}(1) \mbox{ for } i \in I. \ee These elements verify the braid relations attached to the gcm $\As$, so that for any element of the Weyl group  $w \in W:=W(\As)$ with reduced decomposition $w = s_{i_1} \ldots s_{i_k}$ with $i_j \in I$, we may unambiguously define $\dw:= \dw_{a_{i_1}} \cdots \dw_{a_{i_k}} \in \bN. $ Note (see \cite[Lemma 11.2 (i)]{gar:ihes}) that for $a \in \rts_{re}$ and $\mu \in \wts(V)$, if $v \in V_{\mu}$ then there exists $v' \in V_{s_a(\mu)}$ such that $w_a(s) v = s^{- \la \av, \mu \ra} v'$ for any $s \in \R^*$. 
 
\tpoint{The elements $\chi_{\alpha}(\bbsigma)$} \label{subsub:comp-pos} For $\alpha \in \rts_{o}$ and $\bbsigma=\sum_{i \geq i_0} \sigma_i t^i \in \R((t))$ the infinite product  \be{chi:com} \chi_{\alpha}(\bbsigma):= \Pi_{i \geq i_0} \chi_{\alpha + i \delta} (\sigma_i) \ee continues to define an element of $\Aut(V)$ using \cite[Lemma 7.16]{gar:ihes}. Note that if $\bbsigma(t)= st^n$ for $s \in \R$, then $\chi_\alpha(\bbsigma) = \chi_{\alpha + n \delta}(s)$. One again has the one-parameter property (see \cite[(7.20)]{gar:ihes})  \be{chi:add} \chi_{\alpha}(\bbsigma + \bbtau ) = \chi_{\alpha}(\bbsigma) \chi_{\alpha}(\bbtau) \mbox{ for } \bbsigma, \bbtau \in \R((t)), \alpha \in \rts_o \ee satisfying commutation relations as in the finite-dimensional case (see \cite[Lemma 9.7]{gar:ihes}). 

\tpoint{Loop rotation} \label{subsub:loop-rotation} For each $s \in \R^*$, define $\eta(s) \in \Aut(V)$ by the formula \be{eta:defs} \eta(s). v = s^{ \la \dd, \mu \ra } v \text{ for } v \in V_{\mu}. \ee One can then verify that in $\Aut(V)$ we have the identity \be{} \label{eta:conj-1} \eta(s) \chi_a(t) \eta(s)^{-1} = \chi(s^{ \la \dd, a \ra} t) \text{ for } a \in \rts_{re} \text{ and } s \in \R^*, t \in \R. \ee  From this, one can then deduce that \be{} \label{eta:conj-2} \eta(s) \chi_{\alpha}(\bbsigma(t)) \eta(s)^{-1} = \chi(\bbsigma(st)) \text{ for } \alpha \in \rts_{o} \text{ and } \bbsigma(t) \in \R((t)). \ee

\tpoint{Loop groups: definitions. } Finally, we have all the ingredients to define the loop group  \be{} \label{def:G-hat}  \hG:= \hG^{\lambda} = \la \chi_{\alpha}(\bbsigma),  \alpha \in \rts_o, \bbsigma \in \R((t)) \ra \subset \Aut_{\R}(V^{\lambda}). \ee We also set $\hG^e:= \hG^{\lambda, e}$ to be the group generated by $\hG$ together with the automorphisms $\eta(s), \, s \in \R^*$ defined in \eqref{eta:defs}.  Sometimes we wish to fix a value of $s= e^{-r} \in \R^*$, so $r > 0$, and consider  \be{} \begin{array}{lcr} \hG^r:= \{ g \eta(s) \in \hG^e \mid g \in \hG \} & \mbox{ and } &  \hG^+:= \sqcup_{r > 0 } \hG^r \end{array}.\ee One can show that $\hG^+$ is actually a semi-group. 

 \label{remark-on-rational-points} This constructions of $\hG$ is carried out over $\R$, but a similar construction works over a general field $\k$ by working in $V_{\k}:= \k \otimes_{\zee} V^{\lambda}_{\zee}$ in place of $V:= V_{\R}$. We shall later need to work with groups over $\Q$, in which case we shall write $\hg_{\Q}$ for the corresponding object.  Our convention will be that if $H \subset \hG$ denotes some subgroup, then $H_{\Q}:= H \cap \hG_{\Q}.$

 \tpoint{Relation to defining representation $V^{\lambda}$} The group $\hG$ is defined in terms of a representation $V^{\lambda}$ and when we wish to emphasize this dependence we write $\hG^{\lambda}$ and refer to its set of generators as $\chi_{\alpha}^{\lambda}(\bbsigma)$. One can compare the groups $\hG^{\lambda}$ for using the following result: 
 
 \begin{nthm} \cite[Thm. 15.9]{gar:ihes} \label{thm:change-of-weights} Let $\lambda_1, \lambda_2$ be two normal weights and assume that $\Xi_{\lambda_2} \subset \Xi_{\lambda_1}$. Then there exists a unique group homomorphism \be{} \pi(\lambda_1, \lambda_2): \hG^{\lambda_1} \rr \hG^{\lambda_2} \ee which sends $\chi_{\alpha}^{\lambda_1}(\bbsigma)$ to $\chi_{\alpha}^{\lambda_2}(\bbsigma)$ for $\alpha \in \rts_o$ and $\bbsigma \in \R((t))$. \end{nthm}  Combing this result with the remarks from \S \ref{subsub:V-lambda}, we find that if $\mu$ is any normal weight there exists a positive integer $m$ such that we have a map \be{} \label{Gmu:G} \hG^{\lambda} \rr  \hG^{m \, \mu}. \ee In particular, if $\lambda$ is such that $\Xi_{\lambda} = \Lambda$, the group $\hG^{\lambda}$ acts on $V^{\mu}$ or every normal $\mu$. We say that $\hg:= \hg^{\lambda}$ is \textit{simply-connected} in this case.

\tpoint{The torus $\hT^e$ and its connected component $\hA^e$} \label{s:AJ} Let $\hT$ (respectively $\hT^e$) be the subgroup of $\hG$ (resp. $\hG^e$) generated by $h_a(s)$ for $s \in \R^*$ and $a \in \rts_{re}$ (resp., together with $\eta(s), s \in \R^*$). From (\ref{ha:act}) and (\ref{eta:defs}), for $h \in \hT^e, \mu \in \wts(V)$, there exists $h^{\mu} \in \R^*$ such that \be{h:mu} h. v = h^{\mu} v \text{ for any } v \in V^{\lambda}_{\mu}. \ee Since the root lattice is contained in $\Xi_{\lambda}$, we may define $h^{\beta}$ for $\beta$ in this root lattice.  

\begin{nlem} (Structure of $\hT^e$) \begin{enumerate}
		\item (cf. \cite[Lemma 14.20]{gar:ihes}) For each $a \in \rts_{re}$, $h_a(s)$ is a multiplicative function of $s \in \R^*$ and $\hT$ is an abelian group generated by $h_i(s_i):= h_{a_i}(s_i)$ for $i \in I, s_i \in \R^*$; the group $\hT^e$ has an additional generator $\eta(s), s \in \R^*$.
		\item We have $\eta(s) \,  \prod_{i \in I} h_i(s_i) =1 \,$ if and only if $\prod_{i \in I} s_i^{ \la \av_i, \mu \ra} \, s^{ \la \dd, \mu \ra} =1$ for all $\mu \in \Xi_{\lambda}$.
	\end{enumerate}
\end{nlem}  

\noindent Note that if $h = \prod_{i \in I} h_{a_i}(s_i) \eta(s)$ is \textit{any} decomposition of $h \in \hT^e$, part (2) tells us that \be{} h^{\mu} =  \prod_{i \in I} s_i^{ \la \av_i, \mu \ra} \, s^{ \la \dd, \mu \ra}. \ee If $\Xi_{\lambda}= \Lambda$, the Lemma also tells us that every $h \in \hT^e$ has a unique decomposition \be{} \label{h:unique} h = \prod_{i=1}^{\ell+1} h_i(s_i) \, \eta(s) \mbox{ for } s_i \in \R^*, s \in \R^*. \ee In this case we use the natural map to define the standard Euclidean topology on $\hT^e$ rendering it homeomorphic to $(\R^*)^{\ell+2}$.  We are especially interested the subgroups \be{} \hA:= \la h \in \hT^e \mid h^{\mu} > 0 \mbox{ for all } \mu \in \Xi_{\lambda} \ra  =\la h_{a_i}(t) \mid t> 0 \ra .\ee It is the connected component of the identity of $\hT$, and one can also show that if $h \in \hA$, we may write it as $h = \prod_i h_{a_i}(s_i) \eta(s)$ with $s_i, s >0$.  We also define $\hA^e$ in a similar fashion, and for $0< s < 1$, say $s=e^{-r}$ with $r >0$, we shall write
 \be{} \begin{array}{lcr} \hA^s:= \hA^r:= \{ h \eta(s) \mid h \in \hA \} & \mbox{ and } & \hA^+:= \cup_{0< s< 1} \hA^s \end{array}. \ee One can check that the Lie algebra corresponding to $\hA^e$ is $\hfh^e$ and that the subset $\hA^s$ gets sent to $\hfh^r$ with $s=e^{-r}$ under this correspondence.

\tpoint{On the center $\hT_{\cc}$} \label{subsub:central-T} In analogy with the relation $\cc = \av_{\ell+1} + \sum_{i \in I_o} \dv_i \av_i$, we define the element \be{} \label{def:hc} h_{\cc}(s) := h_{\ell+1}(s) \prod_{i \in I_o} h_i(s)^{\dv_i} \mbox{ for } s \in \R^* \ee and note that this element is central in $\hG$. We also define the subgroups \be{} \begin{array}{lcr} \hT_{\cc}:= \{ h_{\cc}(s) \mid s \in R^* \}& \mbox{ and } & \hA_{\cc}:= \{ h_{\cc}(s) \mid s > 0 \} \end{array}. \ee

\newcommand{\heis}{\mathscr{H}}

\tpoint{The infinite Heisenberg subgroup $\mathscr{H}$} \label{subsub:BorelKM}For each $\alpha \in \rts_o$ and  $\bbsigma \in \R((t))^*,$ we introduce \be{wa:ha} \begin{array}{lcr} 
	w_{\alpha}(\bbsigma) = \chi_{\alpha}(\bbsigma) \chi_{-\alpha}(-\bbsigma^{-1}) \chi_{\alpha}(\bbsigma) &\text{ and } & 
	h_{\alpha}(\bbsigma) = w_{\alpha}(\bbsigma) w_{\alpha}(1)^{-1}. \end{array}\ee The elements $h_{\alpha}(\bbsigma)$ no longer act diagonally on $V$, \textit{cf.} \eqref{ha:act}, for a general $\bbsigma \in \R((t)).$ In fact, \be{} \mathscr{H}:= \la h_{\alpha}(\bbsigma), \alpha \in \rts_o, \bbsigma \in \R((t))^*, h_{\ell+1}(s), s \in \R^* \ra \ee is an infinite Heisenberg group, see \cite[\S 12]{gar:ihes}). On the other hand, if we let $(\R[[t]]^*)_1:= \{ \bbsigma \in \R((t))^* \mid \bbsigma \equiv 1 \, \mod\, \,  t \}$, then we do have an abelian subgroup pn $\heis$,  \be{} \label{heis:1} \heis_1:= \la h_{\alpha}(\bbsigma) \mid \bbsigma \in (\R[[t]]^*)_1, \alpha \in \rts_o \ra. \ee
	
\tpoint{The Borel subgroup of $\hG$}  \label{subsub:BorelKM}

 For each $a \in \rts_{re}$, we define the one parameter subgroup of $U_a:= \la \chi_a(s) \mid s \in \R \ra   \subset \hG.$  Let us also define $U$ to be the subgroup generated by  \be{} \begin{array}{lcr}   \chi_{\alpha}(\bbsigma), \, \,  \alpha \in \rts_{o, +}, \bbsigma \in \R[[t]], &  \chi_{-\alpha}(\bbsigma'), \, \,  \alpha \in \rts_{o,+} , \bbsigma' \in t\R[[t]], \mbox{ and } & \heis_1 \end{array}. \ee With respect to a coherently ordered basis of $V$ (see \S \ref{subsub:hw-lambda-chev}), the elements of $U$ act as ``unipotent upper triangular matrices,'' see \cite[(14.15)]{gar:ihes}. There is a natural homomorphism $U_a \rr U$ sending $\chi_a(s)$ for $s \in \R$ and $a= \alpha + n \delta$ to $\chi_{\alpha}(s t^n).$ Defining $\hB \subset \hG$ as the subgroup generated by $\hT$ and $\hU,$   \be{tB:semi} \hB = \hT \ltimes \hU. \ee To simplify notation, we shall often write, when confusion with the finite-dimensional case is unlikely, $ B:= \hB. $ We also write $\hB^e:= \hT^e \ltimes \hU$.

 \renewcommand{\hB}{B}
 
\tpoint{Iwahori-Matsumoto Coordinates on $\hU$}  \label{subsub:iwahori-matsumoto-coordinates} The group $\hU$ can be equipped with the structure of a pro-unipotent algebraic group (see \cite[\S4.4 and \S 6.1]{kumar}). The pro-structure can be defined, with respect to our coherently ordered basis $\Upsilon= \{ v_0, v_1, v_2 , \ldots \}.$ For a positive integer $N$, let $\Upsilon[N]$ consist of all the vectors of depth less than or equal to $N$, and let $V[N] \subset V$ be the $\R$-span of $\Upsilon[N]$. Note that $\hU$ preserves $V[N]$ (since acting by $\hU$ is depth lowering). Introducing the normal subgroup  \be{} \hU^{(N)}:= \{ u \in \hU \mid u \equiv \id \mbox{ on } V / V[N]  \}, \ee one can verify that the properties of a pro-unipotent group are satisfied for the maps \be{} \label{proj-system-U} \hU/\hU^{(1)}  \stackrel{p_1}{\longleftarrow}  \hU/\hU^{(2)} \longleftarrow \cdots \ee     In fact, what will also be useful for us is an explicit set of coordinates that Garland has introduced on $\hU$ (it is motivated by constructions of Iwahori-Matsumoto coordinates in the  $p$-adic groups \cite[Prop. 2.4]{iwahori:matsumoto}). To decribe these, we  first fix an ordering on $\rts_o$ (the root system attached to $\As_o$) and define the expressions   \be{facts}  u_+&:=& u_+( \{ \bbsigma_{\alpha} \} ):=  \prod_{\alpha \in \rts_{o, +}} \chi_{\alpha}(\bbsigma_{\alpha}) \text{ where } \bbsigma_{\alpha}  = \sum_{k \geq 0} \sigma_{\alpha, k} t^k \in \R[[t]],  \\   u_0 &:=& u_0(\{ \bbsigma_i \}) :=  \prod_{i \in I_o} h_{a_i}(\bbsigma_{i}) \text{ where } \bbsigma_i=1+ \sum_{k \geq 1} \sigma_{i, k} t^k \in (\R[t]])^*,  \text{ and } \\ u_- &:=& u_-(\{ \bbsigma'_\alpha \}) \prod_{\alpha \in \rts_{o, +}} \chi_{-\alpha}(\bbsigma'_{\alpha}) \text{ where } \bbsigma'_{\alpha}  = \sum_{k \geq 1} \sigma'_{\alpha, k} t^k \in t \R[[t]]. \ee Given  $\bim:= (\bbsigma_{\alpha}, \bbsigma_i, \bbsigma'_{\alpha})$ where $\alpha \in \rts_o, i \in I$ and $\bbsigma_{\alpha}, \bbsigma'_{\alpha}, \bbsigma_i$ are as above, we write \be{} u(\bim) := u_+(\{ \bbsigma_{\alpha}\}) u_0(\{ \bbsigma_i \}) u_-(\{ \bbsigma'_{\alpha} \}) .\ee 

\begin{nprop} \cite[(18.11) and Lemma 18.12]{gar:ihes} \label{prop:IM-coordinates-U} Every $u \in \hU$ can be written as $u = u(\bim)$ for a unique family $\bim:= (\bbsigma_{\alpha}, \bbsigma_i, \bbsigma'_{\alpha})_{\alpha \in \rts_{o, +}, i \in I_o}$ as above. If $u \in \hU^{(N)}$ for $N \geq 1$, then $\bbsigma_{\alpha}, \bbsigma'_{\alpha}$ are congruent to $0$ modulo 
	$t^N$ and $\bbsigma_i \equiv 1 \, \mod \, t^N.$ \end{nprop} 

\tpoint{Topology on $\hU$} \label{subsub:U-hat-top} 

The projective limit topology on $\hU$ is the coarsest topology which makes each of the map $\hU \rr \hU/\hU^{(N)}$ continuous, where the (finite-dimensional) group $\hU/ \hU^{(N)}$ is given the standard topology inherited by regarding it as a (unipotent) subgroup of from $\GL(V[N])$. The Proposition \ref{prop:IM-coordinates-U} allows us to identify this topology with the natural (metric) topology on $\hU$ which renders it homeomorphic to $\R^{\mathbb{N}}$. More precisely, we send $u = u(\Sigma)$ to \be{}  \Sigma \in \prod_{\alpha \in \rts_o} \prod_{k=0}^{\infty} \R \times \prod_{\alpha \in \rts_o} \prod_{k=1}^{\infty} \R \times \prod_{i \in I_o} \prod_{k=1}^{\infty} \R. \ee  
With respect to this topology, one has the following:

\begin{nprop} \label{prop:continuity-U-mult} Let $u_1:= u(\Sigma_1)$, $u_2:= u(\Sigma_2)$, and $u_1 u_2 = u(\Sigma_3)$ be three elements from $\hU$. Then the coordinates of $\Sigma_3$ are polynomial in the coordinates of $\Sigma_1$ and $\Sigma_2$. Hence the multiplication map in $\hU$ is continuous. \end{nprop}
\begin{proof} This can be deduced from the results in \cite[\S4.1]{kumar}, but let us just explain the idea for the group $\widehat{\SL_2}$ attached to $\As = \begin{pmatrix} 2 & -2 \\ -2 & 2 \end{pmatrix}$. Writing $\alpha$ the unique positive root of $\SL_2$, the Iwahori--Matsumoto coordinates identify $\hU$ with products of the form \be{} \label{sl2:unip} \chi_{\alpha}(\bbsigma) \, h_{\alpha}(\bbmu) \, \chi_{-\alpha}(\bbtau) \ee where $\bbsigma \in \R[[t]], \bbtau(t) \in t \R[[t]]$ and $\bbmu \in 1 + t \R[[t]]$ (so in particular $\bbmu(t)$ is invertible), and the topology on $\hU$ is just the one built out of the product topology based on the coefficients of $\bbsigma(t), \bbtau(t), \bbmu(t).$  Now we may also consider the algebraic group $\mathbf{\SL}_2$ with coefficients in $\R((t))$; there exists a map $\pi: \hG \rr \mathbf{\SL}_2(\R((t)))$ from \cite[p.64]{gar:ihes} which sends the expression \eqref{sl2:unip} in $\hU$ to the corresponding product in $\mathbf{\SL}_2(\R((t)))$ with \be{} \begin{array}{lcr}   \chi_{\alpha}(\bbsigma)  \mapsto \begin{pmatrix} 1 & \bbsigma \\ 0 & 1 \end{pmatrix}, & \chi_{-\alpha}(\bbtau)  \mapsto \begin{pmatrix} 1 & 0 \\ \bbtau & 1 \end{pmatrix},  & h_{\sigma}(\bbmu)  \mapsto \begin{pmatrix} \bbmu(t) & 0 \\ 0 & \bbmu^{-1} \end{pmatrix}  \end{array}. \ee Moreover, using the argument of \cite[Lemma 18.1]{gar:ihes}, one can show that $\pi|_{\hU}$ is an isomorphism.  It follows from the polynomiality of multiplication in $\pi(\hU)$ that the multiplication in $\hU$ is continuous. 
	
	\end{proof}

\subsection{Parabolic subgroups}

\tpoint{BN-pairs and Bruhat decompositions}  \label{subsub:BN-pair} Let $\bN$ be the subgroup introduced in \S \ref{subsub:elements-in-G}. Then the pair $(\hB, \bN)$ satisfies the axioms of a BN pair or Tits system  \cite[Ch. IV]{bour} and \cite[Thm. 14.10]{gar:ihes}. Hence it follows from general properties of BN pairs that \be{} \label{bruhat:dec} \hG = \sqcup_{w \in \affW} \hB \, \dw \,  \hB. \ee We also have a refined version of this Bruhat decomposition 

\begin{nlem} \cite[Lemma 6.4]{gar:lg2}  \label{refined-bruhat-dec} Fix a reduced decomposition of $w \in \affW$, and hence an ordering on $\Delta_w:= \{ a \in \rts_+\mid w(a) < 0 \}.$ Then every element $x \in B \dw B$ has an expression \be{}  x = b \dw u_w \mbox{ where } b \in B \ee and $u_w = \prod_{a \in \Delta_w} \, \chi_a(u_a) $ is written with respect to the fixed ordering on $\Delta_w$ and $u_a \in \R$.  In this expression, $b$ and $u_w$ are both uniquely determined. \end{nlem} 

We shall denote by $U_w:= U \cap \dw U^- \dw^{-1}$, and note that it is the subgroup generated by $\chi_a(u_a)$ with $a \in \Delta_w$ and $u_a \in \R$. In fact (with respect to a fixed reduced decomposition), we have a map $U_w \rr \R^{\ell(w)}$ which turns out to be a a homeomorphism of topological spaces, where $U_w \subset \hU$ is given the subspace topology.

\newcommand{\X}{\mathscr{X}}
\tpoint{Standard parabolic subgroups of $\hG$}
 \label{subsub:parabolic} From \cite[Chapter IV, \S 2.6, Proposition 4 and \S 2.5 Theorem 3(b)]{bour}, the set of all subgroups of $\hG$ which contain $\hB$ parametrized by subsets $J \subset I$. For each such $J$, we define the corresponding \textit{standard parabolic subgroup} as   \be{} \label{def:parabolic:J} \hP_J:= \la \hB, \, \dw_i, i \in J \ra = \cup_{w \in \affW(J)} \, B \dw B. \ee For $J \neq J'$, these groups are not conjugate and we shall write $\stdp$ for the set of such standard parabolic subgroups (always with respect to our fixed choice of $\hB$). We also write $\stdp'$ for the set of proper, standard parabolics. Note that if $J = \es$, then $P_J =B$.

All of the considerations of this section hold for $\hG_{\k}$ over any field $\k$. We shall need the corresponding groups defined with respect to $\Q$, in which we case they will be denoted as $P_{J, \Q}$, etc.  Note that in the refined Bruhat decomposition \eqref{refined-bruhat-dec}, the group $U_w$ is replaced by $U_{w, \Q}$, i.e. the same group but in which $U_a$ is replaced by $U_{a, \Q}\cong \Q$.

\tpoint{Rational parabolic subgroups} We define the set of rational parabolics   \be{} \begin{array}{lcr} \label{def:rat-parabolics} \ratp:= \{ g \hP_J g^{-1} \mid g \in \hG_{\Q}, J \subset I \} & \mbox{ and } & \ratp':= \ratp \setminus \{ \hG \}. \end{array} \ee For any $P \in \ratp,$ its \textit{type} will be defined as the subset $J \subset I$ so that $P$ is conjugate to $P_J$; the type of well-defined \cite[Ch. IV, \S 2.6, Prop 4(b)]{bour} and for a fixed $J$, the set of rational parabolics of type $J$ is in bijection with $\hG_{\Q}/ P_{J,\Q}$.

\tpoint{Decompositions of $P_J$}  \label{s:parabolics}   For each $J \subsetneq I$ define \be{L:J} \begin{array}{lcr} R_{J} := \la \chi_{a_i}(u), i \in J,  \, u \in \k \ra  & \text{ and } & L_{J} := R_{J} \hT,  \end{array} \ee Writing $w_J$ for the long element in $W_J:=W(\As(J))$, let us also set \be{U:th} \hU_J:= w_J^{-1} \hU w_{J} \cap \hU = \cap_{w \in W_J } w \hU w^{-1}  \ee and note that  the second equality in the first expression is due to \cite[Corollary 6.6]{gar:lg2}) . One can introduce a set of coordinates for $\hU_J$ similar to those on $\hU$ described in \S \ref{subsub:iwahori-matsumoto-coordinates}. We shall not describe this in detail, but remark it can be achieved by replacing $\rts_+$ with $\pb_J$ in the constructions in \emph{op. cit.}. Now, the group $L_J$ normalizes $\hU_J$ and by \cite[Theorem 6.1]{gar:lg2} we have \be{} \label{lev:dec}  P_J = L_J \ltimes \hU_J   . \ee

\tpoint{Parabolic subsets $\hA^+_J$} \label{subsub:A-J}   For each $J \subset I$,  define \be{}\label{def:A_Je} \begin{array}{lcr} (\hA^s)_J:= \{ h \in \hA \mid (h \eta(s))^{a_i}=1 \text{ for } i \in J \} & \mbox{ and } & (\hA^+)_J:= \cup_{0 <s < 1} \hA^s. \end{array}  \ee We also define the subset of $\hA^e$ (keeping in mind the placement of the parenthesis)  \be{} \hA^+_J:= \cup_{0 <s<1} \{ h \eta(s) \mid h \in (\hA^s)_J \} \ee corresponding to (subsets) of Lie algebras $(\hfh^r)_J$, $(\hfh^+)_J$ and $\hfh^+_J$ respectively that were introduced in \S \ref{subsub:cartan-decompositions}.  If we further define $ \hA^+(J):= \hA(J):= \la h_{a_i}(s) \mid i \in J \ra$, then $ \Lie \hA(J) = \hfh^+(J)$ and by exponentiation of  \eqref{hr:dec}, we obtain , again with uniqueness of expression, for each  \be{} \label{Ac:dec}  \hA=  \Ac(J) \times (\Ac^s)_J. \ee Note as before that this decomposition depends on $s$ even though the dependence on the left hand side is not displayed, i.e. $h \in \hA$ sent to $(h(J), h_J)$ with $h(J) \in \hA(J)$ and $h_J \in (\hA^s)_J$ where \be{} h \eta(s) = h(J) \left( h_J \eta(s) \right) \mbox{ with }  \left( h_J \eta(s) \right)^{a_i}=1 \mbox{ for } i \in J. \ee In particular, it is \textit{not} (usually) the case that $h(J)$ depends only on $h$.  For $J \subset K \subsetneq I$, we also define  \be{} \Ac(J, K):= \{ \exp(H) \mid H \in \hfh(J, K) \} = \{ \exp(H) \mid H \in \hfh(K)_{\sta{J}} \} = \hA(K)_{\sJ}. \ee In analogy with  \eqref{eq:relative-lie-algebra} and \eqref{dec:hrJ:K}, we have for $J \subset K \subsetneq I$,  \be{} \begin{array}{lcr} \label{Ac:dec relative} (\hA^s)_{J} = \hA(J, K) \times (\hA^s)_{K} & \mbox{ and } & \Ac(K) = \Ac(K)_{\sJ} \times \Ac(J). \end{array} \ee From the former, we also obtain $  \hA^s_J = \hA(J, K)  \times \hA^s_k$ which we often write as \be{}  \label{A+:dec-rel} a_J \eta(s) = \left( a(J, K),  a_K \eta(s)  \right) \ee where $a_J \eta(s) \in \hA^s_J$, $a_K \eta(s) \in \hA^s_K$ and $a(J, K) \in \hA(J, K)$. Finally, if we are given subsets $J \subset K_1 \subset K_2 \subsetneq I$, then we also have the relative decomposition  \be{} \label{a:j-k1-k2}  \hA(J, K_2) &=& \hA(J, K_1) \times \hA(K_1, K_2) \mbox{ for } J \subset K_1 \subset K_2 \subsetneq I \\ &=&  \hA(J, K_1) \times \hA(K_2)_{\sta{K_1}}. \ee

 \tpoint{Langlands decompositions of $P$} \label{subsub:Langlands-P_J} Recalling that $\eta(s)$ normalizes $L_J$ for all $J \subset I$, define  \be{} Z^s({L_J}) := \{ x \in L_J \mid (x \eta(s)) y (x \eta(s))^{-1} =y \mbox{ for all } y \in L_J \}. \ee  One can check that the connected component of the identity of $Z^s(L_J)$ is equal to the group $(\hA^s)_J$ introduced in \eqref{def:A_Je}.  
  Moreover, there exists a finite cover of the group $R_J$, denoted $M_J$ such that \be{} \label{LJ:prod} L_J  \cong M_J \times (\hA^s)_J. \ee In fact writing $\Ratc(L_J)$ for the set of rational character of $L_J$, we may choose \be{} M_J:= \cap_{\chi \in \Ratc(L_J)} \ker(\chi^2). \ee   Combining \eqref{LJ:prod} with \eqref{lev:dec}, we obtain the decomposition, for each $0 < s < 1$, \be{} \label{langlands-dec} P_J  = M_J \times (\hA^s)_J   \times \hU_J \ee If $P=P_J \in \stdp$, we just write $M_P$ for $M_J$, $(\hA^s)_P$ for $(\hA^s)_J$, and $\hU_P$ for $\hU_J$.

 \label{subsub:general-P-Levi} For $P \in \ratp'$, say  $P = g P_J g^{-1}$ for $g \in \hG_{\Q}$, one obtains similar decompositions. Indeed, first check that $\hU_P:= g \hU_Jg^{-1}$ is independent of the class of $g \in \hg_{\Q}/P_{\Q}$, and similarly for $L_P:= g^{-1} L_{P_J} g$.  We also define as before \be{} Z^s(L_P):= \{ x \in L_P \mid (x \eta(s)) y (x \eta(s))^{-1} = y \mbox{ for all } y \in L_P. \} . \ee The connected component of the identity of $Z^s(L_P)$ will be denoted as $(\hA^s)_P$. Let us also set \be{} \label{A+P} \hA^+_P= \cup_{0 <s < 1} ( \hA^s)_P \ee and write $\Ratc(L_P)$ for the set of rational character of $L_P.$ As before if we define \be{} M_P:= \cap_{\chi \in \mathrm{Rat}_{\Q}(L_P)} \ker(\chi^2), \ee then we obtain, for every fixed $s$ a decomposition for each $P \in \ratp'$, \be{} \label{rat-P:langlands-dec} \begin{array}{lcr} L_P \cong M_P \times (\hA^s)_P, & \mbox{ and } & P = M_P \times (\hA^s)_P \times \hU_P. \end{array} \ee

As a notational shortcut, it $J$ corresponds to a parabolic $P:=P_J$, we now often just write $(\hA^s)_P$ for $(\hA^s)_J$, etc. 
 
 \tpoint{Relative decompositions.} \label{subsub:relative-langlands-dec} We will especially be interested in relative versions of the decompositions \eqref{langlands-dec} and \eqref{rat-P:langlands-dec}. Let $Q, P \in \ratp$ with $Q \subset P$.  Denote by $\sta{Q}:= Q \cap M_P \subset M_P;$ it is a parabolic subgroup of $M_P$ and  is equipped with a decomposition \be{} \label{relative-langlands-dec} \sta{Q} = M_{\sta{Q}} \times \hA(P)_{\sta{Q}} \times U_{\sta{Q}}. \ee One can check that if $P= M_P \times (\hA^s_P) \times \hU_P$ as in \eqref{langlands-dec} and  $Q = M_Q \times (\hA^s_Q) \times \hU_Q$  with  \be{} \label{Q:starP}  \begin{array}{lccr} M_Q = M_{\sta{Q}}, & (\hA^s_{Q}) = (\hA^s_{P}) \hA(P)_{\sta{Q}}, & \mbox{ and } & \hU_Q = \hU_P U_{\sta{Q}} \end{array}. \ee Conversely if $\sta{Q} \subset M_P$ is a parabolic with decomposition \eqref{relative-langlands-dec}, then defining $M_Q, (\hA^s_Q), \hU_Q$ as in \eqref{Q:starP}, the group $Q = M_Q \times (\hA^s_Q) \times \hU_Q$ is a rational parabolic subgroup of $\hG$ which is contained in $P$. In general, these correspondences are inverses to one another, \textit{i.e.} we have: 
 
 \begin{nprop} Let $P \in \ratp$. Then there is a bijective correspondence between $\ratp(P)$, the rational parabolic subgroups contained in $P,$ and $\ratp(M_P)$, the rational parabolic subgroups of $M_P$. Denoting this correspondence by \be{} \begin{array}{lcr}Q \mapsto \, \sta{Q}:= Q \cap M_P & \mbox{ for }  Q \in \ratp(P), \mbox{i.e. }  Q \subset P \end{array} \ee  the relations \eqref{Q:starP} hold. \end{nprop}

\newcommand{\tq}{(q)}
\newcommand{\pd}[1]{^{(#1)}}

\subsection{The corner $\cor(\hA^+)$ attached to $\hA^+$}  \label{sub:bord-A} 
We can now formulate the analogue of the results from \S \ref{sub:bord-h} for the groups $\hA^+_P$ with $P \in \ratp'$. Recall the notations introduced in \S \ref{sub:bord-h}-\S \ref{sub:top-corner} as well as that of the previous section, which we adopt freely throughout this subsection.  We have already treated the simply connected case in \S \ref{subsub:affine-corners}. As the general case is not much different, we just state the results.

\tpoint{On $\cor(\hA^+_J)$} \label{subsub:bord-Aplus} 

For each $J \subsetneq I$, define the set \be{} \label{corner-AJ}  \cor( \Ac^+_{J}):= \ \ \cup \ \ \bigsqcup_{J \subset K \subsetneq I} \Ac(J,K) = \Ac^+_{J} \ \ \cup \ \ \bigsqcup_{J \subset K \subsetneq I} \Ac(K)_{\sta{J}}. \ee This set can be equipped with a natural topology using a family of convergent sequences as in Proposition \ref{prop:top-on-cor(A)}. More precisely, we say that:

\begin{itemize}
	\item A sequence $a_n\eta(s_n)\in \hA^+_J$ converges to $a_{\infty}$ with $a_{\infty} \in \hA(J, K)$ for $J \subset K \subsetneq I$ if, when using \eqref{A+:dec-rel}, we write $a_n\eta(s_n) = (a(J, K)_n, a_{K, n} \eta(s_n) )$ with $a(J, K)_n \in \hA(J,K)$ and $a_{K, n} \in \hA^{s_n}_K,$ then \be{} \begin{array}{lcr} (a_{K, n} \eta(s_n))^{a_i} \rr 0 \mbox{ for } i \in I \backslash K & \mbox{ and  } & a(J, K)_n \rr a_{\infty}. \end{array} \ee 
	\item If $J \subset K_1 \subset K_2 \subsetneq I$ then we say that $\{ a_n(J, K_2) \}$ in $\Ac(J, K_2)$ converges to  $a_{\infty}(J, K_1) \in \Ac(J, K_1) $ if, when using \eqref{a:j-k1-k2} to write $a_n(J, K_2)  = a_{n}(J, K_1) \cdot a_n(K_2)_{\sta{K_1}}$  with $a_{n}(J, K_1) \in \hA(J, K_1)$ and $a_n(K_2)_{\sta{K_1}} \in \hA(K_2)_{\sta{K_1}}$, then \be{} \begin{array}{lcr} (a_{n}(K_2)_{\sta{K_1}})^{a_i}   \rr 0 \mbox{ for } j \in K_2 \backslash K_1 &  \mbox{ and  } & a_n(J, K_1)_n \rr a_{\infty}(J, K_1). \end{array} \ee
\end{itemize}

\noindent The argument as in Proposition \ref{basis bordh}  shows that the above endows $\cor(\hA^+)$ with the structure of a Hausdorff topological space. One could also define the topology first on $\corn{\A^+}$ and then argue that $\corn{\hA^+_J}$ is the interior of the closure of $\hA^+_J$ in this space.

 For a rational parabolic subgroup $P \in \ratp'$, we have introduced the set $\hA^+_P$ in \eqref{A+P}. In the notation introduced above, we then have   \be{} \label{pc:APloop} \cor(\Ac^+_{P}) := \Ac^+_{P} \cup \bigsqcup_{P \subset Q} \hA(Q)_{\sP} \ee with a topology introduced as before.

\tpoint{Topological properties of $\cor(\hA_P^+)$}    \label{subsub:AP-bordification} Let $P, Q \in \stdp'$ with 
$Q:=P_K \supset P:=P_J$ with $K \supset J.$ For  $\sigma > 0$ we now define \be{} \label{A:P:c} \hA^+_{P,\sigma} & :=&  \{ x \in \hA^+_P  \mid x^{a_i} < \sigma \text{ for } i \in  I \setminus J \} \mbox{ and } \\ \hAp(Q)_{^*P,\sigma}  &:=& \{h \in \hAp(Q)_{\sP}  \mid  h^{a_i} < \sigma \mbox{ for } i \in K \setminus J   \}. \ee Using the description \eqref{pc:APloop} to extend the functional $h \mapsto h^{a_i}, \ i \in I \setminus J$ to all of $\cor(\hAp_P)$  by allowing $h^{a_i}$ to assume the value $0$,  we now set \be{} \label{bord:A:sigma} \cor(\hA^+_{P,\sigma}):= \{ x \in   \cor(\hA^+_{P}) \mid x^{\alpha} < \sigma \text{ for } \alpha \in \Delta_P\} .  \ee

\begin{nlem}\label{lem:closure-AP:c} For $P \in \ratp,$  we have \be{} \label{AP,t} \cor(\hA^+_{P,\sigma}) = \Ac^+_{P,\sigma} \, \cup \, \bigsqcup_{P \subset Q} \hAp(Q)_{^*P,\sigma} \ee and moreover $\cor(\Ac^+_{P,\sigma})$ is the interior of the closure of $\Ac^+_{P,\sigma}$ inside $\cor(\hAp_P^+).$ \end{nlem} These sets are not compact, since $\hA^+_{P, \sigma} $ is unbounded in the central direction, i.e. if $x \in \hA^+_{P, \sigma}$, then $x^{\lambda_{\ell+1}}$ can assume any value in $\R$.

For $P \in \stdp'$, we can use collections as in Proposition \ref{basis bordh} to obtain a natural, first countable basis. One can also argue that $\corn{\hA^+_P}$ is contractible and homotopic to its interior. 
Finally, for real numbers $E < F$ and $0< s_0 < 1$ and $\sigma >0$ let us define $\hA^+_P([E, F], s_0, \sigma)$ to be \be{} \{ x \in \hA^+_P \mid x^{a_i} \leq \sigma, i \in I \setminus J \} \cap \{ x \in \hA^+_P \mid E \leq x^{\lambda_{\ell+1}}  \leq  F  \} \cap \{ x \in \hA^+_P \mid x^{\delta} \geq s_0 \}. \ee  Its closure in the corner $\corn{\hA^+_P}$ is then seen to be equal to (where $Q=P_K$ and $P=P_J$ as before) \be{} \label{compact-corner-A} \corn{\hA^+_P([E, F], s_0, \sigma)} =  \hA^+([E, F], s_0, \sigma, ) \, \cup \, \bigsqcup_{Q \supset P} \{ h \in \hA(Q)_{\sta P} \mid x^{a_i} \leq \sigma \mbox{ for } i \in K \setminus J \}. \ee  
\begin{nlem} \label{lem:corner-A-cpt} For real numbers $s_0 >0, \sigma > 0$ and $E < F$ the set $\corn{\hA^+_P([E, F],s_0, \sigma)}$ is compact.  \end{nlem} The proof is similar to that of  Lemma \ref{lemma:compact bordh}.

\subsection{Iwasawa and Horospherical Decompositions} 

\newcommand{\xgs}{X_{\hG^s}}

\tpoint{Iwasawa decomposition} \label{subsub:Iwasawa-decomposition}  Using the positive-definite inner product $\{ {\cdot, \cdot} \}$ on $V:=V^{\lambda}$ (the defining, highest-weight representation of $\hG$) from \S \ref{subsub:hermitian}, define   \be{} \label{def:K} \hK:= \{ g \in \hG \mid \{ gv, gw \} = \{ v, w \} \text{ for all } v, w \in V \}. \ee It plays the role of the standard `maximal compact subgroup' from classical Lie theory, though we are not asserting any topological properties about $\hK$. The following is the loop group analogue of the Iwasawa decomposition.

\begin{nprop} \cite[(16.14)]{gar:ihes} \label{prop:Iwasawa} For $g \in \hG, s \in \R^*$,  there exists $k_g \in \hK, u_g \in \hU$ and $a_g \in \hA$ such that \be{} \label{g:iwasawa} g \eta(s) = k_g \, a_g \eta(s) \, u_g \ee where $k_g, u_g, a_g$ are uniquely determined from $g$. \end{nprop} 

\tpoint{Loop group symmetric spaces $\xg$}  \label{subsub:xg-def} Defining the sets sets  \be{}  \label{def:X} \begin{array}{lcr} X_{\hG^s}:= \hK\setminus \hG \eta(s) & \text{ and } & \xg:= \cup_{0 < s < 1 } \xgs, \end{array} \ee the Iwasawa decomposition from the previous paragraph gives bijections
 \be{} \label{horo:borel} \xg =  \hA^+ \times \hU. \ee As was described in  \ref{subsub:U-hat-top}, we can equip $\hU$ with a natural topology and the set  $\hA^+$ (which can be identified with a subset of $\R_{>0}^{\ell+2}$) also carries a natural metric topology. We can then equip $\xg$ with the product topology coming from \eqref{horo:borel}.
 
\tpoint{$P$-horospherical coordinates} \label{subsub:P-horospherical} For $P \in \ratp'$, recall that we have decompositions \eqref{langlands-dec} \be{} \begin{array}{lcr} P  = M_P \times (\hA^s)_P \times \hU_P & \mbox{ or } & P \eta(s) = M_P \times \hA^s_P \eta(s) \times \hU_P \end{array}. \ee  Generalizing \eqref{horo:borel}, we now obtain the following,

\begin{nprop}[$P$-horospherical coordinates] \label{prop:horo} 
	
	\begin{enumerate}
		\item Let $P, Q \in \ratp'$ with $P \subset Q.$ Then  
	\be{} \label{split glob P} \begin{array}{lcr}  \xg = X_{M_Q} \times \hA^+_Q \times \widehat{U}_{Q} & \mbox{ and } & X_{M_Q} = X_{M_P} \times \hA(Q)_{^*P} \times \hU_{^*P} \end{array}, \ee where $X_{M_P} := (\hK \cap M_P ) \setminus M_P.$ 
		\item The projection maps from $\xg$ to each of the factors above is continuous with respect to the topology on $\xg$ introduced in \S \ref{subsub:xg-def}, the natural metric topologies on $X_{M_Q}$, $\hA^+$, and the pro-unipotent topology on $\hU_P$. A similar results holds for $X_{M_Q}$ and the second factorization.
	
	\end{enumerate}

\end{nprop}

\tpoint{ On the functions $H_P$ and $H^+_P$}  Using Proposition \ref{prop:horo}, we can define the functions \be{} \label{def:HP} \begin{array}{lcr} H_P: \xg \rr (\hfh^+)_P & \mbox{ and } &  H^+_P: \xg \rr \hfh^+_P \end{array} \ee in the following manner. Assume first that $P=P_J$ is standard. Write $x \in \xg$ according to the first descomposition in \eqref{split glob P}, say $x = (m, a^+, u)$ with $m \in X_{M_P}$, $u \in \hU_P$, and $h^+ \in \hA_P^+$, we define $H^+_P(x) := \log a^+$ to be the unique element in $\hfh_P^+$ such that if  $H_P(x) = \sum_{j \notin J} c_i \lv_i + f \cc$, then \be{} \begin{array}{lcr} c_i = \log \, \, (h^+)^{a_i} & \mbox{ and } & \log \, \,  (h^+)^{\lambda_{\ell+1}} = f \end{array}. \ee In case only some multiple of $\lambda_{\ell+1}$ lies in $\Xi_{\lambda}$, we just modify the second condition above appropriately to extract the coefficient of a multiple of $\cc$. If $h^+ = h \eta(s)$ with $h \in \hA$ and $s=e^{-r}$, we also define $H_P(x) \in (\hfh^+)_P$ so that \be{} H_P(x) - r \dd = H_P^+(x). \ee   
By conjugation we can also define $H_P$ and $H^+_P$ for any rational parabolic, and via a similar procedure we also obtain maps $H_ {\sP}: X_{M_Q} \rr \hfh(Q)_{\sP}$ for $P \subset Q$ as above. 

\begin{ncor} \label{lem:H-PQ}
     For a pair of parabolic subgroups $P \subset Q$ with $P, Q \in \ratp'$ and $x \in \xg$ we have,
     \be{} \begin{array}{lcr} \label{eq:lem-H-PQ}  H_P(x) = H_Q(x) + H_{^*P}(\pr_{M_Q}(x)) & \mbox{ and } & H^+_P(x) = H^+_Q(x) + H_{^*P}(\pr_{M_Q}(x)) \end{array} \ee where $\pr_{M_Q}: \xg \rr X_{M_Q}$ is the projection from \eqref{split glob P}.  \end{ncor}
\begin{proof} The second result clearly follows from the first, so let us just focus on the first. It follows by going from the $Q$-decomposition of $X$ to the $P$-decomposition using \ref{split glob P} and then comparing the values $H_P(x)$ and $H_Q(x)$. \end{proof}

\subsection{Orthogonal Families}

\label{sec:orthogonal-families}

\renewcommand{\hb}{B}
\newcommand{\bo}{\mc{B}}

\tpoint{Orthogonal families for $\hG$}  \label{subsub:orthogonal-families} Recall the notion of (positive) Borel subsets or $\rts$ and their classification given in Proposition \ref{prop:base-pb-aff}. Denote by $\bo$  the set of all \textit{positive} Borel subsets of $\rts$, or, just Borel subsets for short. Let $\bor_e:= \rts_+$, and note that Proposition \ref{prop:base-pb-aff} asserts that each $\bor \in \bo$ is of the form $\bor:=\bor_w:= w \bor_e w^{-1}$ for a unique $w \in \affW$.  Two elements $\bor_1, \bor_2 \in \bo$ are said to be \textit{adjacent} if $\bor_1  \cap - \bor_2$ consists of a single element, say $\alpha:=\alpha(\bor_1, \bor_2)$.   If $w, \sigma \in \affW$ then we deduce that $\bor_w$ and $\bor_{\sigma}$ are adjacent if and only if $\sigma = w w_{\alpha}$ for $\alpha \in \simp$ (a simple root) such that $w(\alpha)>0$. Moreover, in this case we have$\alpha(\bor_{w}, \bor_{w_{\alpha} w }) = w(\alpha).$

\begin{define} A collection of elements $\sY= \{ Y_\bor\}_{\bor \in \bo}$ where each $Y_{\bor} \in \hfh^e$ is said to form an \textit{orthogonal family} (in finite dimensions, the notion is due to Arthur, see \cite[\S 1.5]{friday-morning}) if for any pair of adjacent Borels  $\bor_1, \bor_2 \in \bo$ we have \be{} \label{orthogonal-family-relation} \begin{array}{lcr} Y_{\bor_1} - Y_{\bor_2} = r(\bor_1, \bor_2) \, \alphav(\bor_1, \bor_2) & \mbox{ for } & r(\bor_1, \bor_2) \geq 0 \end{array}  \ee for $\alphav(\bor_1, \bor_2)$ the coroot attached to $\alpha(\bor_1, \bor_2)$ (see \S \ref{subsub:coroots}). We furthermore say that $\sY$ is a \textit{regular} orthogonal family if $r(\hb_1, \hb_2) >0$  for all $\bor_1, \bor_2$ adjacent. \end{define}

\begin{nrem}   \label{remark:complementary-polyhedron} Let $(V, \Phi)$ be a finite type root system, so $\Phi$ is the set of roots, $V$ is a Euclidean space containing $\Phi$, $\Delta$ a fixed base of simple roots, $\Phi_+$ the corresponding positive roots, etc. Writing  $L_{\alpha}:= \{ \xi \in V \mid \la \xi, \alphav \ra = 0 \}$ for the root hyperplanes,  we let $\cham(V)$, the chambers of $V$, be the set of connected components of the complement of $L_{\alpha}$ for $\alpha \in \Phi$. These are in natural bijection with Borel subsets of $\Phi$, or equivalently, with the Weyl group $W:=W(\Phi)$ of the root system. The chamber corresponding to $w=e$ is called the fundamental chamber, and all other chambers are of the form $C_w:= w C_e$. Write $\Delta \subset \Phi$ for a set of simple roots, enumerated as $ \alpha_1, \ldots, \alpha_r$, and also set $\lambda_1, \ldots, \lambda_r$ for the corresponding set of fundamental weights. If $C \in \cham(V)$, write $\vertices(C)$ for the fundamental weights with respect to $C$, i.e. $\vertices(C)= w^{-1} \{ \lambda_1, \ldots, \lambda_r \}$ if $C:= w C$ with $w \in W.$   Behrend \cite[Def. 2.1]{beh} defines a  \textit{complementary polyhedron} as a collection of elements $\sY= \{ Y_C\}_{C \in \cham(V)},$ where each $Y_C \in V^*,$ satisfying the conditions: 
\begin{enumerate}
	\item[\textbf{CP1.}] If $C, D \in \cham(V)$, with $\lambda \in \vertices(C) \cap \vertices(D)$, then  $\la \lambda, Y_C - Y_D \ra =0.$ 
	\item[\textbf{CP2.}] If $C$ and $D$ are adjacent chambers\footnote{This means that $\vertices(C)$ and $\vertices(D)$ have all but one element in common} with $\alpha \in \Delta$ the unique simple root such that $ D  = w_{\alpha}(C)$,  \be{} \la \alpha, Y_C - Y_D \ra \leq 0. \ee  
\end{enumerate}	
The actual polyhedron is then the convex hull of the points $\{ Y_C \}$ . This definition can be seen to be equivalent to the notion of an orthogonal family (for finite-type root systems). Indeed, as we already mentioned the chambers $\cham(V)$ are in bijection with Borel subgroups. Moreover, each chamber defines a set of positive roots and hence a Borel subset, with $C_e$ corresponding to $\Phi_+$. The condition \textbf{CP1} is equivalent to the similar condition formulated for adjacent chambers, in which case it states that $Y_C - Y_D$ is a multiple of the unique coroot positive for $C$ and negative for $D$. Condition \textbf{CP2} then ensures that this multiple is positive. Using these ideas, one can prove an equivalence between Behrend's notion of a complementary polyhedron and Arthur's notion of an orthogonal family. 
\end{nrem}

\tpoint{Basic Inequality for Orthogonal Families} \label{subsub:orthogonal-family-basic-inequality} Suppose $\sY=(Y_\bor)$ is an orthogonal family and let $\leq_{\bor}$ be the dominance order with respect to $\bor \in \bo,$ i.e. \be{} \begin{array}{lcr} \label{def:dleq}   x \leq_{\bor} y & \mbox{if and only if} & y - x = \sum_{\alpha \in \bor} n_{\alpha} \alphav \mbox{ with } n_{\alpha} \geq 0  \end{array}. \ee If $\bor=\rts_+$, we just omit it from our notation and call this the \textit{standard} dominance order. 

\begin{nlem} (Orthogonal Family Inequality) \label{lem:orthogonal-family-ineq} For $\sY = (Y_{\bor})$ an orthogonal family for $\hG^e$, \be{} \begin{array}{lcr} \label{YB-ineq} Y_{\bor'} \leq_{\bor} Y_\bor & \mbox{ for any } & \bor, \bor' \in \bo \end{array} \ee \textit{i.e.} $Y_{\bor} - Y_{\bor'}$ is a non-negative linear combination of the roots  in $\bor \cap - \bor'.$ \end{nlem}
\begin{proof} If $\bor$ and $\bor'$ are adjacent, then $Y_{\bor} - Y_{\bor'}$ is a positive multiple of the unique root positive for $\bor$ and negative for $\bor'$. We argue by induction on the length of $w$ such that $w \bor = \bor'$ as follows: choose a chain $\bor':= \bor_0, \bor_2, \ldots, \bor_r:=\bor$ such that each $\bor_i$ and $\bor_{i+1}$ are adjacent. Then we have \be{} Y_{\bor'}:=Y_{\bor_0} \leq_{\bor_{1}} Y_{\bor_{1}} \leq_{\bor_{2}} \leq \cdots \leq_{\bor_r} Y_{\bor_r}=Y_{\bor}. \ee Choosing a reduced decomposition $w:= w_{b_1} \cdots w_{b_r}$ with each $b_i \in \Delta$, we may obtain one such a chain by setting $\bor_0:=w \bor'$, $\bor_{1}= w_{b_r} \cdots w_{b_{r-1}} \bor', \ldots, \bor_r= \bor.$  The result follows using \eqref{adjacent-pos-Weyl}. \end{proof}

\tpoint{Example: Weyl group orbits} Suppose $T \in \hfe$ is dominant, i.e $\la T , \alpha_i \ra \geq 0$ for all $i \in I$. Then the collection $\sY:= (Y_w)$ with $Y_w = w \cdot T$ for $w \in \affW$ forms an orthogonal set. Indeed, if $\bor_1=B_w, \bor_2=B_{w w_{\alpha}}$ are two adjacent Borel subsets, i.e. $\alpha \in \simp$, and $w(\alpha) > 0,$ , then \be{} Y_{\bor_1} - Y_{\bor_2} = w T - w w_{\alpha}  T = \la T, \alpha \ra w(\alpha). \ee  If $T$ is regular and dominant, then the orthogonal family is also positive. 

\tpoint{Example: $H^+_B(x)$} Let $x \in \xg$, and consider the collection $\sY_x:= \left( w \cdot  H^+_B(x \, \dw )\right)_{w \in W}.$ We claim that it forms a positive orthogonal family. One may also write more invariantly: $ \sY_x:= \left( H_{\bor}(x) \right)_{\bor \in \bo}.$ Indeed, if $\bor_w:= \dw^{-1} \rts_+ \, \dw$ corresponds to the Borel $\leftidx^w \hB:= \dw^{-1} B \dw$, then  $x = k b$ with $k \in \hK, b \in \hB$, implies $x \cdot \dw = k' \dw^{-1} b \dw$ with $k' := k \dw \in \hK$ and $\dw^{-1} b \dw \in \leftidx^wB,$ \textit{i.e.} \be{} H^+_B(x \dw) = w^{-1} \cdot H^+_B(x). \ee The properties of a positive orthogonal family follow from  \cite[Lemma 6.1]{gar:ms-2}.

\begin{ncor}\label{cor:orthogonal-ineq} Let $x \in \xg$. Then for any $w \in W$ we have \be{} \label{ortho-ineq} w \cdot H^+_B(x \dw)  \leq H^+_B(x) \ee \end{ncor}

\begin{nrem} In \cite[Corollary 3.14]{beh} Behrend states a certain uniquness theorem for any complementary polyhedron that allows him to attach to each such object a unique facet. Our main results in \S\ref{sec:canonical-pairs} suggests that such a generalization should exist, but we have not pursued this point here (Behrend's proof does not seem to generalize in a simple fashion). Note that as the above example suggests, one no longer obtains a bounded convex region in the affine case, rather an infinite parabolic region.  \end{nrem} 

\subsection{The group $\hGam$ and Garland's reduction theory } \label{sec:reduction}

\tpoint{The group $\hGam$} Recall that we have fixed an integral form $V_{\zee}:= V_{\zee}^{\lambda} \subset V_{\R}$ in \S 3.1 while defining the group $\hG:= \hG^{\lambda}$. We now \be{} \label{def:hGam}  \hGam:= \{ \gamma \in \hG \mid \gamma ( V_{\zee})  \subset V_{\zee} \} \ee as well as its subgroup\footnote{There is some evidence this is equal to $\hGam$ actually (we believe this to be the case), see \cite{strong-integrality} }  \be{}  \hGam_0:= \la \chi_{\alpha}(\bbsigma) \mid \alpha \in \rts_o, \sigma \in \zee((t)) \ra. \ee  Noting that $\hG_{\Q}$ the group introduced in \S \ref{def:G-hat} is also equal to  \be{}  \hG_{\Q}= \{ g \in \hG_{\R} \mid g (V_{\Q}) \subset V_{\Q} \} \ee it follows (see \cite[(2.5)]{gar:ms-2}) that \be{} \hG_{\Q} = \hGam_0 \, \hB_{\Q} \ee where $\hB_{\Q} = \hB \cap \hG_{\Q}.$   From this observation and from \cite[Prop. 1]{gar:ms-2}, the natural inclusions $\hGam_0 \hookrightarrow \hGam \hookrightarrow \hG_{\Q}$ induces bijectons \be{} \hG_{\Q} / \hB_{\Q} \stackrel{1:1}{\longleftarrow} \hGam_0 / \hGam_0 \cap \hB \stackrel{1:1}{\longrightarrow} \hGam / \hGam \cap \hB. \ee Similarly, for $P \in \stdp$, writing $P_{\Q}:= P \cap \hG_{\Q}$ that we have bijections \be{} \hG_{\Q} / \hP_{\Q} \stackrel{1:1}{\longleftarrow} \hGam_0 / \hGam_0 \cap \hP \stackrel{1:1}{\longrightarrow} \hGam / \hGam \cap P. \ee We shall henceforth also write \be{} \begin{array}{lcr} \hGam_P:= \hGam \cap P & \mbox{ and } & \hGam_{0 , P} = \hGam_0 \cap P \end{array} \ee

\spoint For $P \in \ratp$, recall the (reductive) group $L_P$ defined in \S\ref{subsub:general-P-Levi} and its subgroup  \be{} M_P:=  \bigcap\limits_{ \chi \in \Ratc(L_P)} \ker(\chi^2). \ee Write $\hgam_{L_P}$ for the image of $\hgam_P:= \hgam \cap P$ in $L_P$ under the natural projection $P = L_P \rtimes \hU_P\rr L_P$. As in \cite[Prop 1.2]{borel-serre}, we can prove

\begin{nlem}\label{lem:GammaL} The group $\Gamma_{L_P}$ is contained in $M_P$ and is an arithmetic subgroup of $M_P$ (a reductive, finite-dimensional group). Hence we shall write $\Gamma_{M_P}:= \Gamma_{L_P}$ from now on.  \end{nlem}

\begin{proof} By conjugation, we may assume $P$ is standard, say $P=P_{J}$ with $J \subsetneq I$. Then for each $i \notin J$, we have an element $\chi_i  \in \Ratc(L_P) $ defined by $x \mapsto x^{\lambda_i}$, where $\lambda_i$ is the corresponding fundamental weight defined in \eqref{Lambda:j}. In fact the $\chi_i, i \notin J$ form a basis of $\Ratc(L_P),$ so it suffices to show, for any $\gamma \in \hgam_{L_P}$ that $\chi_i(\gamma) = \pm1.$ In fact, if we write $\gamma = a_{\gamma} m_{\gamma}$ with $a \in \hA_{P}$ and $m_{\gamma} \in M_P$ according to (the $s=0$ version of) \eqref{LJ:prod}, it suffices to show that $a_{\gamma}^{\lambda_i} = 1$ for all $i \notin J$.

	For  $i \notin J$, there exists a positive integer $d$ and a homomorphism $\pi_i: \hg \rightarrow \hg^{d \lambda_i}$ (see Theorem \ref{thm:change-of-weights}). Under this map, we clearly have $\hGam_0^{\lambda}$ gets sent $\hGam_0^{d \lambda_i}$. Since $L_P$ is a reductive group, we know (see  \cite{strong-integrality} and references therein) that \be{} \begin{array}{lcr} \hGam_0^{d \lambda_i} \cap L_P^{d \lambda_i} = \hGam^{d \lambda_i} \cap L_P=: \hGam^{d \lambda_i}_{L_P} & \mbox{ and } & \hGam^{\lambda} \cap L_P= \hGam_0^{\lambda} \cap L_P =: \hGam^{\lambda}_{L_P} \end{array}. \ee 
	
	Now applying $\gamma= a_{\gamma} m_{\gamma} \in \hGam_{L_P}^{\lambda}$  to a primitive highest weight vector $\hw_{d \lambda_i}$ in $V^{d \lambda_i}$ and we find that 
	 \be{} \gamma . \hw_{d \lambda_i } = (a_\gamma)^{d \lambda_i} \hw_{d \lambda_i }. \ee 
	\noindent Since $\gamma \in \hGam_{L_P}$ we must have $a_{\gamma}^{d \lambda_i} \in \zee_{>0}$. Applying the same to $\gamma^{-1} \in \hGam_{L_P}$, it follows that $(a_{\gamma}^{d \lambda_i})^{-1} \in \zee_{>0}$ as well. Hence $a^{d \lambda_i}_{\gamma}=(a_{\gamma}^{\lambda_i})^d=1$ and so also $a^{\lambda_i}_\gamma=1$. The same works for all $i \notin J$.  \end{proof}
\tpoint{Reduction theory for $\hU$} \label{subsub:reduction-theory-U} Let $\hGam_{\hU}:= \hGam \cap \hU$ and recall the coordinates $\Sigma= ( \bbsigma_{\alpha}, \bbsigma_i, \bbsigma'_{\alpha})$ introduced on $\hU$ from \S \ref{subsub:iwahori-matsumoto-coordinates}. Write \be{} \label{omega-0} \Omega_0:= \{ x \in \R\mid |x| \leq 1/2 \} \ee and consider the sets \be{} \begin{array}{lcr} \Omega_0[[t]]:= \{ \sum_{j \geq 0} q_j t^j \mid q_j \in \Omega_0 \}, &  \mbox{ and } & (\Omega_0[[t]])^*:= \{ \sum_{j \geq 0} \, q_j t^j \mid q_0 =1, q_j \in \Omega_0 \}. \end{array} \ee 

\begin{nlem} \label{u-approx} \cite[Lemma 18.16]{gar:ihes} For any $u \in \hU$, there exists $\gamma \in \hGam_{\hU}$ such that \be{} u \gamma \in \hU_{\Omega_0} \ee where $\hU_{\Omega_0}$ consists of $u(\Sigma)$ with $\Sigma =(\bbsigma_{\alpha}, \bbsigma_i, \bbsigma'_{\alpha})$ with $\bbsigma_{\alpha} \in \Omega_0[[t]]$, $\bbsigma'_{\alpha} \in t\Omega_0[[t]]$ for each $\alpha \in \rts_{o, +}$ and $\bbsigma_i \in (\Omega_0[[t]])^*$ for $i \in I_o$. \end{nlem}

Note that under the topology on the coordinates $\Sigma$ that we have described at the end of \S \ref{subsub:iwahori-matsumoto-coordinates},  the subset  $\hU_{\Omega_0} \subset \hU$ is compact follow's from Tychonoff's theorem.

\tpoint{Reduction theory for $\hG$} \label{subsub:siegel-sets} Fix $t_0:= \frac{2}{\sqrt{3}}$ and $\Omega \supset \Omega_0$, the subset constructed in \eqref{omega-0}. The fundamental theorems of reduction theory proven by Garland assert the following. 
\begin{nthm} \label{thm:reduction-theory-1} Let $t < t_0,$ and $\Omega \supset \Omega_0$ and define the \textit{Siegel set} \be{} \label{siegel-set} \sie_{t, \Omega} :=  \hA^+_t \times \hU_{\Omega}. \ee
	
	\begin{enumerate} 
		\item \cite[Theorem 20.14]{gar:ihes}  For every $x \in \xg$, there exists $\gamma \in \hgam$ such taht $x \gamma \in \sie_{t, \Omega}$. 
		
		\item \cite[Thm. 21.16]{gar:ihes} For each $i \in I$, there exists $N_i:=N_i(t)$ such that if there exists \be{} \label{int:hypothesis} \begin{array}{lcr} x \in \sie_{t, \Omega} \gamma \cap \sie_{t, \Omega} & \mbox{ with } & \la a_i, H_B(x) \ra < N_i \end{array} \ee then in fact $\gamma \in \hGam \cap P_i$ where $P_i$ is the standard maximal parabolic corresponding to $I \setminus \{ i \}.$
		
		\item \cite[Lemma 21.10]{gar:ihes} There exists $0< s_0 < 1$, say $s_0 = e^{-r_0}$ for $r_0 >0$ so that if $x \in \sie_{t, \Omega}$ satisfies $\la \delta, H_B(x) \ra < - r_0,$  then there exists some $i \in I$ such that $\la a_i, H_B(x) \ra < N_i$. 
		
	\end{enumerate}
	
\end{nthm}

\begin{nrem}  \label{remark:after-reduction-theory}
	
	\begin{enumerate} 
		\item  For any subset $J \subset I$ and $j \in I$ we have $P_{I \setminus J } \cap P_{I \setminus \{ i \}} = P_{I \setminus \{ J \, \cup \{ i \} \} }.$ So we can also conclude that if $x \in \sie_{t, \Omega} \gamma \cap \sie_{t, \Omega}$ and $\la a_i, H_B(x) \ra < N_i$ for some subset $\psi \subset I$, then $\gamma \in P_{I \setminus \psi} \cap \hGam$.
		\item Examing the proof of (2) in \cite{gar:ihes} one notes that the hypothesis \eqref{int:hypothesis} can be replaced by \be{} \label{weaker-int-hyp} \begin{array}{lcr} x \in 
			\left( \hA^+_{t} \times \hU \right)  \gamma \,  \cap \, \left( \hA^+_{t} \times \hU \right)  & \mbox{ with } & \la a_i, H_B(x) \ra < N_i. \end{array} \ee
		\end{enumerate}
		
	\end{nrem}

\newcommand{\bsig}{\mathbf{\mf{S}}}

\section{Canonical pairs and partitions of $\xg$}  \label{sec:canonical-pairs}

\newcommand{\zt}{\mathtt{Z}}
\newcommand{\dt}{\mathtt{d}}

\subsection{Preliminaries}

In this subsection, we collect various preliminary results about or a Lie algebraic nature that will be used in the formulation of the main results in \S \ref{sub:main-thm-canonical}.

\newcommand{\csimp}{\check{\simp}}

\spoint As before, we write $\hB \subset \hG$ for the subgroup introduced in \S \ref{subsub:BorelKM} and corresponding to the Borel subset $\rts_+$. For each $P \in \stdp$, corresponding to the parabolic subset $\pb_{J}$ with $J \subset I,$ we have attached the subspace $\hfh^+(P):= \hfh(\es, J) \subset \hfh$ as well as the  sub\textit{set}s $(\hfh^+)_P:= (\hfh^+)_{J} \subset \hfh$ and $\hfh^+_P \subset \hfh^e$, see \eqref{def hatrJ} , \eqref{def:a_P}, and \eqref{hplus:d}. Let us define \be{} \label{root p} \begin{array}{lcr} \Delta_{\hB}^P:= \Delta(P) := \{ a_i|_{\hfh^+(P)}, i \in J \} & \text{and }& \Delta_P:= \Delta_P^{\hG} = \{ a_i|_{\hfh^+_P}, \,  i \notin J \} \end{array} \ee By construction each $a \in \Delta(P)$ restricts to zero on  $(\hfh^+)_P$. We view $\Delta_P$ as functionals on $\hfh^+_P$, not on $(\hfh^+)_P$ since it could be the case that for certain $J$, every $a \in \Delta_P$ is zero on $(\hfh^+)_{P}$, e.g. $J= I_o$.  More generally, if $P \subset Q$ are two standard parabolics, we write $\Delta_P^Q$ for the set of non-zero restrictions of elements from $\Delta_{\hB}^Q$ to the subspace $\hfh^+_P$.

\newcommand{\dhat}{\widehat{\Delta}}
\spoint Let  $\dhat_{\hB}^{\hG}$ denote the set of fundamental weights $\{ \lambda_i \}_{i \in I}.$ If $P= P_{J} \in \stdp'$  let \be{} \begin{array}{lcr} \dhat_{\hB}^P= \{ \lambda_i|_{\hfh^+(P)}, i \in J\} & \mbox{ and } & \dhat_{P} = \{ \lambda_i, i \notin J \} \end{array}. \ee 
\noindent  For $P \subsetneq \hG$, $\hfh^+(P)=\hfh(J)$  is the subalgebra from \S \ref{subsub:h(J)-coweights}. Defining $\omega^{J}_i, i \in J$ as in \textit{loc. cit.},  \be{} \lambda_i |_{\hfh^+(P)} = \omega^{J}_i|_{\hfh^+(P)}. \ee On the other hand, viewing each $\omega^{J}_i$ as a functional on all of $\hfh^+$ (using \eqref{omega_J}), we must have $\omega^{J}_i |_{\hfh^+_P} =0$ (and so in general, it is not equal to $\lambda_i$).

More generally, if $P \subset Q$ and $\sta{P} \subset M_Q$ (see \S \ref{subsub:relative-langlands-dec}) we have  from \eqref{dec:hrJ:K}  \be{} \label{h(Q):decP} \hfh^+(Q) = \hfh^+(P) + \hfh^+(Q)_{\sta{P}}. \ee  In this case, we denote by $\dhat_P^Q$ to be the set of elements $\omega \in \dhat_{\hB}^Q$ which are null on $\hfh^+(P)$. \textit{Example:}
\begin{itemize} 
\item if $P = P_{J}$ and $Q= P_{K}$ for $J \subset K \subsetneq I$, then $\dhat_P^Q  = \{ \omega^{K}_j \mid j \in K \setminus J \};$ and
\item  if $Q=\hG$, we just drop it from the notation and write $\dhat_P:= \dhat_P^{\hG}.$ It agrees with our previous notation. 
\end{itemize}

\tpoint{The functionals $\rho_P$} \label{subsub:rho-p-q}  For each $P \subsetneq \hG$ a standard parabolic, we define \be{} \label{rho-rho_P-rho(P)} \begin{array}{lcr}  \rho(P):= \sum_{\omega \in \dhat_{\hB}^P} \omega \mbox{ and } \rho_P:= \rho - \rho(P) \end{array}. \ee When $P = \hG$, we sometimes also set $\rho(\hG) = \rho$ and $\rho_{\hG}=0$. If $P \subsetneq \hG$ is standard define \be{} \dt_P(\alpha) = \la \rho(P), \alphav \ra \leq 0 \mbox{ for each } \alpha \in \Delta_P. \ee The inequality follows since $\rho(P)$ is a non-negative linear combinations of roots in $\Delta(P)$. As $\rho(P)$ is a linear combination of elements from $\Delta(P)$ which are zero on $\hfh^+_P$, we see that \be{} \rho_P = \rho|_{\hfh^+_P}. \ee We can obtain a more precise description of $\rho_P$ as follows.

\begin{nlem} \label{lem:rho-P-as-sum} Letting $\kappa_P(\alpha):= 1 - \dt_P(\alpha)$ for each $\alpha \in \Delta_P$ and $\dt_P(\dd) = \la \rho(P), \dd \ra,$  we have $\kappa_P(\alpha) \geq 0$ and $\dt_P \geq 0.$ As functionals on $\hfh$ and $\hfh^+$ we then have respectively \be{} \label{rho_P:no-d} \begin{array}{lcr} \rho_P= \sum_{\alpha \in \Delta_P} ( 1 - \dt_P(\alpha))  \lambda_{\alpha}  & \mbox{ and } &  \label{rho-dec-hplus} \rho_P =  \left(  \sum_{\alpha \in \Delta_P} \kappa_P(\alpha)  \lambda_{\alpha} \right) - \dt_P(\dd) \delta. \end{array} \ee 
\end{nlem} 
\begin{nrem} If $P$ is maximal, say $P:=P_i$ for $i \in I$, then we see that $\rho_P$ and $\lambda_i$ are positive, rational multiples of one another as functionals on $\hfh$ (they differ on $\R \dd$ though).  In this case, write: 
 \be{} \label{kappa-def} \rho_{P_i} =  \kappa_{i} \lambda_i - \dt_i \delta, \ee where we set $\kappa_i:= \kappa_{P}(\alpha)$ for $\alpha$ the unique element of $\Delta_{P_i}$ and $\dt_i= \dt_P$. When $i=\ell+1$, we have $\rho_{P_{\ell+1}} = \dcox \, \lambda_{\ell+1}$ since $d_{\ell+1}=0$, and the decomposition $\rho = \rho_{P_{\ell+1}}+ \rho(P_{\ell+1})$ is then just Lemma \ref{lem:rho-vs-rho-classical}.

 \end{nrem} 
\begin{proof}Suppose $P = P_{J}$ for $J \subsetneq I$. For (1), as  $\{ \av_i \}_{i \in I}$ is a basis of $\hfh$, we evaluate $\rho_P=\rho - \rho(P)$ on this basis. Since $\la \rho, \av_i \ra =1$ for all $i \in I$ and $\la \rho(P), \av_i \ra =1$ for all $i \in J$, the result follows. 

	Write $P=P_{J}$ with $J \subset I$. Note that $\rho(P)$ is a linear combination of elements from $\Delta(P)$, so that $\dt_P=0$ if $\ell+1 \notin J$. On the other hand, if $\ell+1 \in J$, then $\rho(P) = \sum_{i \in J} m_i a_i$ with each $m_i \geq 0$. Hence \be{} \label{d-Dpos} \dt_P(\dd):= \la \rho(P), \dd \ra = m_{\ell+1} \geq 0. \ee 
	As for \eqref{rho-dec-hplus}, we note that if $H -r \dd \in \hfh^+$ with $H  \in \hfh$, we  have that \be{} \la \rho_P, H - r \dd \ra = \la \rho_P, H \ra - \la \rho_P,  r \dd \ra =  \la \rho_P, H \ra - r \la \rho - \rho(P), \dd \ra =  \la \rho_P, H \ra + r \la  \rho(P),  \dd \ra .\ee

\end{proof}

\tpoint{The functionals $\rho_P^Q$} Suppose we are given $P, Q \in \stdp$ with $P \subset Q.$ Then recall from \S \ref{subsub:relative-variants} that $(\hfh^+)_P = (\hfh^+)_Q +  \hfh(Q)_{\sP}$. Let us now define \be{}  \label{rho-p-q:def} \rho_P^Q:=  \rho(Q)|_{\hfh^+(Q)_{\sP}}. \ee Note that if $Q= \hg$, we have set $\rho(\hg)= \rho$, so that $\rho_P^{\hg} = \rho_P.$

\begin{nclaim} On $\hfh$ we have $\rho_P^Q$ is a non-negative linear combination of $\dhat_P^Q$ and  on $(\hfh^+)_Q$, we have \be{} \label{rho:dec-PQ} \rho_P = \rho_Q+ \rho_P^Q. \ee \end{nclaim}

\subsection{Canonical Pairs and semi-stability} \label{sub:main-thm-canonical}

The aim of this subsection is to formulate the theorem of canonical pairs (see Theorem \ref{thm:canonical-pair}). To begin, we describe the related notion of semi-stability, following closely \cite[\S 2]{chau}.

\tpoint{Degree of instability} \label{subsub:new-deg-inst} Let $x \in \xg$. For $\hq \in \stdp'$ define the \textit{degree of $Q$-instability of $x$} in terms of the following minimization (which is shown to exist in  \S \ref{subsub:existence-of-minima})  \be{} \label{deg:Q-inst} \deg_{\inst}^{\hq}(x):= \min_{\hp \subset \hq, \delta \in  Q_{\Q}/ P_{\Q}} \la \rho^Q_P, H_{\sta{P}}(m_Q(x \delta))  \ra =  \min_{\hp \subset \hq, \delta \in  Q_{\Q}/ P_{\Q}} \la \rho^Q_P, H^+_P(x \delta)  \ra \ee where $\hp$ ranges over standard, proper parabolic subgroups of $Q$, and where we have written $m_Q(x \delta) = \pr_{M_Q}(x \delta)$. In the second expression above, we regard $\rho_P^Q$ as a functional on $\hfh(Q)$ and the equivalence of the two expressions is due to \eqref{eq:lem-H-PQ}. 
 In case $\hq=\hg$, we shall set \be{} \deg_{\inst}(x):= \deg_{\inst}^{\hg}(x) =  \min_{\hp \subset \hg, \delta \in  \hg_{\Q}/ P_{\Q}} \la \rho_P , H^+_{\hp}(x \delta) \ra. \ee Let us note here that for $z \in \xg$ and $P \subsetneq \hG$, we also have \be{} \la \rho_P, H_B^+(z) \ra = \la \rho_P, H_P^+(z) + H(P)(z) \ra = \la \rho_P, H_P^+(z) \ra \ee since in the expression \eqref{rho-dec-hplus} for $\rho_P$, the term $\delta$ is always zero on $H(P)$, and so is the sum over $\Delta_P$. On the other hand, we can also write \be{} \la \rho_P, H_P^+(z) \ra = \la \rho, H_P^+(z) \ra \ee since $\rho=\rho_P + \rho(P)$ and $\rho(P)$, a linear combination of elements from $\Delta(P)$, is zero on $H_P^+(z)$.

\newcommand{\zr}{\mathbf{c}_{\rho}}

\tpoint{On the invariance of $\deg^Q_{\inst}$} For $Q \subsetneq \hG$, the expression $\deg^Q_{\inst}(x)$ is invariant under right multiplication by $\hA_Q$, i.e. $ \deg^Q_{\inst}(x a_Q) = \deg^Q_{\inst}(x ) $ since for $P \subsetneq Q$, 
\be{}  \la \rho_P^Q, H_P^+(x a_Q \delta)  \ra = \la \rho_P^Q, H_P^+(x \delta a_Q)  \ra = \la \rho_P^Q, H_P^+(x  \delta) +H_{P}^+(a_Q)   \ra \ee and, since $\rho_P^Q$ is a linear combination of elements from $\Delta(Q)$ and $H_P^+(a_Q) \in \hfh_Q^+,$ the expression $\la \rho_P^Q, H_P^+(a_Q) \ra =0$. In particular, in particular, $\deg^Q_{\inst}$ is invariant under translations by the central subgroup (see \S \ref{subsub:central-T}) $\hA_{\cc} \subset \hA_Q$. On the other hand when $Q=\hG$, this no longer holds, i.e. $\deg_{\inst}(x)$ is \textit{not} invariant under the action of $\hA_{\cc}$. In fact, if $H_{\cc}(x)$ denotes the projection of $H_{B}(x)$ onto $\R \cc$,
 \be{} \label{deg-inst-not-invariant} \deg_{\inst}( \zeta x ) = \la \rho, H_{\cc}(\zeta) \ra + \deg_{\inst}(x). \ee

\tpoint{Semi-stability} For $Q \in \stdp$, we say that $x \in \xg$ is \textit{$Q$-semistable} if \be{} \label{Q-ss} \deg^Q_{\inst}(x) \geq 0. \ee  If $Q=\hg$, we just refer to this as semi-stability.  Note  that $x$ is semi-stable if and only if $x \gamma$ is semi-stable for $\gamma \in \hG_{\Q}$. 

\begin{nrem} As $\deg_{\inst}$ is not $\hA_{\cc}$-invariant, the condition of being semi-stable is not invariant under translations by $\hA_{\cc}$. \end{nrem}

Note also that every $x \in \xg$ is $\hB$-semi-stable as the condition $\deg^{\hB}_{\inst}(x) \geq 0$ is trivially satisfied.

\begin{nprop} (\textit{cf.} \cite[Lemma 2.2.1]{chau} ) \label{lem:semi-stability-equivalences} Let $x \in \xg$.
	
	 \begin{enumerate} 
	 	
	 	\item For $\hq \subsetneq \hG$ be a standard parabolic, the following  conditions are equivalent.  
				\begin{enumerate}
					\item $x \in \xg$ is $\hq$-semi-stable;
				
					\item for any standard parabolic $\hp \subset \hq$ and any $\omega \in \dhat^{Q}_{P}$ we have \be{} \la \omega, H_{\sP}(m_Q(x \delta)) \ra = \la \omega, H^+_P(x \delta)  \ra \geq 0 \ee for each $\delta \in Q_{\Q}/ P_{\Q};$ and
		
					\item for any standard, maximal parabolic $\hp \subset \hq$ and any $\omega \in \dhat^{Q}_{P}$ we have \be{} \la \omega, H_{\sP}(m_Q(x \delta)) \ra = \la \omega, H^+_P(x \delta)  \ra \geq 0 \ee  for each $\delta \in Q_{\Q}/ P_{\Q}.$
		
				\end{enumerate}
				
		\item When $\hq=\hG$, we have we have that (c) implies (b), and (b) implies (a). Moreover, $x$ is semi-stable if and only if, in the notation of \eqref{kappa-def}, \be{} \label{ss-ineq} \la \lambda_i, H_B(x \delta)  \ra \geq -r \frac{\dt_i}{\kappa_i}  \mbox{ for all } \delta \in \hG_{\Q}/B_{\Q} \mbox{ and all } i \in I. \ee
	
	\end{enumerate}
	
\end{nprop}
\begin{proof} Let us first argue that (a) implies (c).  Suppose $P \subset Q$ is maximal. Then on $\hfh(Q)$, and in particular on $\hfh(Q)_{\sP}$ we have that the unique element of $\dhat_P^Q$ is proportional to $\rho_P^Q$ by a positive (rational) number. So, condition (3) follows from (1).
	
	 Next we show that (c) implies (b). For $J' \subset J \subset I$, let $P:= P_{J'}, Q=P_{J}$ so that $P \subset Q$. Each $\omega \in \dhat_{P}^{Q}$ is of the form $\omega=\omega_i$ for one of the fundamental weights of $M_Q$ (with the convention that $M_{\hG}=\hG$) with $i \in J \setminus J'.$ 
	 Setting $P_i:=P_{J \setminus \{ i \} }$, we have that $P \subset P_i \subset Q$, and $P_i$ is a maximal parabolic in $P.$ We may regard $\omega_i$ as the unique element of $\dhat_{P_i}^Q$. We would like to argue that \be{}\la \omega_i, H_{\sP_i}(m_Q(x \delta )) \ra = \la \omega_i, H_{\sP}(m_Q(x \delta )) \ra \ee or equivalently that  \be{}\la \omega_i, H^+_{P_i}(x \delta ) \ra = \la \omega_i, H_P^+(x \delta ) \ra. \ee This follows by nothing that \be{} H^+_{P}(x \delta) = H^+_{P_i}(x \delta ) + H_{\sP(P_i)}(x \delta) \ee and also that $\omega_i$ is zero on $H_{\sP(P_i)}(x \delta) \subset \hfh(P, P_i)$.

	 Using the decomposition $(\hfh^+)_{P_i} = (\hfh^+)_Q +  \hfh(Q)_{\sP_i}$ together with the fact that $\omega_i$ is null on $\hfh(Q)_{\sP_i}$ we conclude
	   \be{} \la \omega_i, H_{P}(m_Q(x \delta'))  \ra  = \la \omega_i, H_{\sta{P}_i}(m_Q(x \delta)) \ra,\ee
	  where $\delta \in Q_{\Q}/P_{i, \Q}$ and $\delta'$ is any lift to $Q_{\Q}/P_{\Q}$ under the natural surjection $Q_{\Q}/P_{\Q} \rr Q_{\Q}/P_{i, \Q}.$ From these considerations, (2) follows from (3). 
	 
	  Finally to show that (b) implies (a), one just needs to use  that $\rho_P^Q$ is equal to the restriction of $\rho(Q)$, and hence a non-negative linear combination of elements from $\dhat_P^Q$.

	Let us now turn to part (2). First we note that the proof that (c) implies (b) follows as above. Assuming (b), we note from Lemma \ref{lem:rho-P-as-sum} (2) that $\rho_P$ can be written as positive linear combination of $\widehat{\Delta}_P$ as well as $ - \dt_P \delta$ with $\dt_P \leq 0$. Since we assumed that $r \geq 0$, the result follows. The statement (a) implies (c) does not however hold in general. Instead, if suppose $P \subset \hG$ is maximal, say $P=P_i$ for some $i \in I$. Then $\dhat_{P_i}^{\hG}$ consists of $\lambda_i$. Using \eqref{kappa-def} and the definition of semi-stability, 
	\be{} 0  \leq \la \rho_{P_i}, H_{B}^+(x \delta)  \ra = \la \kappa_i \lambda_i - \dt_i \delta, H_{B}^+(x \delta)  \ra  = \kappa_i \la \lambda_i, H_{B}(x \delta)  \ra + r \dt_i. \ee The stated inequality \eqref{ss-ineq} follows. The converse follows similarly using the fact, just observed, that (c) implies (a) even when $\hq=\hG$.

\tpoint{Theorem of canonical pairs} Consider a pair $(P, \delta)$ consisting of a standard parabolic $P \subsetneq \hG$ and $\delta \in \hg_{\Q}/ \hp_{\Q}.$ Such a $(P, \delta)$ is called a \textit{destabilizing pair} for $x \in \xg,$ if  \be{} \la \rho, H_P(x \delta) \ra = \deg_{\inst}(x). \ee Such a pair $(P, \delta)$ will be called  \textit{extremal} if for any standard parabolic $Q \supsetneq P$  \be{} \la \rho, H_Q(x \delta)  \ra > \la \rho, H_P(x \delta) \ra. \ee  In the above, we interpret $\delta$ as the image of the natural projection $\hG_{\Q}/P_{\Q} \rr \hG_{\Q}/Q_{\Q}$, and also allow $Q= \hG$ with the convention that $H_{\hG}(z)=0$ for any $z \in \xg$.

\begin{nthm} \label{thm:canonical-pair} For each $x \in \xg$ which is not semi-stable there is a unique destabilizing and extremal pair $(P, \delta)$ with $P$ a standard parabolic and $\delta \in \hG_{\Q}/ \hP_{\Q}.$  \end{nthm}

\begin{nrem} \label{rem:neg-destab} Note that if $x \in \xg$ is such that  $\deg_{\inst}(x) \geq 0$ (i.e. $x$ is semi-stable) then it cannot have a destabilizing, extremal pair $(P, \delta)$ with $P \subsetneq \hG$: indeed, pick any $(P, \delta)$ such that $\deg_{\inst}(x) = \la \rho_P, H_P^+(x \delta) \ra \geq 0$. But then $P \subset \hG$ and $0= \la \rho, H_{\hG}^+(x \delta) \ra \leq \la \rho_P, H_P^+(x \delta).$  \textit{The} destabilizing, extremal pair attached to a (non semi-stabe) $x \in \xg$ is called its  \textit{canonical pair} and we write \be{} \label{def:cp}  \cp_{\hG}(x):= \cp(x) = (P, \delta). \ee Abusing notation, if $x$ is semi-stable we also set $\cp_{\hG}(x)= (\hG, 1)$. Note however that in this case $\deg_{\inst}(x)$ need not agree with $\la \rho_G, H_G(x) \ra =0,$ 
 
 \end{nrem}

 The proof of the theorem breaks into two parts: the existence which follows from earlier work of H. Garland (see \S \ref{subsub:existence-of-minima}) and the uniqueness will be proven in in \S \ref{subsub:uniqueness}. For finite-dimensional groups, this result has been proven in different guises by a number of authors, \textit{e.g.} \cite{beh}, \cite{stuhler}, \cite{grayson}, \cite{harder-narasimhan}, and in a form most closely resembling the present work in \cite{chau}. We shall henceforth take this finite-dimensional result as granted (though its proof is \textit{not} used here, we shall need to quote the result later; its proof follows by adapting the proof given below in the affine case to the finite-dimensional one)

\tpoint{An auxiliary result}  We begin by observing (\textit{cf.} \cite[Lemma 2.3.2]{chau})

\begin{nlem}  \label{lem:can-pair-crit}Let $x \in \xg$ and $(P, \delta)$ a pair with $P\in \stdp'$ and $\delta \in \hG_{\Q}/ P_{\Q}.$
	\begin{enumerate}
		\item If $(P, \delta)$ is destabilizing for $x$, then $\deg_{\inst}^P(x \delta) \geq 0$, i.e. $x \delta$ is $P$-semistable.
		\item If $(P, \delta)$ is extremal for $x$, then 
		\begin{enumerate}
			\item $\la \alpha, H^+_P(x \delta) \ra < 0$ for all $\alpha \in \simp_P$
			\item $\la \rho_P, H^+_P(x \delta) \ra < 0$.
		\end{enumerate}

	\end{enumerate}
\end{nlem}

\begin{nrem} A converse for this result will be given as Corollary \ref{cor:can-pair-criterion}. Note that in finite type, if $P$ is a proper parabolic and $H \in \mf{h}_P$ is such that $\la \alpha, H \ra < 0$ for all $\alpha \in \Delta_P$. In fact, it is also true that $\la \alpha, H \ra \leq 0$ for all positive roots $\alpha$. Hence $\la \lambda_i, H \ra \leq 0$ for each fundamental weight, since each $\lambda_i$ is a positive linear combination of the simple roots. One can deduce from this that in fact $\la \rho_P, H \ra = \la \rho, H \ra < 0$, i.e. 2(b) is a consequence of 2(a) in finite type. \end{nrem}
	
\begin{proof}
Let $z:= x \delta$. We would like to argue that if $Q \subset P$ and $\eta \in P_{\Q}/Q_{\Q}$, we have $\la \rho_Q^P, H_{Q}(z \eta) \ra \geq 0.$ Note that from \eqref{rho:dec-PQ}, $\rho_Q^P = \rho_Q - \rho_P$ as functionals on $\hfh^+_Q = \hfh^+_P + \hfh(Q, P)$. Our assumption is that $(P, \delta)$ was destabilizing implies that \be{} \la \rho_P, H^+_P(z) \ra = \la \rho_P, H^+_P(z \eta) \ra \leq \la \rho_Q, H^+_Q(z  \eta) \ra. \ee Hence we have $\la \rho_Q^P, H^+_Q(z \eta) - H^+_P(z \eta) \ra \geq 0$ and this proves part (1).

As for part (2),  let $Q \supsetneq P$ is another standard parabolic and let $z:= x \delta$. If $(P, \delta)$ is extremal for $x$, then $\la \rho, H^+_P(z) \ra < \la \rho, H^+_Q(z) \ra$. It follows that that \be{} \label{HP:HR} \la \rho, H^+_P(z) - H^+_Q(z) \ra  = \la \rho_P^Q , H_{\sta{P}}(m_Q(z)) \ra = \la \rho_P^Q, H^+_P(z) \ra < 0. \ee Suppose first that $Q=\hG$, then as $H_Q( z) =0$, we must have $\la \rho, H^+_P(z) \ra < 0$, which is of course just the condition that $\la \rho_P, H^+_P(z) \ra < 0$. If that $Q= \hG$ is the only parabolic properly containing $P$, then $P$ must be a maximal parabolic, say $P=P_i$ and $\Delta_{P}= \{ a_i \}.$ We shall contend that \be{} \la a_i, Y \ra < 0 \mbox{  for any } Y \in \hfh^+_P. \ee If $i=\ell+1$, then this is clear as $Y = m \cc - r \dd$ for $m \in \R, r >0$; so $\la a_i, Y \ra = -r < 0$. So, now suppose $i \neq \ell+1.$ Since $Y \in \hfh^+_{P_i}$, we have $\la a_j, Y \ra =0$ for $j \neq i.$ Hence from \eqref{def:delta},  \be{} 0 = \la a_{\ell+1}, Y \ra = \la - \vartheta + \delta, Y \ra  = - d_i \la a_i, Y \ra - r, \ee with $d_i >0.$  It follows that $\la a_i, Y \ra = -r/d_i < 0$.

It remains to consider the case that $P$ is not maximal. By the argument above, we again deduce condition 2(b). As for 2(a), pick $\alpha \in \Delta_P$ and let $R$ be the parabolic such that $P \subset R \subsetneq \hG$ and $\{ \alpha \} = \Delta_P^R.$ From \eqref{HP:HR},  we have that $\la \rho_P^R, H_P^+(z) \ra < 0$. On the other hand, $\rho_P^R,$ which is a non-negative linear combination of the simple roots in $\Delta(R)$, is actually just positive (rational) multiple of $\alpha$ on $\hfh^+_P$. This proves condition 2(a).

\end{proof}

\subsection{Existence and Uniqueness of a canonical pair} \label{sub:uniqueness}

\tpoint{Existence of minima}  \label{subsub:existence-of-minima} The existence of a canonical pair will follow from the following result due to Garland. 

\begin{nprop} \cite[Lemma 2.5]{gar:ms-2} \label{prop:existence-of-minima} For $0 < s < 1$, let $x \in \hg \eta(s).$ For $P \subsetneq \hG$ a standard parabolic and any $M > 0$, there exists a finite subset $\Xi:= \Xi(x,P) \subset \hG_{\Q} / \hP_{\Q}$ so that for all $\delta \notin \Xi,$ \be{}  \la \rho_P, H_P^+(x \delta)  \ra = \la \rho, H_P(x \delta) \ra \geq M. \ee \end{nprop}
\begin{proof} The proof is the same as that of \cite[Lemma 2.5]{gar:ms-2}, but let us make a few small remarks so that the result stated in \textit{op. cit.} can be applied. First, we note that without loss of generality we can choose $x \in \mf{S}_t$ since we can always pick $ \gamma \in \hG_{\Q}$ such that $z:=x \gamma \in \mf{S}_t$. If the Proposition holds for $z$, then it certainly also holds for $x$. 
	
Suppose $P= P_{J}$ with $J \subset I$. Then we can apply \cite[Lemma 2.5]{gar:ms-2} to the dominant, integral weight $\rho_P = \sum_{i \in J} \kappa_i \lambda_i - \dt_P \delta$. Note that 
the set of $j \in I$ such that $\la \rho_P, \av_j \ra =0$ is exactly equal to $J$, so that the parabolic attached to $\rho_P$ is again $P$. 

Finally we remark that since we are in a class number one situation, $\hG_{\Q}/P_{\Q}$ is the same as $\hGam/ \hGam \cap P$. With these remarks, the Proposition follows from the aforementioned result of Garland.  \end{proof}

	 Applying the Proposition to $x \in \xg$, we conclude that $\deg_{\inst}(x)$ exists since there are finitely many standard parabolics, and for each, we are reduced to a minimum over a finite set by the Proposition.

\tpoint{An auxiliary inequality} \label{subsub:schieder-ineq} To begin the proof of uniqueness of canonical pairs, we begin with the following inequality inspired by \cite[Thm. 4.2]{schieder}. 

\begin{nprop} \label{prop:w-ineq} Fix $x \in \xg.$ Suppose $P_1  \subsetneq \hG$ is a standard parabolic and $\delta_1 \in \hG_{\Q}/P_{\Q}$ is such that $\deg_{\inst}^{P_1}(x \delta_1) \geq 0$ and $\la \alpha, H^+_P(x \delta) \ra < 0$ for all $\alpha \in \simp_P$. Let $(P_2, \delta_2)$ be an another  pair with $P_2 \subsetneq \hg$ a standard parabolic and $\delta_2 \in \hG_{\Q}/\hP_2(\Q).$ Pick $w \in \affW$ a minimal length representative in $\affW(P_1) \setminus \affW / \affW(P_2)$  so that $\delta_1^{-1} \delta_2 \in P_1 \dw P_2.$ Then we have, \be{} \label{schieder-ineq} H^+_{P_1}(x \delta_1) \dleq w^{-1} \cdot \, H^+_{P_1}(x \delta_1 ) \dleq H^+_{P_2}(x \delta_2). \ee If both of the above are equalities, then $w=1.$

		  \end{nprop} 

\begin{proof}
	
To begin assume that $P \subsetneq G$ is a standard parabolic and $\delta \in \hG_{\Q}/P_{\Q}$. Assume $(P, \delta)$ satisfies the condition that $\la \alpha, H^+_P(x \delta) \ra < 0$ for all $\alpha \in \simp_P$.  By definition, we also have $\la \alpha, H^+_P(x \delta) \ra =0$ for $\alpha \in \Delta(P)$. Hence $H:=H_P^+(x \delta)$ is a anti-dominant element of $\hfh$, and as so  $H \dleq w H$ for any $w \in \affW$ with equality if and only if $w \in \affW(P).$ Applying this to the setup of the Lemma gives us the first inequality in \eqref{schieder-ineq}, so let us turn to the second. 

Let us write $P_1:= P_{J_1}, P_2:=P_{J_2} \in \stdp$ with $J_1, J_2 \subsetneq I$. It is enough to verify the inequality    \be{} \label{w-i:1} I(\lambda): \ \ \  \la \lambda, w^{-1} \cdot \, H_{P_1}(x \delta_1 ) \ra  \leq \la \lambda, H_{P_2}(x \delta_2) \ra \ee for  every fundamental weight $ \lambda_d,  d \in I$.  Fix $d \in I$ and use Lemma \ref{lem:ragh} to write $\lambda_d =  \mu + \tau$ with $\tau$ a positive sum of the roots $a_j$ with $j \in J_2$ and $\mu \in (\hfe)^*$ satisfying \be{} \label{mu-ragh-prop} \begin{array}{lcr} \la \mu, \av_i \ra \geq 0, i \in I & \mbox{ and } & \mu|_{\hfh^+(P_2)}=0. \end{array} \ee It suffices to verify the inequality $I(\lambda)$ \eqref{w-i:1} holds separately for $\lambda=\mu$ and $\lambda =\tau.$

First, suppose that $\lambda = \tau$. In this case, the right hand side of \eqref{w-i:1} is $0$ as $\tau$ is in the span of $a_j$ with $j \in J_2$. By the minimality of $w$ in the coset space $\affW(P_1) \setminus \affW / \affW(P_2)$ we have \be{} \label{kost-ineq}  \begin{array}{lcr} w^{-1}(\av_i) > 0, \ \ i \in J_1, & \text{ and } & w(\av_i) > 0, i \in J_2 \end{array}. \ee  Hence $w\tau \geq 0$ since $\tau$ is a \textit{positive} sum of $a_j$ with $j \in J_2$.  The desired  inequality follows since $\la a_i, H^+_{P_1}(x \delta_1) \ra <0$ for $i \in J_1$ and hence also $\la a_i, H^+_{P_1}(x \delta_1) \ra \leq 0$ for all $i \in I$. 

Next, we suppose that $\lambda = \mu$.  We  begin with the following simple observation.

\begin{nclaim} Let $J \subsetneq I$, $P:= P_J$, and let $\sigma \in \affW$ be such that $\sigma(\av_i) > 0$ for $i \in J.$ If $z \in \xg$ is such that $\deg_{\inst}^{P}(z) \geq 0$, then we have $\sigma \cdot H_P(z) \dleq \sigma \cdot H_B(z).$ \end{nclaim}  

\begin{proof}[Proof of Claim] Note that $Y:= H^+_B(z) - H^+_P(z) \in \hfh^+(P)$ . So it suffices to show that $\sigma (Y) \dgeq 0$, or, by virtue of our assumption on $\sigma$ that $Y \dgeq 0$. As $ \deg_{\inst}^{P}(z) \geq 0$, \textit{i.e.} $z$ is $P$--semi-stable, from Lemma \ref{lem:semi-stability-equivalences}, we have \be{} \begin{array}{lcr} \la \omega , H^+_B(z) \ra \geq 0 & \mbox{ for all } & \omega \in \dhat_{B}^P. \end{array} \ee Hence, $\la \omega, H^+_B(z) - H^+_P(z) \ra \geq0$ since $\la \omega, H^+_P(z) \ra =0$ (as $M_P$ is a finite-type root system, $\omega$ is a linear combination of $\alpha \in \simp(P)$). \end{proof}

Now, by  \eqref{mu-ragh-prop}, we have $\la \mu, H^+_{P_2}(x \delta_2) \ra = \la \mu, H^+_{B}(x \delta_2) \ra$ so we just need to show \be{} \label{w-i:2} \la \mu, w^{-1} \cdot \, H_{P_1}(x \delta_1 ) \ra   \leq  \la\mu, H_{B}(x \delta_2) \ra. \ee   Applying the claim to $J = J_1,$ $P=P_1$, $z = x \delta$, and $\sigma = w$ we conclude \be{}  w^{-1} H^+_{P_1}( x \delta_1) \dleq w^{-1} H^+_B(x \delta_1). \ee From Corollary \ref{cor:orthogonal-ineq} we have  $w^{-1} H^+_{B}(x \delta_1 w w^{-1}) \dleq H^+_{B}(x \delta_2)$ so that in sum  \be{} w^{-1} \, H^+_{P_1}(x \delta_1) \dleq H^+_{B}(x \delta_2).  \ee Pairing both sides with $\mu$, the inequality $I(\mu)$ of \eqref{w-i:2} follows.

Finally note that if both inequalities in \eqref{schieder-ineq} are actually equalities, then $w^{-1} \in \affW(P_1)$ by part (1). But as $w$ was chosen to be a minimal length represenative, we must have $w=1.$

\end{proof}

\tpoint{Proof of unicity in Theorem \ref{thm:canonical-pair}.}  \label{subsub:uniqueness} Suppose $x \in \xg$ and $(P_1, \delta_1)$ and $(P_2, \delta_2)$ are two destabilizing, extremal pairs for $x$. As was already remarked, to exist we need $\deg_{\inst}(x) < 0$.

Let $w \in \affW$ as in Proposition \ref{prop:w-ineq} so that $\delta_1^{-1} \delta_2 \in P_1 \dw P_2.$
As $(P_1, \delta_1)$ satisfies the conditions of part (2) of Lemma \ref{prop:w-ineq}, we may pair the string of inequalities \eqref{schieder-ineq} with $\rho$ to obtain  \be{} \label{pf:unicity-1}  \la \rho, H^+_{P_1}(x \delta_1) \ra  \leq \la \rho, w^{-1} \cdot \, H^+_{P_1}(x \delta_1 ) \ra  \leq \la \rho, H^+_{P_2}(x \delta_2)\ra. \ee  Since $(P_2, \delta_2)$ is also destabilizing, we must have equality in \eqref{pf:unicity-1}. Observe that $w^{-1}H_{P_1}(x \delta_1) = H_{P_1}(x \delta_1) - \betav_+$ for $\betav_+$ a non-negative sum of coroots. Then by the regularity of $\rho$ and the minimality of $w$, we can conclude that  $w=1.$ Hence, $H^+_{P_1}(x \delta_1) = H^+_{P_2}(x \delta_1)$. Pairing the left hand side with $a_i$ for $i \in J_1$ we obtain zero, and hence also $\la a_i, H^+_{P_2}(x \delta_1) \ra =0$. This forces $i \in J_2$ and so $P_1 \subset P_2$. By the extremal property of $(P_1, \delta_1)$ we must then have $P_1=P_2$, and the result follows from this. \end{proof}

\spoint We can also state the following converse to Proposition \ref{lem:can-pair-crit}.

\begin{ncor} \label{cor:can-pair-criterion} Let $x \in \xg$. Suppose $P \subsetneq \hG$ is a standard parabolic and $\delta \in \hG_{\Q}/P_{\Q}$ are such that: \begin{enumerate} \item $\deg_{\inst}^P(x \delta) \geq 0$; and \item   $\la \alpha, H_P(x \delta) \ra < 0$ for $\alpha \in \simp_P,$ and $\la \rho_P, H^+_P(x \delta) \ra < 0.$ \end{enumerate} Then $(P, \delta)$ is the canonical pair for $x$.  \end{ncor}
\begin{proof} Suppose $(P_2, \delta_2)$ is the canonical pair for $x$ with $P_2 \subsetneq \hG$. We can now apply Proposition \ref{prop:w-ineq}(2) and argue as in \S \ref{subsub:uniqueness} to conclude that $(P_2, \delta_2)=(P, \delta)$.\end{proof}
	
\newcommand{\wot}{\leftidx_{{P_1}} W_{P_2}}

\subsection{A partition of $\xg$} \label{subsub:partition-X}

Using the theory of canonical pairs (see Theorem \ref{thm:canonical-pair}) we obtain a partition of $\xg$ indexed by the elements in $\stdp$ that we will now study.

\tpoint{Definition of partition}   For $\hP  \subsetneq \hG$ a standard parabolic we define \be{} \qxg(P) =\{x \in \xg | \, \cp(x) = (P, \delta)  \text{ for some } \delta \in \hG_{\Q} / \hP_{\Q} \}, \ee and if are given $(P, \delta)$ as above, we also set  \be{} \xg(P, \delta) :=  \{x \in \xg | \, \cp(x) = (P, \delta)  \}. \ee When $\hp=\hg$, we shall just write  $\xgss:= \xg(\hG)$ for the set of semi-stable elements in $\xg$. 

\begin{nrem}\label{notation:X-G-S} To a pair $(P, \delta)$ with $P \in \stdp'$ and $\delta \in \hG_{\Q}/P_{\Q}$ we can attach a rational parabolic subgroup $S = P^{\delta} = \delta P \delta^{-1}$ for (any) choice of representative in $\hG_{\Q}$ of $\delta$ which we have written with the same name. Moreover, to each rational parabolic $S \in \ratp'$ conjugate to $P$, we can attach a pair $(P, \delta)$ with $\delta \in \hG_{\Q}/P_{\Q}$ well-defined. For this reason, we sometimes write $\xg(S)$ for $\xg(P, \delta)$. \end{nrem}

\begin{nprop} \label{prop:basic-partition} We have a disjoint union \be{} \xg &=& \sqcup_{\hP \in \stdp} \, \xg(\hp) = \xgss \,  \sqcup_{\hP \in \stdp'} \, \xg(\hp) \\
	&=& \xgss \,  \sqcup_{\hP \in \stdp', \zeta \in \hG_{\Q}/P_{\Q}} \, \xg(\hp, \zeta).  \ee For $P \in \stdp$, $\qxg(P)$ is stable under right multiplication by $\hg_{\Q}$. In fact, we have \be{} \label{cp-inv} \begin{array}{lcr} \cp(x)= (P, \delta)& \text{ if and only if } &  \cp(x \eta) = (P, \eta^{-1} \delta)  \end{array} \ee for  $x \in \xg,$ $P\in \stdp$, $\delta \in \hG_{\Q}/P_{\Q},$ and $\eta \in \hG_{\Q}.$ 
\end{nprop}
\begin{proof} That $\xg$ is equal to such a disjoint union follows immediately from Theorem \ref{thm:canonical-pair}. As for the second claim, if $P \subsetneq \hG$, suppose $x \in \xg(\hP)$ and $\eta \in \hG_{\Q}$. Then for any $Q \in \stdp$,   \be{} \min_{\gamma \in \hG_{\Q}/ Q_{\Q}} \la \rho, H_Q(x \gamma)   \ra = \min_{\gamma \in \hG_{\Q}/ Q_{\Q}} \la \rho, H_Q(x \eta \gamma) \ra. \ee Hence if $\cp(x)= (P, \delta)$ then $\cp(x \eta) = (P, \eta^{-1} \delta)$ so that $x \eta \in \xg(\hP)$ as well. The stability under the right $\hG_{\Q}$-action also follows similarly for $\xgss$. \end{proof}

	\newcommand{\xm}{X_{M_P}}
	\newcommand{\xmss}{X_{M_P}^{ss}}
	\newcommand{\stdpm}{\mathrm{StdPar}_{\Q}(M_P)}
	\newcommand{\sr}{\leftidx^*R}
	\newcommand{\ozeta}{\overline{\zeta}}
	
	\tpoint{Canonical pairs for the Levi} \label{subsub:R-sR} \label{subsub:R-P-notation} Recall that for a standard parabolic $P \subsetneq \hG$, we have the projection $\pr_{M_P} : \xg \longrightarrow X_{M_P}$ stemming from the $P$-horospherical decomposition.  Now $M_P$ is a finite-dimensional algebraic group with Borel subgroup $\sta{B}:= \hB \cap M_P,$ and the parabolics containing $\sta{B}$ will be denoted as  $\stdp(M_P)$. They correspond to standard parabolics $R \subset P$ (see \S \ref{subsub:relative-langlands-dec}) with the correspondence denoted as  \be{}  \sta{R} \in \stdpm    \stackrel{1:1}{\longleftrightarrow}  \{ R \in \stdp \mid   \hr  \subset \hp \}. \ee 	One may check (as in the finite-dimensional case) that $\pr_{M_P}$ induces a bijection \be{} \pr_{M_P}: \label{zeta:o} P_{\Q} / R_{\Q} \rr M_{P, \Q} /\sta{R}_{\Q}, \zeta \mapsto \ozeta, \ee so the element $\ozeta$ will often be identified with an element $\zeta \in P_{\Q}$.

	The group $M_P$ also has a theory of semi-stability, canonical pairs, etc. \footnote{This is worked out in \cite{chau} for $GL_n$, but the same proof works for general finite dimensional groups. Alternatively, one can just adapt our argument for Theorem \ref{thm:canonical-pair} to finite-type root systems}. As before, let $\xmss := \xmss(M_P)$ denote the semi-stable locus, and to each non-semi-stable $z \in X_{M_P},$ attach \be{} \label{canonical-pair-M} \cp_{M_P}(z) = ( \sta{R}, \ozeta)= (\sta{R}, \zeta), \ee where $\sta{R} \subset M_P$, for $R \subset P$ a standard parabolic and $\ozeta \in M_P(\Q)/ \sta{R}_{\Q}$ which is identified with $\zeta \in P_{\Q}/R_{\Q}$. Write $\cp(z)=(M_P, 1)$ is $z \in X_{M_P}$ is semi-stable. Then the analogue of Proposition \ref{prop:basic-partition} for $M_P$ results in a decomposition \be{} \label{MP:partition} \xm = \underset{^*R \in \stdpm}{\sqcup} \, \xm(^*R) = \underset{^*R}{\sqcup}\, \underset{\zeta}{\sqcup}\, \xm(\sta{R}; \zeta) \ee where we define $ \xm(\sta{R}) := \{ z \in \xm \mid \cp_{M_P}(z) = (\sta{R}, \zeta) \mbox{ for some } \zeta \in M_{P, \Q}/ \sta{R}_{\Q} \}.$

	\tpoint{$P$-semistablity vs. $M_P$-semistability.} Keep the notation from \S \ref{subsub:R-P-notation}.

	\begin{nprop} \label{prop:lower-rank-semistable}  \label{can pair P,Id surject} Let $P \in \stdp'$, $z \in \xg$, and let $m:= \pr_{M_P}(z) \in X_{M_P}.$ Then  \be{}  \deg^{P}_{ \inst}(z) =  \deg^{M_P}_{\inst}(m) . \ee In particular, $z$ is $P$-semistable if and only if $m$ is $M_P$-semistable. \end{nprop}
	 
	 \begin{proof} The claim essentially follows from the definitions: to define $\deg^{M_P}_{\inst}$ one takes a minimum over all (standard) parabolics $\sta{R} \subset M_P$ and also over all points in $M_{P, \Q}/\sta{R}_{\Q}$. Using  \eqref{zeta:o} and Proposition \ref{subsub:relative-langlands-dec}, this is the same set over which we take the minimum defining $\deg_{\inst}^P$. Next we observe  using \eqref{lem:H-PQ} that if $z \in \xg$, $\zeta \in P_{\Q}/R_{\Q}$, $m:= \pr_{M_P}(z)$ and $\overline{\zeta}$ as in \eqref{zeta:o}, 
	 	\be{} 	\label{hr:*r,p} H^+_R(z\zeta) = H^+_P(z\zeta) +H_{^*R}(\pr_{M_P}(z \zeta) ) = H^+_P(z ) +H_{^*R} (m \overline{\zeta}). \ee  From this, we can deduce that $\la \rho_R^P, H_{\sta{R}}(m \overline{\zeta}) \ra = \la \rho_R^P, H_R( m \zeta) \ra .$ 
\end{proof}

	\newcommand{\hAres}[1]{\hA^+_{#1, (0, 1)}}

\tpoint{The basic fibration} For $z \in \xg$ and $\zeta \in \hG_{\Q}/P_{\Q}$, we have from \eqref{cp-inv} that $\xg(P, \zeta) \zeta= \xg(P, e).$ Recalling also that $\qxg(P):= \sqcup_{\zeta \in \hG_{\Q}/P_{\Q}} \, \xg(P, \zeta)$, we use the previous Lemma to define the map $\pi: \qxg(P) \rr X_{M_P}^{ss}$ which sends $z \in \xg(P, \zeta)$ to $\pr_{M_P}(z \zeta)$. This map is surjective as we can see from Corollary \ref{cor:can-pair-criterion}. Indeed if $m \in X_{M_P}^{ss}$ and $h \in \hA^+_P$ is such that $\la \alpha, H^+_P(h) \ra <0$ for all $\alpha \in \simp_P$ and $\la \rho_P, H_P(h) \ra < 0$, then $x=mh$ lies in $\xg(P, e)$ and certainly projects onto $m$. Such $h$ exist for any $P \subsetneq \hG$ so the surjectivity follows. 

	So, for each $P \subset \hG$ a standard parabolic and $t > 0$, we define \be{} \label{AP:01} \hA^+_{P, (0, t)}:= \{ a \in \hA^+_P \mid \la \alpha, H^+_P(a) \ra < \log t \mbox{ and } \la \rho_P, H^+_P(a) \ra < 0 \}. \ee The argument from the previous paragraph shows  \be{} \label{XGP:MPss} \xg(P, e) = X_{M_P}^{ss} \times \hAres{P} \times \hU_P. \ee If however $S \in \ratp'$, the natural structure present on $\xg(S)=X_G(P, \zeta)$ in the notation of Remark \S \ref{notation:X-G-S} is not that of a product, but that of a fibration  in terms of the natural projection from \eqref{split glob P}. \be{} \label{fib:S} \begin{array}{lcr} f_S: \xg(S)= X_G(P, \zeta) \rr  X^{ss}_{M_P} \times \hA^+_{P, (0, 1)} & \mbox{ which sends } & x \mapsto (\pr_{M_P}(x \zeta), \pr_{\hA^+_P}(x \zeta) ) \end{array}. \ee The fibers are equal to $\zeta^{-1} \hU_P\zeta$, so homeomorphic to $\hU_P$. Sometimes we put all of the $f_S$ together to obtain the natural map $f_P,$ \be{} \label{def:f_P} \begin{tikzcd}
		  \qxg(P):= \sqcup_{\zeta} \, \xg(P, \zeta) \arrow[dr, "f_P"] \arrow[r, "z \mapsto z \zeta"] & \xg(P, e) \ar[d]  \\ & X^{ss}_{M_P} \times \hA^+_{P, (0, 1)}   \end{tikzcd}  \ee
 	
	\spoint Suppose now that $P, R \in \stdp$, $R \subsetneq P \subsetneq \hG.$ Consider the  diagram \be{} \label{diag:P-R-proj}  \begin{tikzcd} \xg \ar[rr, bend right, swap, "\pr_{M_R}" ] \ar[r, "\pr_{M_P}"] &  X_{M_P} \ar[r, "\pr_{M_{^*R}}"] &  X_{M_R}=X_{M_{^*R}}, \end{tikzcd}. \ee and recall the partition of $X_{M_P}$ from \eqref{MP:partition}. Suppose that we are given $z \in \xg$ such that $m:=\pr_{M_P}(z) \in X_{M_P}(\sta{R}, \ov{\zeta})$. Then note that $H_{\sta{R}}(m \ozeta)$ is the closure of the anti-dominant cone of $\hfh(P)$, i.e. \be{} \label{H-star:anti-dom} \la \alpha, H_{\sta{R}}(m \ozeta) \ra \leq 0 \mbox{ for } \alpha \in \Delta(P). \ee Indeed, the above expression is zero for $\alpha \in \Delta(R) \subset \Delta(P)$ and is strictly negative for $\alpha \in \Delta_R^P$ by the finite-dimensional analogue of Corollary \ref{cor:can-pair-criterion} (where the condition involving $\rho_P$ is redundant). Writing $H_{\sta{R}}(m \ozeta)$ as a non-positive linear combinations of coweights of $M_P$, and hence coroots of $M_P$, it follows that $\la \beta, H_{\sta{R}}(m \ozeta) \ra \geq 0$ for $\beta \in \Delta_P$.

	\begin{nprop} \label{Pr can rel} Suppose that we are given $z \in \xg$ such that $m:=\pr_{M_P}(z) \in X_{M_P}(\sta{R}, \ov{\zeta})$.   If 
	 \be{}  \label{delta-P-assumption} \la \alpha, H_P^+(z) \ra &<&  - \la \alpha, H_{\sta{R}}(m \ozeta)  \ra \mbox{ for all } \alpha \in \Delta_P, \mbox{ and } \\  \la \rho_P, H^+_P(z) \ra &< &0, \ee then $z \in \xg(R, \zeta)$ where, as in \eqref{zeta:o}, $\zeta \in P_{\Q}/R_{\Q}$ the lift of $\overline{\zeta}$.

 \end{nprop}

	\begin{proof} We have  $\deg_{\inst}^{\sta{R}}(m \ozeta) \geq 0$ and so, from Proposition \ref{prop:lower-rank-semistable} (or rather its finite-dimensional analogue), we may conclude that $\pr_{M_{^*R}}(m \ozeta) \in X_{M_R}^{ss}=X_{M_{\sta{R}}}^{ss}$. Using Proposition \ref{prop:lower-rank-semistable} as well as the commutativity of the diagram \eqref{diag:P-R-proj}, we have that $z \zeta$ is $R$-semistable.  To see that $z \zeta \in \xg(R)$, from Corollary \ref{cor:can-pair-criterion}, we need that 
		\be{} \label{what-to-verify} \begin{array}{lcccr} \la \alpha, H^+_R(z\zeta) \ra < 0 & \text{ for all } & \alpha \in \Delta_R & \mbox{ and } \la \rho_R, H^+_R(z \zeta) \ra < 0. \end{array} \ee

Suppose first that $\beta \in \simp_R^P:= \simp_R\setminus \simp_P \subset \simp(P)$. Then $\la \beta, H^+_P(z) \ra =0$ and so from \eqref{hr:*r,p}, \be{} \la \beta, H^+_R(z \zeta) \ra =  \la \beta, H_{\sr}(m \overline{\zeta}) \ra < 0.  \ee But as $m\in X_{M_P}(\sta{R}, \ozeta)$, from the finite-dimensional analogue of Proposition \ref{lem:can-pair-crit},  \be{} \label{beta:ineq} \la \beta, H_{\sr}(m \overline{\zeta})  \ra < 0 \mbox{ for } \beta \in \Delta^{M_P}_{\sr}.\ee Identifying $\Delta^{M_P}_{\sr}$ with $\Delta^{P}_{R},$ we obtain \be{} \label{H:neg} \la \beta, H^+_R(z \zeta) \ra =  \la \beta, H_{\sr}(m \overline{\zeta}) \ra < 0  \ee where in the left hand side $\beta$ is regarded as an element of $\Delta_R^P$. 
			
Next, let $\alpha \in \simp_R \setminus \simp_R^P = \Delta_P,$ and consider 
			\be{} 	\la \alpha, H^+_R(z\zeta) \ra = \la \alpha,  H^+_P(z) \ra  +  \la \alpha, H_{^*R} (m \overline{\zeta}) \ra.	\ee By our hypothesis \eqref{delta-P-assumption}, it follows that this expression is again less than zero.

Finally, we next need to argue that \be{} \la \rho_R, H^+_R(z \zeta) \ra = \la \rho, H^+_R(z \zeta) = \la \rho, H^+_P(z) \ra + \la \rho, H_{\sta{R}}(m \ozeta) \ra \ee is less than zero. By assumption, $ \la \rho, H^+_P(z) \ra < 0$. As we remarked before the proof, $H_{\sta{R}}(m \ozeta)$ lies in the closure of the anti-dominant cone of $\hfh(P)$, hence $\la \rho, H_{\sta{R}}(m \ozeta) \ra=   \la \rho(P), H_{\sta{R}}(m \ozeta) \ra \leq 0 .$ \end{proof}

\spoint The following result resembles the theory of parabolic transformations (see Theorem \ref{thm:reduction-theory-1})(2)-(3), but unlike those results there is no restriction on `how deep' one needs to go into the parabolic ends, i.e. how small $\la a_i, H_B(x) \ra$ needs to be. 

\begin{nprop}  \label{prop:can-pairs-parabolic-trans} For $P \in \stdp'$ , the inclusion $i: \xg(P, e) \hookrightarrow \xg(P)$ induces \be{} \label{bijection:can-pairs-parabolic-trans}   \xg(P, e) / P_{\Q} \stackrel{1:1}{\longrightarrow}  \qxg(P)/\hG_{\Q}. \ee \end{nprop}

\begin{proof} Surjectivity follows since if $x \in \xg(P)$ with $\cp(x)= (P, \delta)$ then $\cp(x \delta)= (P, e)$. So $x \delta\in \xg(P, e)$ and $i( x \delta) = x$ in  $\qxg(P)/\hG_{\Q}.$ As for the injectivity, it follows from,

	\begin{nclaim} \label{lem:can-pair-transf}  If $\cp(x)=\cp(y)= (P, \eta)$ for some $\eta \in \hG_{\Q}$, which is well-defined by right translation by $\hP_{\Q}$,  and $x = y \gamma$ for some $\gamma \in \hG_{\Q}$, then  $\eta^{-1} \gamma \eta  \in \hP_{\Q}.$ In particular, if $\eta=1$, then $\gamma \in P_{\Q}$.  \end{nclaim}   
	
	As for the proof of the claim: the relation $x = y \gamma$ forces $\cp(x) = (P, \gamma^{-1} \eta)$. By the uniqueness of canonical pairs,  $\gamma^{-1} \eta = \eta \hP_{\Q}$, \textit{i.e.} $\eta^{-1} \gamma^{-1} \eta \in \hP_{\Q}$.   \end{proof}

\begin{nrem} We often use the following consequence of the above statement: for $\gamma \in \hG_{\Q}$, \be{} \xg(P,e) \gamma \cap \xg(P,e) \neq \es \mbox{ if and only if } \gamma \in P_{\Q}. \ee If $\gamma \in P_{\Q}$, then in fact $\xg(P, e) \gamma = \xg(P, e)$. \end{nrem}

\begin{ncor} \label{X/Gam:canonical-partition} We have a disjoint union \be{}  \xg / \hGam  =  \xgss/ \hGam \ \  \cup \ \  \sqcup_{P\in \stdp'} \,  \xg(P, e) / \hGam_P. \ee   \end{ncor}

	\tpoint{On the ends $\xg(P, e)/\hgam_P$} Recall the description for $\xg(P, e)$ in \eqref{XGP:MPss} which gives consider the natural projection  $\xg(P, e) \rr X_{M_P}^{ss} \times \hAres{P}$ sending $x \mapsto \left(\pr_{M_P}(x), \pr_{\hA^+_P}(x)\right).$
	
	\begin{nprop} \label{xgpgamma des} The projection above descends to a continuous surjection, again denoted by $f_P$,  \be{} \label{def:f_P-quotient} \begin{tikzcd}  \xg(P, e) / \hgam_P \ar[r, "f_P" ] &  X_{M_P}^{ss} / \hgam_{M_P} \times \hAp_{P, (0, 1)} \end{tikzcd} \ee whose fibers are homeomorphic to $\hU_P/ \hgam \cap \hU_P.$
	\end{nprop} 
	\begin{proof} For $\gamma\in \hgam_P$, write $\gamma = r_{\gamma} u_{\gamma}$ with $r_{\gamma} \in \hgam_{M_P}$ and $u_{\gamma} \in \hgam_P$ (see Lemma \ref{lem:GammaL}). Then if $z=(m ,a, u)$ is the $P$-horospherical coordinates of an element $z \in \xg(P, e)$ with $m \in M_P, a \in \hAp, u \in U_P$, we have by definition that $z \gamma = (m r_{\gamma} , a, u^{r_{\gamma}} u_{\gamma}).$ The result follows from this. \end{proof}

\subsection{Miscellanea on $\xg^{ss}$} 

\tpoint{The semi-stable locus is closed}  We begin by observing the following result. 

\begin{nprop} \label{prop:xss-closed} We have $\xgss$ is closed in $\xg$. \end{nprop}

\begin{proof} Let $x_n \in \xgss$ be a sequence converging to $x_{\infty} \in \xg$. Recall that the topology on $\xg$ is the product topology coming from the identification $\xg = \hA \times \hU$ (and this is the same as the topology coming from a $P$-horospherical decomposition as well). It suffices to show that $\deg_{\inst}(x_{\infty}) \geq 0.$ Suppose $P \subsetneq G$ is a standard parabolic and $\gamma \in \hG_{\Q}/P_{\Q}$. Then we know that $\la \rho, H_P( x_n \gamma)  \ra \geq 0$ as $x_n \in \xgss.$ On the other hand, by continuity of $H^+_P: \xg \rr \hfh^+_P$ and the linearity (hence continuity) of the pairing with $\rho$, we obtain  \be{} \la \rho, H_P(x_{\infty} \gamma) \ra = \lim_{n \rr \infty} \, \la \rho, H_P(x_n \gamma) \ra \geq 0.  \ee As this holds for each $P$ and $\gamma$, we have seen that $\deg_{\inst}(x_{\infty}) \geq 0.$ \end{proof}

\tpoint{On the closure of $\xg(P, e)$}  Unlike $\xgss$, the pieces $\qxg(P)$ are no longer closed when $P$  is a proper, standard parabolic. Indeed, first recall from \eqref{XGP:MPss} that   $\xg(P, e) = X_{\mp}^{ss} \times \hAp^+_{P, (0, 1)} \times \hU_P.$ Since we have a homeomorphism $\xg = X_{M_P} \times \hAp^+_P \times \hU_P$ where the right hand side is equipped with the product topology, and since, as in the previous Proposition $X_{M_P}^{ss} \subset X_{M_P}$ is closed, the closure of $\xg(P, e)$ in $\xg$ is equal to  \be{} \label{closure:Xpe-G}  \cl_{\xg}( \xg(P, e) )  = X_{\mp}^{ss} \times \hAp^+_{P, (0, 1]} \times \hU_P \ee
where $ \hAp^+_{P, (0, 1]}$ is the closure of $\hA^+_{P, (0, 1)}$ in $\hA^+_P$, \textit{i.e.},
\be{} \hAp^+_{P, (0, 1]} =\{a \in \Ac^+_P \mid \la \alpha, H^+_P(a) \ra \leq 0, \alpha \in \Delta_P \mbox{ and } \la \rho, H^+_P(a) \ra \leq 0\}.\ee

\begin{nlem} \label{lem:closure-P-in-X} 
	For $P \in \stdp'$ and $x_\infty \in \cl_{\xg}(\xg(P,e))$
	\be{} \deg_{\inst}(x_\infty) = \la \rho, H^+_P(x_\infty)  \ra \leq 0  \ee
\end{nlem}
\begin{proof}
	Pick a sequence $(x_n)$ with $x_n \in \xg(P,e)$ such that $x_n  \rr x_{\infty}$. It follows that $H_B(x_n) \rr H_B(x_{\infty}).$ As multiplication by any $\gamma \in \hG_{\Q}$ is continuous, it also follows that $H_B(x_n \gamma) \rr H_B(x_{\infty} \gamma)$, and as the projections $\hfh^+ \rr \hfh^+_Q$ are continuous for any $Q \in \ratp,$ we also have $H_Q(x_n \gamma) \rr H_Q(x_{\infty} \gamma)$.

	Now, since $x_n \in \xg(P,e)$, \textit{i.e.} $\cp(x_n)=(P, e)$ with $P \neq \hG$, we have 
	\be{} \begin{array}{lcr} \deg_{\inst}(x_n) = \la \rho, H_P(x_n) \ra < 0 & \mbox{ and } & \la \rho, H_P(x_n) \ra  \leq \la \rho, H_R(x_n\eta) \ra \end{array} \ee for  $R \in \stdp$ and an $\eta \in \hG_{\Q}$, 
	Now taking the limit $x_n \rr x_\infty$ and using the remarks above, 
	\be{} \begin{array}{lcr}  \la \rho, H_P(x_\infty) \ra = \lim_{n \rr \infty} \, \la \rho, H_P(x_n) \ra \leq 0&  \mbox{ and } &  \la \rho, H_P(x_\infty) \ra  \leq \la \rho, H_R(x_\infty \eta) \ra   \end{array}. \ee Hence $\deg_{\inst}(x_{\infty}) = \la \rho, H_P(x_{\infty}) \ra \leq 0$.
	
\end{proof} 

\tpoint{Construction of elements in $\xgss$} \label{subsub:construction-semistable points} We may use the result from the previous section to construct a family of elements in $\xg^{ss}$. Let $I_o:= \{ 1, \ldots, \ell \}$ and $P_o:= P_{I_o}$. Recall also  \be{} \begin{array}{lcr} \hfh^+_{P_o}= \{ m \cc - r \dd \mid m \in \R, r > 0 \}& \mbox{ and } & \hfh^+_{P_o, (-\infty, 0)} = \{ m \cc - r \dd \mid m < 0, r > 0 \}  \end{array}. \ee  For $n  > 0$ and $r >0$, consider \be{} \begin{array}{lcr} H_n  = \frac{-1}{n}\cc - r\dd \in \hfh^+_{P_o, (- \infty, 0)} & \mbox{ and } & a_n:= \exp(H_n). \end{array}. \ee  We claim that $a_n \in \hA^+_{P_o, (0, 1)}$: indeed, $\la \rho, H_n \ra  = \frac{-\rho(\cc)}{n} <0$ and $\la \alpha_{\ell+1}, H_n \ra = - r <0$. So, for any $m_{\infty} \in X_{M_{P_o}}^{ss}$ and $u_{\infty} \in \hU_{P_o}$, we can now consider the elements \be{} x_n:= (m_{\infty}, a_n, u_{\infty}) \in X(P_o, e) = X^{ss}_{M_{P_o}} \times \hA^+_{P_o, (0, 1)} \times \hU_{P_o}. \ee By what we said above, $x_n \rr x_{\infty}:=(m_{\infty}, \eta(s),  u_{\infty})$ where $s=e^{-r}$. From Lemma \ref{lem:closure-P-in-X}, we conclude that  $x_{\infty} \in \cl_{\xg}(\xg(P,e))$ satisfies
\be{} \deg_{\inst}(x_\infty) = \la \rho,H_P(x_\infty) \ra  = 0\ee
Hence $x_\infty \in \xg^{ss}$, and in this way we have established: 

\begin{nprop} \label{prop:semi-stable-elements} For any $r > 0$, if $s:= e^{-r}$, then we have $X^{ss}_{M_{P_o}} \times \{ \eta(s)) \} \times \hU_{P_o} \subset \xg^{ss}$ \end{nprop}  

\begin{nrem} \label{rem:semi-stable-not-closed-in-bordification} It follows from the above claim that $\xg^{ss}$, which is closed in $\xg$, is \textit{not} closed in the bordification $\bord{\xg}$ that will be introduced in \S \ref{sec:bordification} (we use the notation from that section freely now). Indeed, pick a sequence of numbers $r_n > 0$ such that $r_n \rr \infty$, and write $s_n= e^{-r_n}$. Then for any $m_{\infty} \in X_{M_{P_o}}^{ss}$ and $u_{\infty} \in \hU_{P_o}$, the elements $x_n:= (m_{\infty}, \eta(s_n), u_{\infty}) \in \xg^{ss}$ by the Proposition. However, it is easy to see that it converges to $(m_{\infty}, u_{\infty}) \in \e(P_o)$ since $\eta(s_n)^{\alpha_{\ell+1}} \rr 0$ as $n \rr \infty$. \end{nrem}

\tpoint{Semi-stability in the presence of bounded central directions} We note that if $x$ is semi-stable then `positive' central translations $\zeta x$, i.e. $\zeta \in \hA_{\cc}$  such that $\la \rho, H_{\cc}(\zeta) \ra \geq 0,$ will again be semi-stable. Fix $t, \Omega$ as in Theorem \ref{thm:reduction-theory-1}, and consider for $r, C >0$ \be{} \sie^r_{t, \Omega}[C] := \sie_{t, \Omega} \cap \hg^r \cap \{ x \in \xg \mid  -C \leq \la\lambda_{\ell+1}, H_B(x) \ra \leq C \}. \ee  If $x \in \xg$ is semi-stable, then we must have (see Lemma \ref{lem:rho-vs-rho-classical} ) \be{} \la \rho, H^+_B(x) \ra = \la \rho_o + \dcox \lambda_{\ell+1}, H^+_B(x) \ra \geq 0.\ee Now, as $x \in \sie_{t, \Omega}$, we have $\la a_i, H^+_B(x) \ra < t$ for all $i \in I$.  Since $\rho_o$ is a positive linear combination of the $a_i$ for $i \in I_o$, there exists a constant $M:=M(C, t) > 0$ so that if $x \in \sie_{t, \Omega}^r[C]$ is chosen so that \textit{any} of the $\la \alpha_i, H_B^+(x) \ra < -M$ for $i \in I_o$, then  $\la \rho, H_B(x) \ra < 0.$ On the other hand, from Theorem \ref{thm:reduction-theory-1} (3), given any $M > 0$, there exists $r >0$ so that if $x \in \sie_{t, \Omega}^r$, then $\la \alpha_i, H^+_B(x) \ra < -M$ for some $i \in I$. In conclusion, we thus obtain

\begin{nlem} \label{lem:semi-stable-and-r} Let $t, \Omega$ as in Theorem \ref{thm:reduction-theory-1}(1). Given $C >0 $, there exists $M:=M(C, t) > 0$ so that if $r:=r(C)$ is made sufficiently large, then for any $x \in \sie^{r}_{t, \Omega}[C]$ that is also semi-stable, we must have $\la a_{\ell+1}, H_B^+(x) \ra < -M$. \end{nlem}

	\part{Bordifications and their arithmetic quotients}
	
	For each rational parabolic $P \in \ratp'$, we have constructed in \S \ref{sub:bord-A} the corner $\cor(\hA^+_P)$ of $\hA_P^+$. Building upon this, in section \S \ref{sec:bordification} construct a bordification $\bord{\xg}$ of $\xg$ by adding a `boundary component' $\e(P)$ for each $P \in \ratp'$. We then analyze the topological properties of  $\bord{\xg}$ in the rest of this section, and then take up the study of $\bord{\xg}/\hGam$ in Section \S \ref{sec:bord-mod-Gamma}.

	\newcommand{\Gc}{\widehat{G}}
	\newcommand{\Kc}{\widehat{K}}
	
	\section{The bordifcation of $\bord{\xg}$ } \label{sec:bordification} 
		
		\subsection{On $\bord{\xg}$ and its topology.} 
		
		\tpoint{Boundary components} For $P \in \ratp'$, define a corresponding \textit{boundary component}  \be{} \label{def:e} \e(P) = X_{M_P} \times \hU_P, \ee where $X_{M_P}$ was the symemtric space attached to the finite-dimensional group $M_P$ and $\hU_P$ the pro-unipotent radical attached of $P$. Giving $X_{M_P}$ the usual metric topology and $\hU_P$ the pro-unipotent toplogy, $\e(P)$ acquires the product topology. Sometimes we also shall write $\e(\hG):= \xg.$ Let now $P, Q \in \ratp'$ with   $Q \supset P.$ Write $\sta{P} \subset M_Q$ for the corresponding parabolic and recall the decompositions from  \S \ref{subsub:relative-langlands-dec}, especially \eqref{Q:starP}. Noting that the last two decompositions in the latter are compatible with the product topology, we can easily show the following. 
		
		\begin{nlem} 
			\begin{enumerate}
				\item For $P, Q \in \ratp'$ with $Q \supset P$, then  \be{} \label{e:Q:P} \e(Q) = X_{M_P} \times \hA(Q)_{\sP} \times \widehat{U}_P, \ee  
				\item For $P, Q_1, Q_2 \in \ratp'$ with $P \subset Q_1 \subset Q_2$, we have  \be{}\label{eq2:q1:p}   \e(Q_2) = X_{M_P} \times \hA(Q_2)_{\sQ_1} \times \hA(Q_1)_{\sP} \times \hU_P.  \ee 
			\end{enumerate}	
		\end{nlem}
	\begin{proof} The first result follows from \eqref{split glob P} and the remarks we made preceeding the Lemma. The second follows similarly by also using $\hA(Q_2)_{\sP} \cong \hA(Q_2)_{\sQ_1} \cdot  \hA(Q_1)_{\sP}$, see \eqref{Q:starP}. \end{proof}

\tpoint{$\bord{\xg}$ and its topology} \label{subsub:bord-top} We now introduce the \textit{rational bordification} of $\xg$ \be{} \label{def:bord-X} \bordX:= \xg \cup  \, \sqcup_{P \in \ratp' } \e(P). \ee  We may equip it with a topology extending the existing ones on $\xg$ and $\e(P)$ by specifying an additional class of convergent sequences (again using \S \ref{sub:moore-smith}).
\begin{enumerate}

	\item We say $\{ x_n \}$ with $x_n \in \xg$ converges to $x_{\infty}= (x_P, u_P) \in \e(P)$ for $P \in \ratp'$ where $x_P \in X_{M_P}, u_P \in \hU_P$, if writing  $x_n = ( x_{P, n}, a_{P, n}, u_{P, n})$ with $x_{P, n} \in X_{M_P}, a_{P, n} \in \hA_{P}^+, u_{P, n} \in \hU_P$ (see \eqref{split glob P}), \be{} \begin{array}{lccr} (x_{P, n}, u_{P, n}) \rr (x_P, u_P) & \mbox{ and } &   a_{n, P}^{\alpha} \rr 0, \,  \alpha  \in \Delta_P & \text{ for } n \rr \infty, \end{array} \ee and where the first convergence is with respect to the natural product topology on $X_{M_P} \times \hU_P$.
	
	\item If $P, Q \in \ratp'$ with $\hg \supsetneq Q \supsetneq P$ we say $x_n(Q) \rr x_{\infty}(P)$ where $x_{n}(P) \in \e(Q)$ and $x_{\infty}(P) \in \e(P)$ if writing $x_n(Q) = ( x_{n, P}, a_{n, \sP}, u_{n, P} )$ where $x_{n, P} \in X_{M_P}, a_{n, \sP} \in \hA(Q)_{\sP}, u_{n, P} \in \hU_P,$ (see \eqref{e:Q:P}),   \be{} \begin{array}{lcr} (x_{n, P},  u_{n, P} ) \rr x_{\infty}(P) & a_{n, \sP}^{\alpha} \rr 0 & \mbox{ for } \alpha \in \Delta_P^Q \mbox{ and when } n \rr \infty \end{array}. \ee 
			
\end{enumerate}

Using an argument as in \S \ref{sub:bord-h} and \S \ref{sub:bord-A}.

\begin{nprop} \label{prop:bord-top} The above specifies a well-defined convergence class on $\bordX$ and hence defined a topology on this set. \end{nprop} 

\begin{nrem}
    \begin{enumerate}
    	\item Each element in $\e(P)$ can be obtained as the limit of a convergent sequence from $\xg$ and, if $Q \supset P$, each element of $\e(P)$ can be obtained from a convergent sequence in $\e(Q).$

\item  \label{subsub:convergence-along-boundary} It will be useful to make more explicit convergence among the boundary components. Let $P, Q \in \ratp$ with $P \subset Q \subsetneq \hG$ and pick $y_j \in \e(Q)$ which we write as  $y_j = (z_j, n_j) \in \e(Q)$ with $z_j \in X_{M_Q}, n_j \in \hU_Q$. Assume $n_j \rr n_{\infty}$ for $n_{\infty} \in \hU_Q$ with respect to the standard topology on $\hU_Q$.  Suppose that when we decompose $z_j = (z_j^*, a_j^*, n_j^*)$ in terms of the $\sP$-horospherical coordinates of $X_{M_Q}$, i.e. with respect to  $X_{M_Q}  = X_{\sP} \times \hA(Q)_{\sP} \times U_{\sP}$ of \eqref{split glob P}, we have \begin{itemize} \item $z^*_j \rr z^*_{\infty}$ and $n_j^* \rr n_{\infty}^*$ with $z^*_\infty \in X_{M_P}$ and $n^*_{\infty} \in \hU_{\sta{P}}$ in the topology of $X_{\sP}=X_{M_P}$ and $\hU_{\sP}$ respectively; \item  $(a_j^*)^{\alpha} \rr 0$ for $\alpha \in \Delta_P^Q$. \end{itemize} Then $y_j \rr y_{\infty} \in \e(P)$ with  $y_{\infty} = ( z_{\infty}^*, n_{\infty}^* n_{\infty} ).$   

\end{enumerate}
\end{nrem}

	\tpoint{Embedding extension result}   Recall that $\cor(\hA_P^+) = \hA_P^+ \cup \, \bigsqcup_{Q\supset P} \hA(Q)_{\sP},$ so that \be{} X_{M_P} \times \cor(\hA^+_P) \times \widehat{U}_P  = \bigsqcup_{Q\supset P} X_{M_P} \times \hAp(Q)_{\sP} \times \hU_P. \ee We can define a map $\iota:=\iota_P$ from this set to $\bord{\xg}$  by sending the point $(x_P, h_{\sP}, u_P)$ with $x_P \in X_{M_P}, h_{\sP}\in \hA(Q)_{\sP}, u_P \in \hU_P$ to the corresponding point of $\e(Q)$ specified by  \eqref{e:Q:P}. 
	
	\begin{nlem} \label{lem:embedding-extension-lemma} The map $\iota=\iota_P$ is a topological embedding \be{} \iota_P: X_{M_P} \times \cor(\hA^+_P) \times \widehat{U}_P \hookrightarrow \bord{\xg} \ee where the set on the left is equipped with the product topology. \end{nlem} 

\begin{nrem}  \label{remark:corner} The image of the map $\iota_P$ described above will be denoted as $\cor(P);$ it is an open set in $\bord{\xg}$ and is called the \textit{corner attached to $P$}. It is easy to see that  \be{} \label{corner:dec} \cor(P) = \xg \cup \, \bigsqcup_{Q\supset P} \e(Q) \ee where the (finite) union is over $Q \in \ratp'$.   So for example if $P=\hB$, then $\cor(\hB)$ is the union of $\xg$ and $\e(P)$ for all (proper) standard parabolics.
\end{nrem} 

 \begin{proof}  In light of the topology we have introduced in \S \ref{sub:bord-A} on $\cor(\hA_P^+)$, to verify that $\iota$ is continuous, it suffices to check the following two facts:  \begin{enumerate}
\item If $y_n \in X_{P} \times \hA^+_P \times \widehat{U}_P \subset X_{P} \times \cor(\hA_P^+) \times \widehat{U}_P$ is a convergent sequence in the product topology with limit $y_{\infty} \in X_{P} \times \hA(Q)_{\sP} \times \widehat{U}_P$  for some $Q \supset P$, then the image sequence $\iota(y_n) \in \xg$ converges in the  topology of $\bord{\xg}$ to  $\iota(y_{\infty}) \in \e(Q)$. 
\item If $y_n \in X_{P} \times \hAp(Q)_{\sP} \times \widehat{U}_P $ is a sequence, convergent in the product topology to  $y_{\infty} \in X_{P} \times \hA(Q')_{\sP} \times \widehat{U}_P$  for some $Q'\in \ratp(Q), Q' \supset P$, then the image sequence $\iota(y_n) \in \e(Q)$ converges in the  topology on $\bord{\xg}$ to $\iota(y_{\infty}) \in \e(Q')$. 
\end{enumerate}

\noindent As the proof of the second assertion is similar to the first, we just concentrate on that one. Suppose we are given a sequence of points  $y_n \in \xg$, say $y_n:= (x_{P, n},  h_{P, n},  u_{P_n})$ with $x_{P, n} \in X_{M_P}, h_{P, n} \in \hA_P^+, u_{P, n}\in \hU_P$ which converges in the product topology of $X_{M_P} \times \cor(\hA_P^+) \times \widehat{U}_P$ to a point \be{} y_{\infty}:=(x_P, h_{\sP},u_{P, n}) \in X_{M_P} \times \hA(Q)_{\sP} \times \hU_P \ee  for some $Q \supset P$. 
 In other words, we are assuming that  \be{}\begin{array}{lcr} h_{P, n }  \rr h_{\sP} , &   x_{P, n} \rr x_P&  u_{P, n} \rr u_{P} \end{array} \ee for $h_{\sP} \in \hA(Q)_{\sP}$, $x_P \in X_{M_P}$,  $u_P \in \hU_P$ and with respect to the topologies of  $\cor(\hA_P^+),$ $X_{M_P}$ and $\hU_P$ respectively.
We need to verify that $\iota(y_n) \rr (x_P, h_{\sP}, u_P) \in \e(Q)$ in the topology we have introduced on $\bord{\xg}.$ Recalling from (\ref{split glob P}) that for $Q \supset P$ \be{} \label{recall-PQ}  \begin{array}{lcr}  X_{M_Q} = X_{M_P} \times \hA(Q)_{\sP}\times U_{\sP}, & 
		\hA_{P}^+ = \hA(Q)_{\sP}\times \hA_{Q}^+, & \hU_{P} = \hU_{\sP} \hU_Q \end{array}. \ee  Writing $h _{P, n} = (h_{P, n})_{\sP} \times (h_{P, n})_{Q} $  and $ u_{P_n} = (u_{P_n})_{\sP} (u_{P_n})_{Q} $, a $Q$-horospherical decomposition  of $y_n$ is  
		 \be{} y_n = \underbrace{( x_{P, n} , (h_{P, n})_{\sQ(P)}, (u_{P_n})_{\sP} )}_{X_{M_Q}} \times \underbrace{(h_{P, n})_{Q}}_{\hA^+_Q} \times \underbrace{(u_{P_n})_{Q},}_{\hU_Q}. \ee  
		\noindent By definition of convergence in $\cor(\hA_P^+)$ if $h_{P, n} \rr h_{\sP}$ we must have $(h_{P, n})_{\sP}^{a} \rr 0$ for $a \in \Delta_P^Q$. As the topology on $\hU_P$ is compatible with the product decomposition in \eqref{recall-PQ}, the result  follows from the construction of the topology on $\bord{\xg}$ specified in \S \ref{subsub:bord-top}  
		
		We leave the remainder of the proof to the reader. 
		
		\end{proof}

\tpoint{A neighborhood basis in $\bord{\xg}$} If $P = \hG$, one can take a neighborhood basis of $\xg:= \e(\hg)$ coming from a basis of the topologies of $\hU$ and of $\hA^+$. As for the boundary, we have

\newcommand{\Nb}{\mathcal{N}}

\begin{nlem} \label{sp nbhd} Let $P \in \ratp'$ and  $(x,u) \in X_{P} \times \hU_P =\e(P).$ A neighbourhood basis of $(x,u)$ in $\bordX$ can be taken to be of the form 
	\be{} \label{nhbd:W-sigma-V} \Nb_P(W, \sigma, V) := W \times \cor(\Ac^+_{P,\sigma}) \times V, \ee  where $W \subset X_{M_P}$ and $V \subset \hU_P$  vary over a basis of neighbourhoods for $x \in X_{M_P}$ and $u \in \widehat{U}_P$ respectively, and $\sigma$ varies over $\R_{>0}$ (recall that  $\cor(\hA^+_{P, \sigma})$ was defined in \eqref{bord:A:sigma}).
	
	 \end{nlem}

\begin{proof} For a given $x_{\infty}:=(x, u) \in \e(P)$, we will argue that if $x_n$ is a convergent sequence with limit $x_{\infty}$, then there exists a basis element $ \Nb_P(W, \sigma, V)$ as above which contains infinitely many $x_n$. As usual, we have two cases. First, we assume that each $x_n \in \xg$. Writing a $P$-horospherical decomposition $x_n = (x(P)_n , a_{P, n}, u_{P, n}) \in X_{P} \times \Ac_{P}^+ \times \widehat{U}_P$, this means we have  \be{} \begin{array}{lcr} u_{P, n} \rr u, & x(P)_n \rr x, & a_{P,n}^{\alpha} \rr 0  \hspace{0.1cm}\mbox{ for all }  \alpha \in \Delta_P \end{array}.\ee
	So, in particular we can choose an $N_0 \in \mathbb{N}$, and open sets $W $ and $U$ in $X_{M_P}$ and $\hU_P$ containing $x$ and $u$  respectively, along with a $\sigma >0$ small enough such that for $n \geq N_0$ the elements $x_n$ lie in $W \times \Ac^+_{P,\sigma} \times U \subset  \Nb_P(W, \sigma, V)$. For these choices infinitely many $x_{n}$ will lie in $ \Nb_P(W, \sigma, V)$.
	
	The other case we need to consider is that $x_n \rr x_{\infty}$ along the boundary. Without loss of generality, we may assume there exists $Q \in \ratp$ with $P \subset Q$ such that $x_n \in \e(Q).$  If $x_n \rr x_{\infty}$, we follow the procedure in the second remark after Proposition \ref{prop:bord-top}: first decompose $x_n \in e(Q)$ using  $e(Q) = X_{P} \times \hAp(Q)_{^*P} \times \widehat{U}_P$ by writing $x_n  = (x_{P,n}, a_{^*P,n},  u_{P, n})$, then we must have have
	\be{} \begin{array}{lcr} x_{P,n} \rr x &  u_{P,n} \rr u & \text{ and } a_{^*P,n}^{\alpha} \rr 0 \hspace{0.1cm} \text{ for all }  \alpha \in \Delta_{P}^Q. \end{array} \ee Fix $N_o \in \mathbb{N}$ and choose  $\sigma >0$ small enough such that for $n \geq N_0$ the elements $a_{\sP, n} \in \hAp(Q)_{\sP, \sigma}.$ Then if $W$ and $U$ are as in the previous paragraph, since $W \times \hAp(Q)_{\sP,\sigma} \times V \subset \Nb_P(W, \sigma, V)$ we see that the latter contains infinitely many $x_n$. 
\end{proof}

\subsection{The right action of $\hG_{\Q}$ on $\bord{\xg}$ }

Right multiplication by $\hG_{\Q}$ on $\hG^+$ induces  \be{} \begin{array}{lcr} \label{def:R-act} R_g: \xg \longrightarrow \xg, &  z \mapsto z . g & \mbox{ for } z \in \xg, g \in \hG_{\Q}. \end{array} \ee This aim of this subsection is to first show that $R_g$ is continuous with respect to the product topology on $\xg$ obtained by regarding it as $\xg = \hA^+ \times \hU$ and then extend this result to $\bord{\xg}$. 

\tpoint{Continuity of $R_g$}

\begin{nprop}

    The right $\hG_{\Q}$- action $R_g$ on $\xg$  is continuous.
\end{nprop}

\begin{proof} In light of the homeomorphism $\xg = \hA^+ \times \hU$ and the  Bruhat decomposition $\hG_{\Q} = \displaystyle \bigcup_{w \in \widehat{W}} \hB_{\Q}w\hU_{\Q}$ (see \S \ref{subsub:BN-pair}), to verify the continuity of $R_g$ it is enough to check verify the following. 
    \begin{enumerate}
    \item The multiplication map $\hU\times \hU\rr \hU$ is continuous. This was already explained in Proposition \ref{prop:continuity-U-mult} and we adopt the same notation as in that proof. In particular, we will just explain the key idea for the group $\widehat{\SL_2},$ as the generalization to higher rank is straightforward.

    \item If $g  = t_0 \in \hT_{\Q}$ and we pick a point $(a,u) \in \xg$ with $a \in \hA^+$ and $u \in \hU$ and  note  $R_g((a,u))  = (at_0, t_0^{-1}ut_0).$  So, the continuity of the map $R_g$ is reduced to showing the continuity of conjugation of $\hT_{\Q}$ on $\hU$ as well as that of right multiplication of $\hT_{\Q}$ on $\Ac^+.$ The latter is clearly continuous. We will next argue that for any $t \in \hT$, it's conjugation on $\hU$ is continuous. Consider $u$ as in \eqref{sl2:unip} and let $t \in \hT$ be written as $t = \eta(s) h_{\cc}(t) h_{\alpha}(u)$ with $s, t, u \in \R^*$ and $h_{\cc}(t)$ is as in \S \ref{def:hc}.  As $h_{\cc}(t)$ is central, we treat the case of conjugation by $\eta(s)$ and $h_{\alpha}(u)$ separately. From direct computation (in the notation of \S \ref{subsub:U-hat-top}) using \cite[Lemma 11.2]{gar:ihes} shows \be{} \begin{array}{lcr} h_{\alpha}(u) \, \chi_{\alpha}(\bbsigma) h_{\alpha}(u)^{-1} =  \chi_{\alpha}(s^{2} \bbsigma) & \mbox{ and } & h_{\alpha}(u) \, \chi_{-\alpha}(\bbsigma) h_{\alpha}(u)^{-1} =  \chi_{-\alpha}(s^{-2} \tau(t)) \end{array}     \ee Further $h_{\alpha}(u)$ commutes with $h_{\alpha}(\bbmu)$ using the computation of symbols in \cite[\S12]{gar:ihes}. Hence conjugation by $h_{\alpha}(u)$ on $\hU$ is continuous.  As for $\eta(s)$, we first note that \be{} \eta(s) \chi_{\alpha} ( \bbsigma(t)) \eta(s)^{-1} = \chi_{\alpha}(\bbsigma(s t)).\ee On the other hand, the conjugation of $\eta(s)$ on $h_{\alpha}(\bbmu)$ can be described using \cite[Lemma 3.3.3]{pat:birkhoff}, which again shows its polynomiality and hence continuity.

    \item Finally, we argue that the right action by $g = \dw_i : = w_{a_i}(1)$ is continuous. Using the description of the topology given in \S \ref{subsub:U-hat-top}, one show that $\pi^i$ and $\pi_i$, projections stemming from  \be{} 
    \begin{array}{lcr} \hU= \hU^i \ltimes  U_i & \mbox{ where }  U_i = \{ \chi_{a_i}(u) \mid u \in \R \} & \mbox{ and }  \hU^i = \hU\cap \dw_i ^{-1} \hU \dw_i \end{array},\ee are both continuous.  To compute the action of $R_{\dw_i}(a, u)$, we note that writing $u=u^-u_i$ with respect to the decomposition above, we have \be{} a u \dw_i = a u_i u^i \dw_i = a^{w_i} u_{-i} (u^i)^{\dw_i}  \ee where  $u_{-i}:= \dw_i \chi_{a_i}(u) \dw_i = \chi_{-a_i}(u)$. As the map $a \mapsto a^{\dw_i}$ is clearly continuous, it remains to verify that the maps: (i) $(h, 1) \mapsto (h, 1) R_{u_{-i}}$ for $h \in \hA^+$; and (ii) the conjugation action of $\dw_i$ on $\hU^i$ are continuous. Note that  (i) can be reduced to a computation in $\SL_2(\R)$. As for (ii), a formula for this conjugation can be obtained by working in the group $\SL_2(\R((t)))$ (so modulo the center) where it is cleary seen to be polynomial.

\end{enumerate}
\end{proof}

\tpoint{An identification} \label{subsub:eP-action} Let $P \in \ratp$ and let $g \in \hG_{\Q}.$ Decompose $g$ as  
\be{}  \label{g:iwa} g=   p \, k  \ \mbox{ where } k \in \hK, p \in P. \ee Although $k$ is not well-defined (since $\hK \cap P$ is non-empty), we do have that $P^k  = P^g$, i.e.  $P^k$ does not depend on the precise choice of  $k$ in the factorization \eqref{g:iwa}.Next, we note that \be{} \begin{array}{lcr} M_{P^k} = (M_P)^k & \mbox{ and } & U_{P^k} = (U_P)^k \end{array}. \ee Hence we obtain a canonical identification between $\e(P)$ and $\e(P^k)=\e(P^g)$ as follows  \be{} \label{ePk:P} \begin{array}{lcr} \e(P) . k = \e(P^k) & \mbox{ sending} & (x_P, n_P).k \mapsto (x_P^k,  n_P^k) \end{array}. \ee By the way in which the topologies on $X_{M_{P^k}}$ and $\hU_{P^k}$ are constructed, this map is a homeomorphism.

\tpoint{Extension of $R_g$ to $\bord{\xg}$} \label{subsub:def-Rg-bord} To motivate the right action of $R_g$ on $\e(P)$ let us express $R_g: \xg \rr \xg$ in $P$-horospherical coordinates. Suppose $x = (m_x, a_x, n_x)$ with $m_x \in M_P, a_x \in \hA^+_P, n_x \in \hU_P$. Writing $g$ as in \eqref{g:iwa} and decomposing $p=n \,m \, a$ for $m \in M_P, a \in \hAp_P, n \in \hU_P$,   \be{} \label{Rg:act-interior}  x g = (m_x m)^k \,  (a_x a)^k \, \left( (n_x n)^{m a} \right)^k.  \ee  Note that the product $ma$ is  well-defined (i.e., depends just on $g$), even though $m$ and $a$ may not be ($k$ and $n$ are also well-defined).  So for $P \in \ratp$ and $g \in \hG_{\Q}$, let us now define $R_g: \e(P) \rr \e(\leftidx^k P)$ by writing $z \in \e(P)$ as $z = (  m_z, n_z) $ with $m_z \in X_{M_P}, n_z \in \hU_P$ and setting  \be{} \begin{array}{lcr} \label{Rg:act-boundary} R_g (z)= (m_z m  , (n_z n)^{ma} ).k& \mbox{ for } & g = p k = ma n k \end{array} , \ee using the identification of  \eqref{ePk:P}. It is easy to check that this process defines a right action of $g$ on $\bord{\xg}$, that we henceforth also denote as \be{} \label{def:Rg} R_g: \bord{\xg} \rr \bord{\xg}. \ee

\begin{nprop} \label{prop:cont-Rg}  For each $g \in \hG_{\Q}$, the map $R_g: \bordX \rr \bordX$ is continuous. 
\end{nprop}

The argument is almost identical to \cite[Prop. III.9.15]{borel-ji}, but we reproduce it here for completeness.  There are two cases to consider, depending on whether the convergence is entirely along the boundary of $\bordX$ or not. 

\tpoint{Proof of Proposition \ref{prop:cont-Rg}, Step 1: convergence into the boundary } \label{proof-continuinity-into-boundary} Suppose a sequence $\{y_j \}$ of points in  $\xg$ converges to a point $y_{\infty} \in \e(P)$ for $P \in \ratp$. Write $y_{\infty}= ( m_{\infty}, n_{\infty})$ with $m_{\infty} \in X_{M_P}$ and $n_{\infty} \in \hU_{P}$. The convergence $y_j \rr y_{\infty}$ means that if we write in $P$-horospherical coordinates $y_j = (m_j, h_j, n_j)$ with $m_j \in X_{M_P}, h_j \in \hA_P^+, n_j \in \hU_P$, we have $m_j \rr m_{\infty}$, $n_j \rr n_{\infty}$ in the topology of $X_{M_P}$ and $\hU_P$ respectively, and also that $h_j^{\alpha} \rr 0$ for $\alpha \in \Delta_P$. 

We need to check that $R_g(y_j) \rr R_g(y_{\infty}) \in \e(P^g)$. From \eqref{Rg:act-boundary}, we compute \be{} \begin{array}{lcr} R_{g}(y_{\infty}) = z_{\infty}.k & \mbox{ where } & z_{\infty}:= (m_{\infty} m, (n_{\infty}n)^{ma} ) \end{array} \ee and where $g = p k = man k$ is as above. On the other hand, using \eqref{Rg:act-interior}, we have that \be{} R_g(y_j)=  m_j m \,  h_j a \, (n_j n)^{ma} k = z_j k,  \ee where $z_j \in \xg$ has $P$-horospherical coordinates $( m_j m, h_j a, (n_j n)^{ma}).$ By the previous paragraph, we have $z_j \rr z_{\infty}$, so we can finish the argument using the next result, which follows from the observation made earlier that right multiplication by $k$ is a homeomorphism from $\e(P)$ to $\e(P^k)$ 

\begin{nclaim}Suppose $z_j$ is a sequence of elements in $\xg$ which converge to $z_{\infty} \in \e(P)$ for some $P \in \ratp$. If $g \in \hG_{\Q}$ with $g=kp$ as in \eqref{g:iwa}, then $z_j k \rr z_{\infty}. k$ as elements in $\e(P^k)$ \end{nclaim}

\tpoint{Proof of Proposition \ref{prop:cont-Rg}, Step 2: convergence along the boundary.}  Without loss of generality we may assume that we are in the setting of \S \ref{subsub:convergence-along-boundary} and we adopt all of the notation as in \emph{loc. cit.}, \textit{i.e.} we have $y_j=(z_j, n_j) \in \e(Q)$ converging to $y_{\infty}=(z_{\infty}, n_{\infty}) \in \e(P)$, etc. We are aiming to show \be{} R_g(y_j) \rr R_g(y_{\infty}) \ee and we do this by explicitly computing both quantities above. 

\begin{enumerate}
\item To compute $R_g(y_j)$ we first factor $g \in \hG_{\Q}$  as in \eqref{Rg:act-boundary}, say $g = m a n k$ with $k \in \hK, m \in M_Q, a \in \hA_Q, n \in \hU_Q$. Then from \eqref{Rg:act-boundary} , we find \be{} R_g(y_j) = \left(z_j m, (n_j n)^{a m} \right). k, \ee where we note that $q_j:= (z_j m, (n_j n)^{am}) \in \e(Q)$. Let us first compute the limit of $q_j$, for which we follow the receipe of \S \ref{subsub:convergence-along-boundary}  and first factor $z_j m$ using the $\sP$-horospherical decomposition of $X_{M_Q}$. We leave it as an exercise to verify that if $z_j=(z_j^*, a_j^*, n_j^*)$ and $m=(m^*, a^*, n^*)$ are the $\sP$ factorizations of $z_j$ and $m$ respectively, then \be{} z_j m = ( z_j^* m^*, a_j^* a^*, (n_j^*)^{m^* a^*} n^*) \ee is the corresponding factorization for $z_j m$. Applying Lemma \ref{lem:boundary-convergence}  we find that \be{} q_j \rr q_{\infty}:= (z_{\infty}^*m^*, (n^*_{\infty})^{m^*a^*}\, n^* n_{\infty}), \ee and hence \be{} \label{limit-Rgy:j} R_g(y_j) \rr (z_{\infty}^*m^*, (n^*_{\infty})^{m^*a^*}\, n^* n_{\infty}).k \in \e(P^k). \ee

\item We next need to compare \eqref{limit-Rgy:j} with the expression $R_g(y_{\infty})$. For this we need to write a $P$-horospherical decomposition for $g$ in terms of its $Q$-horospherical decomposition $g=m a n k$ given above. Writing $m= (m^*, a^*, n^*)$ as in the previous step, one may check that the $P$-decomposition of $g$ becomes \be{} g = m^* (a^*a) (n^*)^a n) k. \ee 
Then we have using \eqref{Rg:act-boundary} that \be{} R_g(y_{\infty}) =\left( z_{\infty} m^*,  \left(  n^*_{\infty} n_{\infty} \, ((n^*)^a n) \right)^{m^* a^* a} \right).k \ee

Comparing with \eqref{limit-Rgy:j}, we are done.

\end{enumerate}

\tpoint{Hausdorff property}  Recall that the topologies on $\xg$ and $\e(P), \, P \in \ratp$ are Hausdorff. We have also seen that $\cor(\hA_P^+)$ has a Hausdorff topology (see Prop. \ref{basis bordh}) and so by Lemma \ref{lem:embedding-extension-lemma}, the corner $\cor(P)$ has one as well. Extending this further, we now prove

\begin{nprop} \label{prop:Hausdorff} The topology on $\bordX$ is Hausdorff. \end{nprop}

\begin{proof}	Since the topology on $\bordX$ is defined through a convergence class of sequences, using Proposition \ref{lem:top-from-seq} (1), it is enough to show that a convergent sequence in  $\bordX$ has a unique limit. As usual, there are two cases to consider depending on whether the convergence is from the interior to the boundary or whether it is entirely along the boundary. The argument for both cases is similar, and proceeds by showing that the convergence may be taken inside a space already known to be Hausdorff. 
	
	We begin with the case of convergence along the boundary. Let $Q \in \ratp$ be proper and suppose we have sequence $x_n \in \e(Q)$ converging to points $x_{\infty} \in \e(P_1)$ and also to $y_{\infty} \in \e(P_2)$, with $P_i \subset Q$ for $i=1, 2.$ 
	
	Note that $M_Q$, a finite-dimensional group, with attached symmetric space $X_{M_Q}$ has a bordification as in \cite[\S 7]{borel-serre}, which we denote as $\bord{X_{M_Q}}$. It is constructed in the same fashion as  $\bord{\xg}$ (see \cite[\S III.9 ]{borel-ji}), and with boundary components by $\e_{M_Q}(\sP), \, P \subset Q$. We have \be{} \label{bord:MQ} \bord{X_{\mq}} = X_{\mq} \cup \, \sqcup_{P \in \ratp(Q)} \e_{\mq}(\sP) \ee 
	
	Denote by $\pi_Q$ the natural projection by $\pi_{Q}: \e(Q) \longrightarrow X_{M_Q}$ and consider $\e_{\mq}(\sP_i) \subset \bord{X_{\mq}}$ for $i=1, 2$ the boundary components corresponding to $P_i.$ By the description of convergence along the boundary given in Remark (2) after Proposition \ref{prop:bord-top}, it follows that there exists $x^*_{\infty} \in \e_{\mq}(\sP_1)$ and $y^*_{\infty} \in \e_{\mq}(\sP_2)$ such that $\pi_Q(x_n) \rr x^*_{\infty}$ and $\pi_Q(x_n) \rr y^*_{\infty}$, where the convergence is in the topology of $\bord{X_{M_Q}}$. Since $\bord{X_{\mq}}$ is Hausdorff (\cite[Theorem 7.8]{borel-serre}) we have $x^*_{\infty} = y^*_{\infty}$ and consequently $\sP_1 =  \sP_2$, which in turn implies $P_1 =P_2.$ By picking a minimal parabolic $S \subset P_1$, the convergence takes place entirely within the corner $\cor(S)$, which we already commented was Hausdorff, and so $x_{\infty} = y_{\infty}$.

	Now we consider the case of of a sequence $\{ x_n \}$ of points from $\xg$ which converges to two points in the boundary, say $x_{\infty} \in \e(P)$ and to $y_{\infty} \in \e(Q)$ where $P,Q \in \ratp$. Without loss of generality, we can assume $P, Q^{\gamma} \in \stdp$  for some $\gamma \in \hGam$. From Lemma \ref{sp nbhd}, we find that if $R \in \stdp$ and $x_n \in \xg$ is a sequence such that $x_n \rr x_{\infty} \in \e(R)$, there exists $t >0$ such that for $n \gg 0$, we must have $x_n \in \hAp^+_t \times \hU.$

	Applying this observation to $R=P=P_J$ and $x_n \rr x_{\infty} \in \e(P)$, there exists $t > 0$ such that almost all $x_n$ lie in $ \hAp^+_t \times \hU$, and for $i \notin J$, we have $x_n^{a_i} \rr 0$. On the other hand, the convergence of $x_n \rr y_{\infty}$ implies that $R_{\gamma}(x_n) \rr R_{\gamma}(y_{\infty})$, i.e. we have $x_n \gamma \rr R_{\gamma}(y_{\infty}) \in \e(Q^{\gamma}).$ So, again we find (increasing the previous $t$ if necessary) that for almost all $n$, we also have $x_n \gamma \in \hAp^+_t \times \hU.$  In sum we have \be{} \begin{array}{lcr} x_n, x_n\gamma \in \hAp^+_t \times \hU  & \text{ and } & x_n^{a_i} \rr 0 \text{ for } i \in I \setminus J \end{array}.\ee
	
	 Applying the remarks after Proposition \ref{thm:reduction-theory-1}, we conclude that $\gamma \in \hGam \cap P_J$, which gives us $P_J \cap Q \supset \hB^{\gamma}$. Hence both $\e(P_J)$ and $\e(Q)$ lie in the corner $\cor(\hB^{\gamma})$, which is Hausdorff, so  $x_{\infty} = y_{\infty}$.

\end{proof}

\subsection{On the boundary of $\bordX$ and the affine Rational Tits Building}

\tpoint{Closures and intersections } If $Q \in \ratp$, write $\ratp(Q)$ for the (infinite) set of parabolics contained in $Q$.

\begin{nprop} \label{prop:closure-e} Let $P, Q \in \ratp'$. 
   \begin{enumerate}
   	\item We have a disjoint union, 
    \be{} \label{e(Q) clos decom} \overline{\e(Q)} = \bigsqcup_{P \in \ratp(Q)} \e(P) \ee
	\item  The intersection $\overline{\e(P)} \cap \overline{\e(Q)}$  is equal to $\overline{\e(P \cap Q)}$ if $P \cap Q \in \ratp$ and is empty
	otherwise. In particular, $\overline{\e(P)} = \overline{\e(Q)}$ if and only if $P =Q$.
	\item  \label{prop:corners-open-closure-e}  We have the equivalences: $\e(P) \cap \overline{\e(Q)} \neq \emptyset \iff \e(P) \subset \overline{\e(Q)}  \iff P \subset Q$.
	\item Let $Q \in \ratp$, then an open neighbourhood of $\overline{\e(Q)}$ in $\bordX$ is given by 
	\be{} Y(Q) = \bigcup_{P \in \ratp(Q)} \cor(P), \ee where $\cor(P)$ was the corner introduced in \eqref{corner:dec}.
	
	\end{enumerate}

\end{nprop}
\begin{proof}
Part (1)  follows from the choice of the convergence class used to construct $\bordX$ and the remark after \eqref{def:bord-X}.
As for (2), for an $R \in \ratp$, using \eqref{e(Q) clos decom} we see that $\e(R) \subset \overline{\e(P)} \cap \overline{e(Q)}$ if and only if $ R \subset P \cap Q.$ However,  this is only possible if $ P \cap Q \in \ratp$, since any subgroup containing a parabolic is also a parabolic (this follows fom the axioms of a BN pair). Hence 
   \be{}  \overline{\e(P)} \cap \overline{\e(Q)} = \bigsqcup_{R \in \ratp(P \cap Q)}\e(R) = \overline{\e(P \cap Q)}\ee
Part (3) follows from \eqref{e(Q) clos decom} As $\cor(P) \subset \bord{\xg}$ is open, (4) the result follows from (1).

\end{proof}

\spoint Recall the bordification $\bord{X_{M_Q}}$ for $Q \in \ratp'$ from \eqref{bord:MQ}.

\begin{nprop} \label{prop:eQ-closure-fin-dim}
	For $Q \in \ratp'$ we have a homeomorphism 
	\be{} \label{e(Q) bar des}\overline{\e(Q)} \simeq \bord{X_{M_Q}} \times \widehat{U}_Q \ee
	where the set on the right is equipped with product topology. 
\end{nprop}

\begin{proof}
	
	For $P \in \ratp(Q)$, as $\e_{M_Q}(\sP) = X_{M_{\sP}}\times U_{\sP} = X_{M_P} \times U_{\sP}$ and  $\hU_P = U_{\sP} \times \hU_Q$ (the last is topological, not group theoretic, direct product; see \eqref{Q:starP}), we have 
	\be{} \e(P) = X_{M_P} \times \hU_P = X_{M_P} \times U_{\sP} \times \hU_Q = \e_{\mq}(\sP) \times  \widehat{U}_Q.  \ee We can use this to define projections
	\be{}  \begin{array}{lcr} \pr^P_{s,Q}: \e(P) \longrightarrow \e_{M_Q}(^*P) & \mbox{ and } &  \pr^P_{u,Q}:\e(P) \longrightarrow \widehat{U}_Q  \end{array}. \ee Using the description \eqref{e(Q) clos decom}, we can define maps \be{} \begin{array}{lcr} \Psi_s: \overline{\e(Q)}  \longrightarrow \bord{X_{\mq}} & \mbox{ and } & \Psi_u: \overline{\e(Q)}  \longrightarrow \widehat{U}_Q  \end{array} \ee by requiring \be{} \begin{array}{lcr} \Psi_s|_{\e(P)} = \pr^P_{s,Q} & \mbox{ and }   \Psi_u|_{\e(P)} = \pr^P_{u,Q} \end{array}. \ee 
	
These in turn give a map $\Psi:=\Psi_s \times \Psi_u: \ov{\e(Q)} \rr \bord{X_{M_Q}} \times \widehat{U}_Q $
	fitting into the diagram \be{} 
	\begin{tikzcd}[row sep=huge]
		& \overline{\e(Q)}\ar[dl,"\Psi_s",sloped] \ar[dr,"\Psi_u",sloped] \ar[d,dashed,"{ \Psi= \Psi_s \times \Psi_u }" description] & \\
		\bord{X_{\mq}} & \bord{X_{\mq}} \times \widehat{U}_Q\ar[l,"\pr_{M_Q}"] \ar[r,"\pr_{U_Q}",swap] & \widehat{U}_Q
	\end{tikzcd} \ee The result will follow if we can show
	\begin{enumerate}
		\item The maps $\Psi_s$ , $\Psi_u$, and hence $\Psi$ are continuous.
		\item The map $\Psi$ has a continuous inverse.
	\end{enumerate}
	
	To verify the continuity of $\Psi_s$ and $\Psi_u$, it is enough to check that if $y_j \rr y_{\infty}$ in the topology of $\overline{\e(Q)}$ then $\Psi_s(y_j) \rr \Psi_s(y_{\infty})$ and similarly for $\Psi_u$. By the definition of the convergence class, we may further assume there exist $P_1, P_2 \in \ratp(Q)$ such that  $y_j \subset \e(P_2) $ and $y_{\infty} \in \e(P_1)$. Using the description in \S \ref{subsub:convergence-along-boundary} and considering the definition we just gave for $\Psi_s, \Psi_u$,  the required claims follows.

	Next we construct a map $\Phi: \bord{X_{\mq}} \times \widehat{U}_Q \longrightarrow \overline{\e(Q)}:$  for $P \in \ratp(Q),$ the restriction of $\Phi$ to $\e_{M_Q}(^*P) \times \widehat{U}_Q$ is given by 
	\be{} \begin{array}{lcr} \e_{M_Q}(^*P) \times \widehat{U}_Q \rr \e(P) & \mbox{  } &  ((z_{^*P},u_{^*P}),u_Q) \mapsto (z_{^*P},u_{^*P}u_Q)  \end{array}, \ee where  $z_{\sP} \in X_{M_{\sP}}=X_{M_P}$, $u_{\sP} \in U_{\sP}$, and $u_Q \in \hU_Q$. It is immediate that $\Phi$ is  inverse to $\Psi$. 
	
	To check the continuity of $\Phi$ we will show that $\Phi^{-1}(B)$ is open for every element in a basis for the topology of $\overline{\e(Q)}.$ Using Lemma \ref{sp nbhd}, we may take $B$ of the form 
	\be{} B_P(W,\sigma, V) = \Nb_P(W, \sigma, V) \cap \overline{\e(Q)} = W \times \bigsqcup_{Q \supset R \supset P}\, \hAp(R)_{^*P,\sigma} \times V \ee
	where $P \in \ratp(Q)$, $\sigma > 0$, and $W \subset X_{M_P}$ and $V \subset \widehat{U}_P$ open subsets. Since we have a topolgical product $\widehat{U}_P = U_{^*P} \times \widehat{U}_Q$ we may pick $V$ of the form $V  = V_{^*P} \times V_{Q}$ with $V_{^*P}$ open in $U_{^*P}$ and $V_{Q}$ open in $\widehat{U}_Q$. Then we have
	\be{} \Phi^{-1}(B_P(W,V,\sigma))  = (W \times \bigsqcup_{Q \supset R \supset P}A(R)_{^*P,\sigma} \times V_{^*P} ) \times V_Q. \ee
	This is open in $\bord{X_{\mq}} \times \widehat{U}_Q$ equipped with the product topology as $W \times \bigsqcup_{Q \supset R \supset P}A_{^*P(R),\sigma} \times V_{^*P}$ is open in $\bord{X_{\mq}}$ by an analogue of Lemma \ref{sp nbhd}. 
	
\end{proof}

\spoint We refer to the Appendix \ref{appendix:abs-retracts} for a brief review of the notion of Absolute Retract (AR) and Absolute Neighborhood Retract (ANR). 

\begin{nprop} \label{eq ar}
	For $Q \in \ratp'$, we have that $\overline{\e(Q)}$ is an Absolute Retract (AR).
\end{nprop}

\begin{proof}

	As $\bord{X_{\mq}}$ is contractible (see \cite[Lemma 8.6.4]{borel-serre}) and an ANR  (\cite[Lemma 8.3.1]{borel-serre}), it is an AR using Proposition \ref{prop:ANR-AR}. Likewise $\widehat{U}_Q$, which is homeomorphic to $\R^{\mathbb{N}}$, is an AR using Proposition \ref{prop:ANR}(1) and (2). Since a product of two AR spaces is again AR (see Proposition \ref{prop:ANR}(2)), we are done. \end{proof}

\tpoint{The affine, rational Tits building} Recall the affine Weyl group and its Coxeter structure defined in \S \ref{subsub:Weyl-group} $(\affW, S)$ as well the parabolic subgroups $\affW(J)$ for $J \subset I$ introduced in \textit{loc. cit.}. Let $\Sigma(\affW, S)$ be the poset consisting of cosets of the form $w \affW(J)$ where $w \in \affW$ and $\affW(J):= \la s_i \mid i \in J \ra$ for $J \subset I$, ordered by reverse inclusion, \textit{i.e.} $B \leq A$ in $\Sigma(\affW, S)$ if and only if $B \supset A$ as subsets of $\affW$. In this case we say $B$ is a face of $A$. 
It is known (see \cite[Thm. 3.5, p. 116]{brown:buildings}) that $\Sigma(\affW, S)$ is a simplicial complex with vertex set identified with $\affW$. Regarded as a simplicial complex, we call $\Sigma(\affW, S)$ the \textit{affine Coxeter complex}. 

An \textit{affine building} is a simplicial complex $\mc{T}$ which can be expressed as a union of subcomplexes $\Sigma$ called \textit{apartments}, and which satisfy the following conditions. \begin{enumerate}
	\item Each apartment $\Sigma$ is (equivalent, as a simiplical complex) to the \textit{affine} Coxeter complex. 
	\item For any two simplices $A, B \in \mc{T}$, there is an apartment containing both of them.
	\item If $\Sigma, \Sigma'$ are two apartments containing simplices $A$ and $B$, there is an isomorphism $\Sigma \rr \Sigma'$ fixing $A$ and $B$ pointwise.
\end{enumerate}

The main example of an affine building for us will be the (positive) \textit{ rational Tits building $\build(\hG)$} whose vertices are given by maximal elements,\textit{ i.e.} the set of maximal parabolics $\maxratp \subset \ratp.$ The $k$-simplices are then formed by intersections $P_{1} \cap \cdots \cap P_{k}$ for $P_{i} \in \maxratp$ which also lie in $\ratp$. We equip $\build(\hG)$ with the weak topology (cf.  \cite[p. 41]{Lundell2012-rs}).  Since $\hG_{\Q}/ P_{J, \Q}$ is in bijective correspondence with the rational parabolics conjugate to the standard parabolic $P_J$, the set of all  $k$-simplices are in bijection with $\cup_{J \subset I, |J| = k+1} \, \hG_{\Q}/ P_{J, \Q}$ and hence \be{} \build(\hG) = \cup_{\es \neq J\subset I} \; \hG_{\Q}/ P_{J, \Q} \ee

\tpoint{A criterion for homotopy equivalence} \label{subsub:homotopy-equiv-criterion} Let $V$ be a topological space along with a locally finite covering $\{V_j\}_{j \in J}$ by non-empty closed subsets $V_j \subset V$. The $\textit{nerve}$ $T$ of this cover is the simplicial complex whose vertices are $J$, and whose simplexes are finite subsets of $S \subset J$ such that $V_S = \cap_{j \in S} \, V_j$ is non-empty. Denote by $\Ss$ the set of simplices in $T$, and let $\underline{T}$ be the geometric realization of $T$ (see \cite[Section 8.2]{borel-serre}).

\begin{nprop} \label{prop:homotopy eq gen res}  \cite[Theorem 8.2.1]{borel-serre} In the notation above, assume the following,
	
	\begin{itemize}
		\item There exists a positive integer $M$, such that $\text{card}(S) < M$ for all $S \in \Ss$.
		
		\item for each $S \in \Ss$, the space $V_S$ is an Absolute Retract (AR).
	\end{itemize}   
	Then the spaces $V$ and $\underline{T}$ are of the same homotopy type.
\end{nprop}

\tpoint{Homotopy type of $\partial  \, \bord{\xg} $} \label{subsub:homotopy-type-boundary}  Our aim will now be to show the following.

\begin{nthm} \label{thm:boundary-homotopy}
	The spaces $\partial  \, \bord{\xg}$  and $\build(\hg)$ are homotopy equivalent.
\end{nthm}
\begin{proof}  We apply the previous proposition to the case when $V=\partial  \, \bord{\xg}$  and the covering is $\{ \overline{\e(Q)} \}_{Q \in \maxratp}$. This cover is locally finite since for any given point $z \in \partial  \, \bord{\xg}$ , there is a unique $P \in \ratp$ such that $z \in \e(P) \subset \partial  \, \bord{\xg}.$ On other other hand, if $\e(P) \cap \overline{\e(Q)}$, then from Proposition \ref{prop:closure-e}(3),  $Q \supset P$. However, for a fixed $P$, there are only finitely many such $Q$. As a simplicial complex, the nerve of this covering   is $\build(\hg)$ since for a finite subset $S \subset \maxratp$,
\be{} \bigcap_{Q \in S} \overline{\e(Q)}  = \overline{\e(\bigcap_{Q \in S} Q)}  \mbox{ if and only if } \bigcap_{Q \in S} Q \in \ratp. \ee Hence for any such set $S$ we must have $| S | \leq \ell+2$, verifying the first assumption of Proposition \ref{prop:homotopy eq gen res}. The second assumption holds by applying Proposition \ref{eq ar}.

\end{proof}

\subsection{$\bord{\xg}$ and partitions}

In this subsection, we would like to extend the semi-stability partition of $\xg$ to one for $\bord{\xg}.$

     \newcommand{\Xf}{\mathfrak{B}}

     \tpoint{On the sets $\pbord{S}$}   \label{subsub:bord-partition} For a rational parabolic $S \in \stdp$, we shall write \be{} \begin{array}{lcccr} \label{S:Q-notation}  S \sim (Q, \delta)  & \text{ with } & Q \in \stdp, \,  \delta \in \hG_{\Q}. & \mbox{ to mean} & S = Q^{\delta}  \end{array} \ee  Of course $\delta$ is only well-defined modulo  $Q_{\Q}$. For the remainder of this section, we shall reserve the letters $P, Q, R$ for a standard parabolics, so $P \sim (P, e)$, etc.  For $P \subsetneq \hG$ a standard parabolic, define
     \be{} \label{bord:Pe-desc} \pbord{P}&:=& X^{ss}_{M_P} \times \cor(\hAp^+_{P,(0, 1)}) \times \hU_P \\
     &=&   \xg(P,e) \ \ \cup \ \  \label{bord:Pe-desc-2}\bigsqcup_{Q \supset P} X_{M_Q}(^*P,e) \times \widehat{U}_Q, \ee where $\cor(\hAp^+_{P,(0, 1)}):= \hAp^+_{P, (0, 1)} \cup \ \ \sqcup_{Q\supset P} \, \hAp(Q)_{\sP, (0, 1)}$ is the interior of the closure of $\hA^+_{P, (0, 1)}$ in the corner $\corn{P}$. The second equality \eqref{bord:Pe-desc-2} follows from  \eqref{XGP:MPss} and implies  \be{} \label{bord cap eq}\pbord{P} \cap \e(Q) = \begin{cases} \es & \mbox{ if } Q \not\supset P \\ X_{M_Q}(^*P,e) \times \widehat{U}_Q &  \mbox{ for } Q \supset P. \end{cases} \ee 
     
    \noindent  Now for $S \sim (P, \delta)$ a general rational parabolic, in the absence of a simple relation like \eqref{XGP:MPss}, we start from the expression $R_{\delta^{-1}} \,\left(   \xg(P, e)\right) = \xg(P, \delta) $ and set \be{} \label{def:bord-S} \begin{array}{lcr} \pbord{S}:=  R_{\delta^{-1}} \,  \left( \pbord{P} \right) & \mbox{ for } & S \sim (P, \delta) \end{array}. \ee To check this is well-defined, we need the following result.

  \begin{nlem} \label{Par trans for bord} Let $P \in \stdp'$  and $\delta \in \hG_{\Q}.$ Then  $\pbord{P} \cap R_{\delta} \, \left( \pbord{P} \right)  \neq \emptyset$ if and only if $\delta \in P_{\Q}$, and moreover in this case $R_{\delta} \,  \left( \pbord{P} \right) = \pbord{P}$
  \end{nlem}   
    \begin{proof} Using Remark \ref{lem:can-pair-transf} it suffices to prove the lemma `along the boundary': given $x \in \pbord{P}, \delta \in \hG_{\Q}$ such that $x\delta \in \pbord{P}$ as well, and assuming there exist standard parabolics $Q,R$ containing $P$ and such that $x \in \pbord{P} \cap \e(Q)$ and $x\delta \in \pbord{P} \cap \e(R)$, we want to show that $\delta \in P_{\Q}$. To see this, first note that $\e(Q)\cdot \, \delta \cap \e(R) \neq \emptyset$ which is true if and only if $\e(Q^{\delta}) \cap \e(R) \neq \emptyset$, i.e. $Q^{\delta} =R.$ Since $Q,R$ are standard parabolics we must have $Q = R$ and
       $\delta \in Q_{\Q}$.  From \eqref{bord cap eq}, 
       \be{} x,x\delta \in  \pbord{P} \cap \e(Q) = X_{M_Q}(^*P,e) \times \widehat{U}_Q\ee
       Denoting by $\pr_{M_Q}: \e(Q) \longrightarrow X_{M_Q}$ the natural projection map and writing $\delta  = \delta_m \, \delta'$  
       with $\delta_m \in M_{Q, \Q}$ and $\delta' \in \hA_{Q} \, \hU_{Q} \subset \hP_{\Q}$, we have by the formula for the action of $\hG_{\Q}$ on $\e(Q)$ (see \eqref{Rg:act-boundary})
       \be{} \pr_{M_Q}(x),\pr_{M_Q}(x)\delta_m \in X_{M_Q}(^*P,e) \ee
       which, by a finite dimensional version of Remark \ref{lem:can-pair-transf}, implies that $\delta_m \in \sP_{\Q}.$ Hence $\delta \in P_{\Q}$. 
   \end{proof}

\tpoint{A partition of $\e(Q)$} Generalizing \eqref{bord cap eq}, we have the following result.

\newcommand{\odelta}{\overline{\delta}}

\begin{nprop} \label{prop:e(Q)-part} Let $Q, S \in \ratp'.$ 
	\begin{enumerate}
		\item Suppose $S\sim (P, \delta)$ with $P \in \stdp$ and $\delta \in \hG_{\Q}/P_{\Q}.$ Then $\pbord{S} \cap \e(Q)$ is non-empty if and only if $Q \supset S$ in which case $\delta \in Q_{\Q}$ and 
		\be{} \pbord{S} \cap \e(Q) =  X_{M_Q} (^*P, \odelta) \times \hU_Q \ee for $\odelta$ the element in $M_{Q, \Q}/^*P_{\Q}$ corresponding to $\delta$.  
		\item We have a partition 
		\be{} \label{bordx part e(Q)} \e(Q) = \bigsqcup_{S \in \ratp(Q)} \e(Q) \cap \pbord{S} ,\ee where we recall that $ \ratp(Q)= \{ S \in \ratp, S \subset Q \}$.   \end{enumerate}
\end{nprop}

\begin{proof}
	As for (1), first observe that by conjugation we may assume $Q \in \stdp$ as well, which we do from now on. Assume the intersection $\pbord{S} \cap \e(Q) \neq \es$. Then as $S=P^{\delta^{-1}}$, we may apply $R_{\delta}$ to  obtain $\pbord{P} \cap \e(Q^{\delta}) \neq \es$.  Applying \eqref{bord cap eq}, we find that $Q^{\delta} \in \stdp$ and  $Q \supset S,$ so that $\delta \in Q_{\Q}$ as well. The result follows from \eqref{bord cap eq}.

	As  for (2), note that as $\delta$ ranges over $Q_{\Q}/P_{\Q}$, the subgroups $S=P^{\delta^{-1}}$ ranges over all parabolics in $\ratp(Q).$ So we may write, using part (1) and the canonical pair partition of $X_{M_Q}$, 
	\be{}  \bigsqcup_{S \in \ratp(Q)} \, \e(Q) \cap \pbord{S}\,  = \bigsqcup_{\odelta \in M_{Q,\Q}/\sta{P}_{\Q} } X_{M_Q}(^*P,\odelta) \times \widehat{U}_Q = X_{M_Q} \times \widehat{U}_Q = \e(Q) \ee \end{proof}

 \tpoint{A partition of $\bord{\xg}$} Abusing notation slightly, let us set $\pbord{\hG}:= \xgss. $
 \begin{nprop} \label{prop:bord:partition}
      There is a partition of $\bord{\xg}$,     \be{} \label{bord:partition} \bordX =\bigsqcup_{S \in \ratp} \pbord{S} \ee
 \end{nprop}
\begin{proof} The fact that the right hand side covers $\bord{\xg}$ follows from \eqref{prop:e(Q)-part}.
	     If $S \in \ratp$ is proper, then $\pbord{\hg}=\xg^{ss}$  and $\pbord{S}$ are disjoint by Proposition \ref{prop:basic-partition}, so we may suppose, by way of contradiction, that $S,T \in \ratp$ are two proper parabolics such that $\pbord{S} \cap \pbord{T} \neq \es$. Moreover, without loss of generality we may assume $S=P$ and  $T = Q^{\zeta^{-1}}$ with $P,Q \in \stdp$ and $\zeta \in \hG_{\Q}/Q_{\Q}$. 
     Pick an $x \in \pbord{P} \cap \pbord{Q^{\zeta^{-1}}}.$  By the uniqueness in the theorem of canonical pairs, the intersection must lie along $\partial \, \bord{\xg}$, \textit{i.e.} there exists $R \in \ratp$ such that \be{} \label{P,Q,R}x \in \pbord{P} \cap \pbord{Q^{\zeta^{-1}}} \cap \e(R).\ee  
     Furthermore, by \eqref{bord cap eq}, we must have $R \supset P$ and thus $R \in \stdp$ as well. By the same argument, we also have  $Q^{\zeta^{-1}}  \subset R$ which implies that $R^{\zeta} \in \stdp$ and hence $\zeta \in R_{\Q}/Q_{\Q}$. Letting $\pr_{M_R} :\e(R) \longrightarrow X_{M_R}$ be the natural projection map, from \ref{P,Q,R} and \ref{bord cap eq} we can then conclude that
     \be{} \label{bd pqr}\pr_{M_R}(x) \in X_{M_R}(^*P,e) \cap X_{M_R}(^*Q,\ozeta) \ee
where $\ozeta \in M_{R, \Q}/^*Q_{\Q}$ is the unique element corresponding to $\zeta$. This forces $^*P =\leftidx{^*}Q,$ and $\ozeta \in ^*Q_{\Q}$, which implies $P=Q$ and $\zeta \in Q_{\Q}$. Hence, $S  = T$. \end{proof}

 \spoint Collecting parabolics by their type, we now define \be{} \begin{array}{lcr} \bord{P} := \bigsqcup_{S \sim (P,\delta)} \, \pbord{S} & \mbox{ and } & \bord{\hG} = \pbord{\hG} = \xg^{\ss} \end{array}, \ee where $P \in \stdp'$ with the union over all $\delta \in \hG_{\Q}/P_{\Q}.$ 

 \begin{nprop} \label{par trans bord} \begin{enumerate}
 		\item $\bord{\xg} = \sqcup_{P \in \stdp} \, \bord{P}$
 		 
 		\item For $P \in \stdp'$, the inclusion $i: \pbord{P} \hookrightarrow \bord{P}$ induces a bijection \be{} \label{bijection:can-pairs-parabolic-trans bord}   \pbord{P} / P_{\Q} \stackrel{1:1}{\longrightarrow}  \bord{P}/\hG_{\Q}. \ee \end{enumerate}     
 \end{nprop}
 \begin{proof} The first result is a restatement of Proposition \ref{bord:partition}. As for the second, 
 	if $x \in \bord{P}$, there exists a unique $\delta \in \hG_{\Q}/P_{\Q}$ such that $x \in \pbord{S}$ with $S \sim (P,\delta).$ Hence $x\delta \in \pbord{P}$ and $i(x\delta) = x$ in $\bord{\qxg(P)}/\hG_{\Q}$. Injectivity follows from \ref{Par trans for bord}.
 \end{proof}

  \newcommand{\nbhd}{\mathcal{N}}

\section{On the structure of $\bord{\xg} / \widehat{\Gamma}$ } \label{sec:bord-mod-Gamma}

\subsection{On the parabolic ends and Siegel sets for $\bord{\xg}$}

\spoint  From Propositions \ref{prop:bord:partition} and \ref{par trans bord},we obtain a disjoint union
\be{} \label{bord:quotient-Gamma} \bord{\xg} / \hGam = \xgss / \hGam \, \cup \, \bigsqcup_{P \in \stdp'} \, \pbord{P} / \hGam_P. \ee For $P \in \stdp'$ we have  from \eqref{bord:Pe-desc}  $\pbord{P}= X_{M_P}^{ss} \times \corn{\hAres{P}} \times \hU_P$ and hence a projection  \be{} \mathfrak{f}_P: \pbord{P} \rr X_{M_P}^{ss} \times \corn{\hAres{P}}, \, \, x \mapsto \left(\pr_{M_P}(x), \pr_{\hA^+_P}(x)\right),  \ee extending the map $f_P$ from \eqref{def:f_P}. Moreover we also obtain 

\begin{nprop} \label{prop:parabolic-ends-revisted} For each $P \in \stdp',$ the projection $\mf{f}_P$ above descends to a continuous surjection that is denoted by the same name and which fits into the diagram where $i$ denotes natural inclusion   \be{} \label{comm-diag-fp} \begin{tikzcd} 
		\xg(P,e) / \hgam_P \ar[r, "i"] \ar[d, "f_P"] &  \pbord{P} / \hgam_P \ar[d, "\mf{f}_P" ] \\ 
		X_{M_P}^{ss} / \hgam_{M_P} \times \hAp^+_{P, (0, 1)} \ar[r, "i"]  & X_{M_P}^{ss} / \hgam_{M_P} \times \corn{\hAp^+_{P, (0, 1)}} \end{tikzcd}. \ee The fibers of both $f_P$ and $\mf{f}_P$ are homeomorphic to $\hU_P/ \hgam_{U_P}.$
\end{nprop}

\begin{nrem} \label{remark:harder-approach} We can obtain a similar diagram for $P \in \stdp'$ replaced with a general $S \in \ratp'$, and more generally (using the same notation as introduced earlier) \be{} \label{comm-diag-fp} \begin{tikzcd} 
		\qX_G(P) / \hgam \ar[r, "i"] \ar[d, "f_P"] &  \bord{P} / \hgam \ar[d, "\mf{f}_P" ] \\ 
		X_{M_P}^{ss} / \hgam_{M_P} \times \hAp^+_{P, (0, 1)} \ar[r, "i"]  & X_{M_P}^{ss} / \hgam_{M_P} \times \corn{\hAp^+_{P, (0, 1)}} \end{tikzcd}. \ee

	As remarked in the introduction and taking into account Proposition \ref{par trans bord}, we obtain an alternate approach to construct $\bord{\xg}/\hGam$:  start with the `core' $\xg^{ss}/\hGam$ and then add in the boundary points along the $\hA^+_P$ directions to the parabolic ends $\qX_G(P)/\hgam.$ \end{nrem}

\tpoint{Analogues of Siegel sets for $\bord{\xg}/\hGam$}

\newcommand{\corS}{\mathscr{S}}

Recall the Siegel sets $\sie_{t, \Omega}$ constructed in \S \ref{subsub:siegel-sets} as well the neighborhood basis $\Nb_P(W, t, V)$ from \eqref{nhbd:W-sigma-V}.  For $t >0$ and $\Omega \subset \R$ we then define 
	\be{} \label{siegel W}  \corS_{t,\Omega} = \Nb_B(\es, t, \hU_{\Omega}) = \corn{\Ac^+_t} \times \hU_{\Omega}. \ee 
Let us note that $\corS_{t, \Omega} \cap \xg = \sie_{t, \Omega}$. More generally,

\begin{nprop} \label{prop:bord-Siegel-corner} Let $t > 0$ and $\Omega \subset \R$ a bounded subset of $\R$. For each $P \in \stdp'$, there exists a bounded set $\Omega' \supset \Omega$ so that \be{} \label{bord-Siegel-corner} \corS_{t, \Omega} \cap \e(P) \subset  \mf{S}(M_P)_{t, \Omega'} \times \hU_{P, \Omega'}, \ee where $\mf{S}(M_P)_{t, \Omega}$ is a Siegel set for $M_P$, defined as in \eqref{siegel-set} but with respect to $A_{\sta{B}}$ and $U_{\sta{B}}$ where $\sta{B} \subset M_P$ is the Borel corresponding to $B \subset P$. \end{nprop}

For any $P \in \stdp'$, we have a decomposition $\hU = \hU_{\sta{B}} \rtimes \hU_P.$ We would like to understand the image of $\hU_{\Omega}$ under the natural projections associated to the semi-direct product. The main thing to notice is that the description we gave for $\hU_{\Omega}$ depends on a specifc choice ordering of the roots of $\rts_o$. This may not be compatible with the semi-direct prodct decomposition above, but after a finite set of commutations, which will change  the set $\Omega$ to a set $\Omega'$, it can be made so.

	\spoint Generalizing Theorem \ref{thm:reduction-theory-1}(1), we now have 

\begin{nprop} \label{prop:bord-siegel} Let $t_0, \Omega_0$ as in Theorem \ref{thm:reduction-theory-1} and pick  $\Omega \supset \Omega_0 $ a bounded subset and $t < t_0.$ There exists a bounded set $\Omega' \supset \Omega$ such that for any $x \in \bordX$, we can find $\gamma \in \hGam$ such that $x\gamma \in \corS_{t,\Omega'}.$
\end{nprop}
\begin{proof}
	For any $x \in \bordX$, there exists $\gamma \in \hGam$ such that $x\gamma \in \cor{(\hB)}$, the corner attached to $\hb$, see \eqref{corner:dec}. As we will work within this corner, without loss of generality we assume $x \in \corn{\hB}$. Using Theorem \ref{thm:reduction-theory-1} , there exists $\gamma \in \hGam$ so that, if $x \in \xg \cap \corn{\hB}$, then $x\gamma \in \sie_{t, \Omega}=  \corS_{t, \Omega} \cap \xg$.

	 Next we treat the case  $x \in \e(P)$, say $x=(z, u)$ with respect to the decomposition $\e(P)= X_{M_{P}} \times \widehat{U}_{P}.$ If  
	 $\gamma \in \hGam \cap P$ is decomposed as $\gamma = r_{\gamma}u_{\gamma}$ with $r_{\gamma} \in \Gamma_{M_P}$, $u_{\gamma} \in \hGam \cap \widehat{U}_P$ (\textit{cf.} Lemma \ref{lem:GammaL} ), we have   
	\be{} \label{xgam:P}  x\gamma = (z r_{\gamma}, u^{r_{\gamma}}u_{\gamma}) \ee using  \eqref{Rg:act-boundary}. As $\Omega$ was assumed bounded, there exists (bounded) $\Omega' \supset \Omega$ as in Proposition \eqref{prop:bord-Siegel-corner}. By reduction theory for $M_P$, we can find $r_{\gamma} \in \hGam \cap M_P$ so that $ z r_{\gamma} \in \hA_{\sta{B}, t} \times \hU_{\sta{B}, \Omega'},$ i.e. $z r_{\gamma} \in \mf{S}(M_P)_{t, \Omega'}$. Having fixed $r_{\gamma}$ we can then choose $u_{\gamma} \in \hU_P\cap \hGam$ so that $u^{r_{\gamma}}u_{\gamma} \in \hU_{P, \Omega'}$ by applying the reduction theory for $\hU_P$. We are then done by applying Proposition \eqref{prop:bord-Siegel-corner} once again.
	To finish the proof, we may need to extend $\Omega$ for each standard parabolic $P$, but there are only finitely many such extensions required, each one producing a new bounded set. 	\end{proof}

\subsection{Some compactness results}

\tpoint{Compactness along the boundary} As for the boundary of $\bord{\xg}$, we have:
\begin{nprop} \label{prop:compact-boundary} The space $\partial \, \bord{\xg} /\hGam$ is compact. \end{nprop}

\begin{proof} From Proposition \ref{prop:closure-e}, we deduce that $\partial \, \bord{\xg}$ has a closed cover $\{\overline{\e(Q)}\}$ for $Q$ ranging over maximal parabolics. Letting $\stdmaxp$ be the subset of standard, maximal parabolics, 
	\be{} \label{partbdqt} \partial \bord{\xg} /\hGam =  \bigcup_{Q \in \stdmaxp}\overline{\e(Q)}  / \hGam_Q. \ee
	
	Indeed, if $x \in \overline{\e(Q)} \cap \overline{\e(Q)}\gamma$ for some $\gamma \in \hGam$, then from Proposition \ref{prop:closure-e} there exist rational parabolics $ S$ and $T$ which are contained in $Q$ such that $x \in \e(S) \cap \e(T)\gamma$. Hence $\e(S) = \e(T^{\gamma})$, and so also $S  = T^{\gamma}.$ This in turn implies $Q \cap Q^{\gamma} \neq \emptyset$ and finally $\gamma \in \hGam_Q$.

	Using \eqref{partbdqt}, it is enough to show that $\overline{\e(Q)}/\hGam_Q$ is compact for a standard, proper parabolic $Q \in \ratp$. To see, this, recall from (\ref{e(Q) bar des}) the homeomorphism,
	\be{} \label{closure-e-prod} \overline{\e(Q)}\simeq \bord{X_{M_Q}} \times \widehat{U}_Q \ee
	The compactness of $\overline{\e(Q)}/\hGam_Q$ then follows from the compactness of $\bord{X_{M_Q}}/\hGam_{M_Q}$ (cf. \cite[ Theorem 9.3]{borel-serre}) and the fact that $\widehat{U}_Q/\hGam_{\widehat{U}_Q}$ is compact by arguing as follows:  from \eqref{closure-e-prod}, we obtain a projection $\pi: \overline{\e(Q)}/\hGam_Q \longrightarrow \bord{X_{M_Q}}/\hGam_{M_Q}$ which makes $\overline{\e(Q)}/\hGam_Q$ the total space of a fiber bundle with base $\bord{X_{M_Q}}/\hGam_{M_Q}$ and fibers homeomorphic to $\widehat{U}_Q/\hGam_{\widehat{U}_Q}$. Since the fibers are compact, it follows that $\pi$ is locally closed, which can be verified by picking a trivializing neighborhood around any point in the base $\bord{X_{M_Q}}/\hGam_{M_Q}$. Hence $\pi$ is a closed map with compact fibers and so a proper map (\cite[Proposition 4.93(c)]{Lee2010-gn}), from which the desired result follows.

\end{proof}

\tpoint{Compactness in the interior} \label{subsub:compact-interior} The space $\bord{\xg}/\hGam$ cannot be expected to be compact, since $\bord{\hA^+_{P, (0, 1)}} $ is itself non-compact. In analogy with Lemma \ref{lemma:compact bordh}, we can state the following.

\begin{nprop} \label{prop:compact-interior} Pick $t \geq t_0$ and $\Omega$ be a closed, bounded subset so that Proposition \ref{prop:bord-siegel} holds.  Let $Z$ be a $\hGam-$invariant closed subset of $\bord{\xg}$ and suppose there exists $M_0 , r_0 >0$ such that for all $ x \in Z \cap \corS_{t,\Omega} \cap \xg$, we have $H_B(x) \in \mc{K}(\hfh^+; M_0, r_0, \log t)$, the set defined in \eqref{C:m-t-head}. Then the image of $Z$ in the quotient $\bord{\xg}/\hGam$ is compact.  
    \end{nprop}

  \begin{proof}   	By our assumption, if $s_0:= e^{-r_0}$ we have $Z \cap \corS_{t, \Omega}   \subset  \corn{\hA^+([-M_0, M_0], s_0, t)} \times \hu_{\Omega}$
   where $\corn{\hA^+([-M_0, M_0], s_0, t)}$ is as in \eqref{compact-corner-A}. As $\hU_{\Omega}$ is also compact, 
	we have $\corn{\hA^+([-M_0, M_0], s_0, t)} \times \hu_{\Omega}$  is compact in $\bord{\xg}$, and hence so is its image in the quotient by $\hGam$. But $Z \cap \corS_{t, \Omega}$ is a closed subset which projects surjectively onto this image of $Z$ by $\hGam,$ so the Proposition follows. \end{proof}

\subsection{On the separability of $\bordX / \hGam$} \label{sub:separability}

\tpoint{Properness and Separability} \label{subsub:prop-separability}

Let $Z$ be a Hausdorff topological space. A continuous action of a group $H$ on $Z$ is said to be \textit{proper} if for any two points $x,y \in Z$ there exists open neighbourhoods $U_x \ni x$, $U_y \ni y$ in $Z$ such that the following set is finite, 
 \be {} \label{proper nbd}   \{\gamma \in H \mid  U_x \cap  U_y \gamma\neq \emptyset \}.   \ee In this case, $Z/ H$ is again Hausdorff in the quotient topology (cf. \cite[Ch. III, \S 4.2, Proposition 3, \S 4.4, Proposition 7]{Bourbaki1998-qu}).   We can then state the following criterion
 
 \begin{nlem} \label{proper g action} The $\hGam$ action on $\bord{\xg}$ will be proper, and hence $\bord{\xg}/\hGam$ Hausdorff, if the following condition holds: for $t$ and $\Omega$ chosen as in Proposition \ref{prop:bord-Siegel-corner}, 
 	\be{}   \label{SF-assumption} \begin{array}{lcr} \mbox{ \textbf{Siegel Finiteness (SF)} } & & | \{\gamma \in \hGam \mid \mathscr{S}_{t,\Omega} \cap \mathscr{S}_{t,\Omega} \gamma \neq \emptyset \} | < \infty \end{array} \ee 
 \end{nlem}
 
 \begin{proof}
 	Fix $\Omega$ bounded but open, so that $\mathscr{S}_{t,\Omega} \subset \bord{\xg}$ is open. For a given pair $x,y \in \bord{\xg}$, by Proposition \ref{prop:bord-Siegel-corner}, there exists $\gamma_x,\gamma_y \in \hGam$ such that $x\gamma_x,y\gamma_y \in \mathscr{S}_{t,\Omega}$. Choosing the open neighborhoods of $x$ and $y$ in $\bord{\xg}$, say  $ U_x = \mathscr{S}_{t,\Omega} \gamma_x^{-1}$ and $U_y=\mathscr{S}_{t,\Omega}\gamma_y^{-1}$ respectively, the finiteness of the set 
 	$\{\gamma \in \hGam| U_x \cap U_y \gamma  \neq \emptyset\}$ follows from assumption \textbf{(SF)}.
 \end{proof}

\tpoint{On Siegel Finiteness} \label{subsub:siegel-finiteness} For a fixed $0 < s < 1$, let use define   \be{} \corS^s_{t, \Omega}:=  \left(  
\sie_{t,\Omega} \cap X_{\hG^s} \right) \sqcup \bigsqcup_{J \subsetneq I} \corS_{t, \Omega} \cap \e(P_J). \ee

\begin {nprop} \label{siegel finite} Assume that we have the following
\begin{center} 
\textbf{(SF-U)} For any bounded subset $\Omega \subset \R$, $| \{ \gamma \in \hGam_{U} \mid \hU_{\Omega} \cdot \gamma \cap \hU_{\Omega} \neq \es \} | < \infty$  \end{center}

\noindent Then there exists $0< s_0<1$ so that if $0 < s < s_0 <1$, if  $t, \Omega$ are as in Proposition \ref{prop:bord-siegel}, 
    \be{}    | \{\gamma \in \hGam \mid \corS^s_{t,\Omega} \cap \corS^s_{t,\Omega} \gamma \neq \emptyset \} | < \infty.  \ee

\end{nprop}
\begin{nrem} The validity of condition \textbf{(SF-U)} in finite-dimensions is straightforward and we will plan to investigate its affine analogue elsewhere \end{nrem}

\begin{proof}
    If $\mathscr{S}^s_{t,\Omega} \cap \mathscr{S}^s_{t,\Omega} \gamma \neq \emptyset$ at least one of the following must occur: \textit{(a)} $\widehat{\mathfrak{S}}^s_{t,\Omega} \cap \widehat{\mathfrak{S}}^s_{t,\Omega} \gamma \neq \emptyset$; or \textit{(b)} 
$\left(\mathfrak{S}(M_{P_J})_{t,\Omega'} \times \hU_{\Omega'} \right) \cap \left(\mathfrak{S}(M_{P_J})_{t,\Omega'} \times \hU_{\Omega'} \right) \gamma \neq \emptyset$ for $\Omega'$ as in Proposition \ref{prop:bord-Siegel-corner}. If \textit{(a)} holds, choosing $s_0$ as in Theorem \ref{thm:reduction-theory-1}(2)-(3), we conclude that $\gamma \in \hGam_P$ for some (maximal) parabolic $P$. 
In \textit{(b)} holds, since $P_J$ and $P_K$ are not $\hG_{\Q}$-conjugate unless $J  = K$ and since  $\e(P)\gamma \subset \e(P^{\gamma})$ for $\gamma \in \hGam,$ we must have $\gamma \in \hGam \cap P_J.$ Hence we are reduced, after choosing $s_0$ sufficiently small, to proving the finiteness under the assumption that $\gamma \in \hGam_P$ for some $P \in \stdp'$. 

Suppose the intersection falls under case \textit{(b)} for $P=P_J$. Write $\gamma = r_{\gamma} u_{\gamma}$ with $r_{\gamma} \in \hGam_{M_P}$ and $u_{\gamma} \in \hGam_{U_P}$. Arguing as in the proof of Proposition \ref{prop:bord-siegel},  and using the Siegel finiteness result of Borel and Harish--Chandra for finite-dimensional groups \cite[Thm 15.2]{borel:red}, the possibilities for $r_{\gamma}$ are finite. The assumption \textbf{(SF-U)} would then tell us that the possibilites for $u_{\gamma}$ are also finite. A similar (but easier argument) allows us to deal with the situation in which case \textit{(a)} holds. \end{proof}

\begin{appendices}
	
\section{Topological Notions}

\subsection{Moore-Smith Convergence} \label{sub:moore-smith} 

The idea of constructing a topology on a set by specifying a class of convergent sequences is due to Moore--Smith (see \cite[Ch.2]{kelley:top}).

\spoint \label{subsub:conv-seq} Let $X$ be a set and $\mathcal{C}$ be a subset of $X^{\mathbb{N}} \times X$ consisting of pairs $(\{x_n\}_{n=1}^{\infty},x_0)$ with a sequence $\{x_n\}_{n=1}^{\infty}$ in $X$ and a point $x_0 \in X$. If a pair $(\{x_n\}_{n=1}^{\infty},x_0) \in \mathcal{C}$, we say the sequence $\{x_n\}_{n=1}$ $\mathcal{C}$-converges to $x_0$, and write this as
$x_n\xrightarrow{\mathcal{C}}x_0$ otherwise $x_n\nrightarrow x_0$.

\tpoint{Convergence class} \label{cclass}
A class $\mathcal{C}$ of pairs $(\{x_n\}_{n=1}^{\infty},x_0)$   is called a convergence class of sequences on the set $X$s if the following conditions are satisfied

\begin{enumerate}

	\item If $\{x_n\}_{n=1}^{\infty} $ is a constant sequence, \textit{i.e.}
	$x_n = y$ for $n \geq 1$ where $y \in X$, we have $x_n\xrightarrow{\mathcal{C}} y$
	
	\item If $x_n\xrightarrow{\mathcal{C}}x_0$, so does every subsequence of 
	$\{x_n\}_{n=1}^{\infty} $.
	
	\item If $x_n\nrightarrow x_0$, then there exists a subsequence 
	$\{x_{n_k}\}$ of $\{x_n\}$, such that for any further subsequence
	$\{x'_{n_k}\}$ of $\{x_{n_k}\}$; $x'_{n_k}\nrightarrow x_0$
	
	\item For a double sequence $\{x_{m,n}\}_{m,n=1}^{\infty} $, suppose for each fixed $m$, $x_{m,n} \xrightarrow{\mathcal{C}} x_{m,0}$ and the sequence $x_{m,0} \xrightarrow{\mathcal{C}} x_{0,0}$. Then there exists $f:\mathbb{N} \longrightarrow \mathbb{N}$ such that $\displaystyle \lim_{n \to \infty} f(n) = \infty$ and
	$x_{n,f(n)} \xrightarrow{\mathcal{C}} x_{0,0}$.
\end{enumerate}

\spoint If $\mathcal{C}$ is a convergence class of sequence on a space $X.$ We define for every subset $A$ of $X$, 
\begin{equation*}
	\bar{A} = \{x_0 \in X\ |\text{There exists a sequence} \{x_n\}_{n=1}^{\infty} \hspace{0.1cm} \text{in}\hspace{0.1cm} A \hspace{0.1cm} \text{such that} \hspace{0.1cm}  x_n\xrightarrow{\mathcal{C}}x_0 \}
\end{equation*}
The map $A \longrightarrow \bar{A}$ is a closure operator (see \cite[p.43]{kelley:top}) that defines a topology on $X$ whose closed sets are of the form $\bar{A}$ for $A \subset X$, see \textit{loc. cit.}

\begin{nlem} \label{lem:top-from-seq} Let $X$ be a set equipped with  a convergence class $\mc{C}$. Define the topology on $X$ coming from this convergence class as above.
	
	\begin{enumerate}
		\item \cite[Thm. 2.3]{kelley:top} $X$ is Hausdorff if and only if every convergent sequence has a unique limit.
		\item \cite[Chap. II, Prop. 1.8.13]{borel-ji} $X$ is compact if and only if every sequence has a convergent subsequence.
	\end{enumerate}
\end{nlem}

\subsection{Absolute Retracts}  \label{appendix:abs-retracts}

As a reference for the material in this section, we recommend \cite{hu:retracts}.
\tpoint{Absolute retracts}

A metric space $Z$ is called an \textit{Absolute Retract} (AR) if whenever $Z$ is embedded as a closed subspace in a metric space $Y$, there exists a retraction $r: Y \longrightarrow Z$, i.e. $r$ is continuous and its restriction to $Z$ is the identity.   This is equivalent (see \cite[Thm 3.1, 3.2]{hu:retracts}) to the condition that $Z$ is an \textit{Absolute Extensor} (AE):  if $A$ is any closed subspace of a metric space $Y$, and $f:A \longrightarrow Z$ is a continuous map, then there exists a continuous map	$F: Y \longrightarrow Z$ extending $f$. 


\tpoint{Absolute neighbourhood retracts}

A metric space $Z$ is called an \textbf{Absolute Neighbourhood Retract} (ANR) if whenever $Z$ is embedded as a closed subspace in a metric space $Y$, there is an open neighbourhood $U$ of $Z$ in $Y$ such that $Z$ is a retract of $U$. Again, this is equivalent  (\textit{loc. cit.}) to the condition that $Z$ is an \textit{Absolute Neighbourhood Extensor} (ANE): if $A$ is a closed subspace in a metric space $Y$, and $f:A \longrightarrow Z$ is a continuous map, then there exists a neighbourhood $U$ of $A$ in $Y$ and  a continuous map $F: U \longrightarrow Z$ extending $f$. An AR space is an ANR, and in fact:

\begin{nprop} \label{prop:ANR-AR} A metric space $X$ is an AR if and only if X is an ANR and contractible. \end{nprop}

\tpoint{Examples}

\begin{nprop} \label{prop:ANR}
	\begin{enumerate}
		\item \cite[Cor. 6.4]{hu:retracts} Any convex subset of a Banach space is an AR.
		\item \cite[Prop. 7.5]{hu:retracts} Any countable product of ARs is again an AR (with the product topology)
		\item  \cite[Prop. 7.5]{hu:retracts} Any finite product of ANRs is again an ANR (with the product topology)

	\end{enumerate}
\end{nprop}

\end{appendices}

	\newcommand{\F}{\ratp}
	
	\bibliographystyle{abbrv}
	\bibliography{LBS.bib}

\end{document}